\documentclass[11pt]{article}

\usepackage{sectsty}

\usepackage{amsmath,amsthm,bm,mathrsfs,color,comment}

\usepackage{amsmath,amssymb}
\usepackage{epsfig}
\usepackage{setspace}
\usepackage{rotating}
\usepackage{caption}
\usepackage{amsfonts}
\usepackage{mathtools}
\usepackage{float}
\usepackage{subcaption}

\usepackage{multirow}
\usepackage{mathrsfs}
\usepackage[percent]{overpic}

\usepackage[english]{babel}
\usepackage[utf8]{inputenc}
\usepackage{amsmath}
\usepackage{csquotes}

\usepackage[style=authoryear ,backend=biber]{biblatex}
\numberwithin{equation}{section}

\begin{document}

\title{Deterministic and Stochastic in-host Tuberculosis Models for Bacterium-directed and Host-directed Therapy Combination}
\author{ {\sc Wenjing Zhang}\\[2pt]
Department of Mathematics and Statistics,
Texas Tech University\\
Lubbock, Texas \ \ 79409-1042\\
wenjing.zhang@ttu.edu\\
}
\maketitle
\begin{abstract}
 {Mycobacterium tuberculosis infection can involve all immune system components and can result in different disease outcomes. The antibiotic TB drugs require strict adherence to prevent  both disease relapse and mutation of drug- and multidrug-resistant strains. 
 To overcome the constraints of pathogen-directed therapy, host-directed therapy has attracted more attention in recent years as an adjunct therapy to enhance host immunity to fight against this intractable pathogen.
 The goal of this paper is to investigate in-host tuberculosis models to provide insights into therapy development.
Focusing on therapy-targeting parameters, the parameter regions for different disease outcomes are identified from an established ODE model. 
Interestingly, the ODE model also demonstrates that the immune responses can both benefit and impede disease progression, depending on the number of bacteria engulfed and released by macrophages.
We then develop two It\^{o} SDE models,
which consider the impact of  demographic variations at the cellular level and  environmental variations during therapies along with demographic variations.
The SDE model with demographic variation suggests that stochastic fluctuations at the cellular level have significant influences on (1) the T-cell population in all parameter regions, (2) the bacterial population  when parameters located in the region with multiple disease outcomes, and (3) the uninfected macrophage population in the parameter region representing active disease. Further, considering environmental variations from therapies, the second SDE model suggests that disease progression can slow down if therapies (1) can have fast return rates and (2) can bring parameter values into the disease clearance regions.}

{In-host tuberculosis model, deterministic model, stochastic model, bifurcation analysis}
\end{abstract}
\section{Introduction}

Tuberculosis (TB) is a major cause of death from infectious diseases globally.
It infects one quarter of the world population and claims 1.4 million lives per year (\cite{world2020global}). 
Even though it is an ancient disease, TB is still formidable due to the spread of drug-resistant strains, which can lead to re-emergence in the regions of the world with smaller TB burdens. 
The world is at a critical tipping point in the elimination of TB (\cite{world2015end}). 
The primary strategies to fight against TB are rapid identification, appropriate treatment of infected cases, and continued surveillance to manage outbreaks.
It is necessary to incorporate interventions at the population level with a decrease in the disease progression at the individual level.

On an individual basis, the disease outcome varies among clearance, active disease, and latent tuberculosis infection (LTBI) (\cite{gammack2005understanding},\cite{meermeier2018early}). 
In more detail, after the exposure with \emph{Mycobacterium tuberculosis} (Mtb),  a fraction of individuals can eradicate the infecting Mtb pathogens naturally and achieve early clearance. Evidence of early clearance in Uganda (\cite{bark2017identification,}) ranges from $3.4\%$ to $23.8\%$ (\cite{cobat2009two}).
Another roughly $5\%$—$10\%$ of infected individuals rapidly progress to active disease (fast progression) (\cite{sud2006contribution}) after initial exposure, 
while the majority of infected individuals stay in latent tuberculosis infection status (\cite{lin2010understanding}).
The global prevalence of LTBI is estimated about $25\%$ (\cite{cohen2019global}).
Patients with LTBI status have Mtb pathogens present but show no clinical symptoms.
However, endogenous reactivation (slow progression) can still induce disease progression to active TB for the rest of their life at about $5–10\%$ (\cite{world2020global}).
Moreover, even though only about one in ten LTBI cases become active, this progression contributes to approximately $80\%$ of active disease cases globally (\cite{schwartz2020tuberculosis}).
One of the causes of different disease outcomes and disease progression is an impaired host immune system, particularly in the case of immune suppression by HIV and malnutrition. 
 For HIV infected individuals, the damage to the innate immunity  increases macrophage turnover (\cite{hasegawa2009level}), and the impairment of adaptive immune system reduces CD4+ T cell counts and the baseline monocyte to lymphocyte ratio (\cite{naranbhai2013blood}).
 For malnourished TB patients (such as vitamin D deficiency),
delayed recovery and higher mortality rates are reported  (\cite{gupta2009tuberculosis}).

The development of the antibiotic TB medication streptomycin started in 1943 (\cite{britannica2013}). Nowadays, LTBI and active TB disease can be treated. 
The anti-TB drug treatments include four first-line drugs,
isoniazid (INH), rifampin (RIF), pyrazinamide (PZA), and ethambutol (EMB), which either act as a bactericide or a bacteriostatic agent. 
Even though the anti-TB regimens have an up-to-$95\%$ efficacy (\cite{conde2011new}),  treatment effectiveness differs according to the adherence of strict routine conditions (\cite{mitnick2009tuberculosis}). If patients stop using medications too soon, take incorrect medication, or use irregular medication, the Mtb bacteria may still live, which can cause disease relapse and mutation of antibiotic-resistant  strains. The antibiotic-resistant Mtb strains are much more difficult to treat. 
Frustrated by antimicrobial resistance and strict treatment adherence, researchers have focused on host-cell factors to inhibit Mtb pathogens' replication or persistence and promote host immune responses against the invading pathogens.

Host-directed therapies (HDTs) are emerging as adjunctive therapies lately (\cite{zumla2020host}).  HDTs aim at the enhancement of long-term treatment outcome and the functional cure of persistent Mtb bacteria. Examples of adjunctive HDTs include (1) nonsteroidal anti-inflammatory drugs, such as ibuprofen, to prevent host inflammatory responses, (2) vitamins and dietary supplements, such as vitamin D, to promote immune response, (3) and other drugs to reduce intracellular bacillary load of Mtb (\cite{zumla2015towards}).
To facilitate the development of novel therapy strategies, identification of the key mechanisms of the interactions between Mtb bacteria and host immune response is crucial. 
This serves as the motivation of this project.

Mathematical modeling has been a fundamental research tool for studying complex dynamics and identifying the underlying causal factors. 
The papers (\cite{wigginton2001model}, \cite{marino2004human}, \cite{gammack2005understanding} and \cite{sud2006contribution}) are starts of in-host modeling regarding immunology of Mtb infections and host-immune response.
These previous works successfully combined both bacterial and immune response mechanisms and developed models of Mtb-host interactions to identify the key factors for various disease outcomes through deterministic models.
The study of how stochastic variations affecting in-host dynamics has still fallen short.
This may be due to the high model complexity.
The previous in-host modeling work contains a large number of equations and parameters, which put a big challenge on the identification of significant model parameters driving the disease progression.

In this project, we adopt an established in-host TB infection model (\cite{du2017simple}), which incorporates the essential bacterial and immune response mechanisms in a four-dimensional deterministic system with eighteen parameter values and  successfully generates all disease outcomes  including clearance, LTBI, and active disease.
In this contribution, however, Du et al. employed and analyzed only an asymptotic version of the model that neglects the effects of the CD4 T-cell population. 
Then, the papers (\cite{zhang2020analysis}, \cite{Zhang2020TB}, and \cite{zhang2021investigation}) present  rigorous mathematical analyses on a reduced model and the original full model, respectively. 
In this paper, we focus on the driving factors and consider the stochastic variations on demographic variables (including Mtb pathogen and immune cells populations) and on environmental variables (including the identified driving factors) to explore their influences on disease progression.

  The paper is organized as follows. In the rest of section 1, we present an established deterministic Mtb-host model and its basic mathematical properties. In section 2, we formulate two stochastic models considering demographic variations and environmental variations. In section 3, we delimit the parameter regions according to different disease outcomes through bifurcation analysis. We further demonstrate the stochastic influences on each parameter region. In section 4, we investigate the speed of the disease progression under combination therapies.
The paper ends with a conclusion and discussion in section 5.

\subsection{Mtb-host Dynamics}
Mtb infection  most frequently happens in
the respiratory tract, particularly the lung (pulmonary TB), and the regional lymph nodes.
When Mtb are inhaled into the lungs and taken up by
resident alveolar macrophages, these bacteria start to multiply and make alveolar macrophages their main target.
If receiving adequate stimulation for activation, 
macrophages can effectively ingest and destroy their phagocytized bacteria.
Otherwise, the phagocytized bacteria reproduce inside their host macrophages.
The host macrophages are eventually unable to be activated due to the increased phagocytized bacteria load, and become chronically infected. 
The chronically infected macrophages, which are unable to kill their intracellular bacteria, ultimately either undergo programmed cell death ( i.e., apoptosis or necrosis) due to an excessive intracellular bacterial load
or are destroyed through T-cell-mediated immune responses. 
Both processes lead to the death of chronically infected macrophages, which release intracellular bacteria to the extracellular environment. 
These extracellular bacteria again are engulfed by activated macrophages and result in either bacterial elimination or chronically infected macrophages. 
The initiation of the T-cell-mediated immunity starts from
the infected front-line innate immune cells, including macrophages and dendritic cells. 
These innate immune cells migrate from the lung to the
draining lymph node and activate naive T cells with the present of Mtb.
Activated effector and memory T cells then  travel back to the lung infection site, engage in granuloma formation, and control the infection.
The major elements of the host immune response against Mtb infection involve  macrophages, T lymphocytes, and Mtb bacteria.

\subsection{In-host Deterministic Model}
Our model describes
the host-pathogen dynamics for human Mtb infection in lung tissue.
Since the site of infection is the human lung, we focus on the interactions among Mtb bacteria (the pathogens), macrophages (Mtb ideal target cells), and lymphocytes (especially T cells for cell-mediated, cytotoxic adaptive immunity).
The measurement for all cells is taken in units per milliliter.
The 4-dimensional model \eqref{eqn1}  describes the dynamics
of the uninfected $M_u$ and infected $M_i$ macrophages,
the Mtb bacteria $B$, and the  CD4+ T cells $T$. In this model, we consider the effects from CD8+ T cells (cytotoxic T cells or CTL) and cytokines indirectly.
The model was developed by \cite{du2017simple}, analyzed by \cite{Zhang2020TB} and \cite{zhang2021investigation}.
  We present the model in \eqref{eqn1} followed by model descriptions.
Moreover, we write the parameter descriptions and values in Table \ref{tab1}.
\begin{equation}
\begin{array}{cl}
\dfrac{d M_u}{dt}&= s_M - \mu_M M_u - \beta M_u B \\[3.0ex]
\dfrac{d M_i}{dt}&= \beta M_u B - b M_i - \gamma M_i \dfrac{T/M_i}{T/M_i + c} \\[3.0ex]
\dfrac{dB}{dt} &= \delta B \left(1-\dfrac{B}{K}\right) + M_i \left( N_1 b + N_2 \gamma
\dfrac{T/M_i}{T/M_i +c} \right ) - M_u B (\eta + N_3 \beta) \\[3.0ex]
\dfrac{dT}{dt} &= s_T + \dfrac{c_M M_i T}{e_M T + 1} + \dfrac{c_B B T}{e_B T + 1} - \mu_T T.
\end{array}
\label{eqn1}
\end{equation}

{\color{blue}
\begin{table}[h!]
    \centering
    \begin{tabular}{|c c c |}
\hline
 Symbol & Description (Unites) & Value \\
\hline
$s_M$& recruitment rate of $M_u$ (1/ml day) &  $5000 $  \\
$s_T$& recruitment rate of $T$ (1/ml day) &  $6.6 $ \\
$\mu_M$& death rate of $M_u$ (1/day) &  $0.01 $  \\
$b$ & loss rate of $M_i$ (1/day) &  $0.11 $  \\
$\mu_T$ & death rate of $T$ (1/day) &  $0.33 $  \\
$\beta$ & infection rate by $B$ (1/day) &  $2\times 10^{-7} $  \\
$\eta$ & bacteria killing rate by $M_u$ rate (1/ml day) &  $1.25 \times 10^{-8}$   \\
$\gamma$ & cell-mediated immunity rate (1/day) &  $1.5 $  \\
$\delta$ & proliferation rate of $B$ (1/day) &  $5 \times 10^{-4}$  \\
$c_M$ & expansion rate of $T$ induce by $M_i$ (1/day) & $ 10^{-3} $  \\
$c_B$ & expansion rate of $T$ induce by $B$ (1/day) & $ 5 \times 10^{-3} $ \\
$e_M$ & saturating factor of $T$ expansion related  to $M_i$ &
$ 10^{-4} $  \\
$e_B$ & saturating factor of $T$ expansion related  to $B$ &
$ 10^{-4} $  \\
$c$ & half-saturation ratio for $M_i$ lysis ($T/M_i$) &
$3 $ \\
$K$ & carrying capacity of $B$ (1/ml) &
$ 10^8 $  \\
$N_1$ & max MOI of $M_i$  ($B/M_i$) &  $50 $  \\
$N_2$ & max No. of $B$ released by apoptosis (and necrosis)  ($T/M_i$) &  $20 $ \\
$N_3$ & $N_3=N_1/2$  ($B/M_i$) &  $25 $  \\ 
\hline
\end{tabular}
    \caption{Parameter Symbol, Descriptions, and Values (source: \cite{wigginton2001model}, \cite{marino2004human}, \cite{gammack2005understanding}, \cite{du2017simple}).}
    \label{tab1}
\end{table}
}

\paragraph{Macrophages.} Uninfected $M_u$ macrophages arrive at
the infection site at a constant recruitment rate of $s_M$ and 
have a finite lifespan with a natural death rate $\mu_M$.
With the presence of Mtb, uninfected macrophages engulf the bacteria $B$ and become chronically infected at a rate $\beta$  without receiving sufficient stimulation for activation.
The apoptosis and necrosis of infected macrophages could occur naturally, 
by bursting induced by the excessive intracellular bacterium load
caused by the phagocytized bacteria multiplication, or by receiving the death signal from CD8+ cells. The loss rate of infected macrophages due to natural death and bursting is denoted as $b$.

Since  the activation and proliferation of CD8+ T cells require the signal from CD4+ T cells, we consider the cytotoxic action (i.e., CD8+ T cells target intracellular bacteria) is proportional to CD4+ T cell function.
Therefore, the ratio of CD4+ T cell to infected macrophage ($T/M_i$) determines the rate of CD8+ T cells killing infected macrophages ($M_i$).
This rate reaches to its half-maximum at the value of $c$ and its maximum at $\gamma$.

\paragraph{Bacteria.} 
The bacterial population in TB infection is comprised of intracellular and extracellular bacteria.
Extracellular bacteria become intracellular or phagocytized through
engulfment by a host macrophage.
Intracellular bacteria are released and become extracellular
when infected macrophages undergo programmed cell death or killed by T cells.
We define extracellular bacteria as $B$, whose
 division is governed by a logistic term $\delta B(1-B/K)$
with a constant growth rate $\delta$ and a maximal carrying
capacity $K$.
The death of an infected macrophage releases $N_1$ intracellular bacteria per cell due to natural death and cell bursting, and
releases $N_2$ intracellular bacteria per cell due to CD8+ T cell cytotoxic action.
With adequate stimulation, uninfected macrophage engulfment can kill the phagocytized bacteria
 ($\eta M_u B$, assuming the killing rate  is $\eta$).  Otherwise, phagocytosis turns uninfected macrophages into chronically infected ($N_3 \beta B M_u$). We assume the engulfment process takes $N_3$ bacteria per  macrophage cell.

\paragraph{CD4+ T lymphocytes.}
CD4+ T cells play a role in signaling CD8+ T cell activation and replication to eliminate infected macrophages via cytotoxic action.
CD4+ T cells also produce cytokines for granuloma formation and infection control.
During the TB infection, professional antigen-presenting cells, including macrophages and dendritic cells, travel  to lung-draining lymph node to present Mtb antigens to naive T cells and activate them.
We assume that CD4+ T cells are activated by (1) infected macrophages with phagocytized bacteria at rate $c_M$ with the saturating factor $e_M$, and (2) dendritic cells (assumed in proposition to the extracellular bacteria) at rate $c_B$ with the saturating factor $e_B$. The natural recruitment and death rates of CD4+ T cells are denoted as $s_T$ and $\mu_T$.

\subsection{Basic Properties of the Deterministic Model}
As an airborne disease, TB is transmitted from person to person through airborne particles containing Mtb (i.e. droplet nuclei). 
The invasion of the infecting Mtb produces immune responses in the host.
Early clearance can be achieved if an effective innate immune response successfully eradicate the infecting Mtb before the development of an adaptive response. Otherwise, the infection stays and 
develops to either an active TB or LTBI. A delayed clearance can still be achieved thanks to adaptive immune responses.

 In the previous paper (\cite{Zhang2020TB}), we proved that the solution trajectories stay in a bounded positive quadrant for all positive time and carried out the steady state analysis.
We denote the steady state as $(\Bar{M}_u,\,\Bar{M}_i,\,\Bar{B},\,\Bar{T})$.
If uninfected macrophages phagocytose all extracellular  Mtb
and kill all the engulfed bacteria through phagocytosis, then  no chronically infected macrophages will be generated, and the early infection is eradicated. 
 That means the infection term ($\beta\,M_u B$) is $\frac{s_M}{\mu_M+\beta B}B \beta \frac{1}{b+\gamma}=M_i=0$. We then have $B=0$.
Here, $s_M/(\mu_M+\beta B)$ denotes the uninfected macrophage population on a quasi-steady state level, $1/(b+\gamma)$ is the minimum
average lifetime of an infected macrophage ($b$ the loss rate of infected macrophages and $\gamma$ the maximum infected macrophage elimination rate by T-cell mediated immune response).
The infection then reaches a trivial steady state, $(\Bar{M}_{u0},\,\Bar{M}_{i0},\,\Bar{B}_0,\,\Bar{T}_0) = (s_M/\mu_M, \,0,\,0,\,s_T/\mu_T)$. It denotes the clearance of resident bacteria by macrophages with the help of T-cell mediated immune response.
If $\frac{s_M}{\mu_M+\beta B}B \beta \frac{1}{b+\gamma}=M_i\neq 0$,
the existence of bacteria located anywhere other than inside infected macrophages is possible. 
This infected steady state is
       $\Bar{M}_u=\frac{s_M}{\mu_M+\beta B}$,  
      $\Bar{T}=\frac{\beta B \Bar{M}_u  - b \Bar{M}_i}{\beta B \Bar{M}_u - \left(b+\gamma \right) \Bar{M}_i} \, c \,\Bar{M}_i$,
      $\Bar{M}_i = \left(\frac{B\delta}{K}-\frac{\delta+s_M[(N_2-N_3)\beta-\eta]}{\beta B+\mu_M}\right)\frac{B}{b(N1-N2)}$, and
$ F(B)=-e_B e_M \mu_T \Bar{T}^3+ \left[(c_M \Bar{M}_i + e_M s_T - \mu_T) e_B + e_M (c_B B - \mu_T)\right]$ $\Bar{T}^2
+ \left[c_B B + c_M \Bar{M}_i + (e_B + e_M) s_T - \mu_T\right] \Bar{T} + s_T=0.$

\section {Two Stochastic Models}
 Randomness is a typical feature in immunology (\cite{wodarz2007killer}).
It prevents us from foreseeing the exact state at a given time during the process of immune responses.
However, if the same experiment is repeated a large number of times, we can obtain a certain trend with similar outcomes. For example, the solution of a linear deterministic model represents the expectation of the corresponding stochastic model. However, the solutions of deterministic models fail to provide information about the intrinsic variability in stochastic models. These random variations in immune responses could induce various disease and treatment outcomes in a model with multiple stable equilibriums.
We apply multivariate It\^{o} SDEs to model the dynamics of Mtb bacteria and immune cells' interaction by using the modeling algorithm based on a multivariate continuous Markov chain model (\cite{allen1999stochastic},  \cite{allen2007modeling},  \cite{allen2017primer},\cite{allen2008construction}, and \cite{allen2020real}).

Two stochastic models are formulated based on demographic variations and environmental variations, respectively. 
 The study of population changes of pathogens and immune cells can be viewed as the study of demography in population dynamics. 
The dynamics between pathogens and the immune system demonstrate predator-prey interactions. The immune cells, acting as predators, capture and kill their pathogen prey. During this interaction, the sizes of immune cell populations grow. This induces a faster decline in the pathogen population. In the meantime, Mtb pathogens proliferate and evolve survival strategies to avoid the predation from immune cells.
Moreover, the infection process turns healthy macrophage cells into infected cells.

We model the changes over time within the cell population processes, such as proliferation, death, immigration, and transition of both immune cells and Mtb bacteria as demographic variations. 
An SDE model with demographic variations captures the random nature of this cell population birth-death-immigration process.
Moreover, the noise-induced demographic variations can cause species collapse/recover between high- and low-abundance equilibriums (\cite{meng2020tipping}) in the ecology system.
Analogously, demographic variations could induce the transition between a low-infection state (such as LTBI) and a high-infection state (such as active disease) in the interplay of pathogens and immune system.
In addition, the inhabiting environment (such as the individual's physical factors) also changes and very likely affects the pathogen and host immune cell populations in other random manners, such as the infected macrophage loss rate, cell-mediated immunity rate, macrophages bacteria killing rate, and bacterial proliferation rate. 
Instead of including additional variables in the existing ODE model, an alternative way to include these environmental variations is modeling the affected parameters as random variables following certain stochastic processes (\cite{allen1999stochastic}). 
These environmental variations happen especially in drug therapy and alter the cellular environment with random manner due to the heterogeneity in individual cell response.
This motivates us to develop the second SDE model with both demographic and environmental variations.

\subsection{SDE Model with Demographic Variations}
To formulate stochastic processes from the underling deterministic model \eqref{eqn1}, 
the first step is to identify the possible interactions resulting in population changes
and their probabilities in a given small-time interval $\Delta t$.
The summary of these changes is in Table \ref{tab2}. 
Let $\overrightarrow{X}(t)=(M_u,\,M_i,\,B,\,T)^{tr}$, $\Delta \overrightarrow{X}=\overrightarrow{X}(t+\Delta t)-\overrightarrow{X}(t)$, and $t\in(0,\,\infty)$.
Notice that for the SDE model, it is not necessary to convert the units of continuous random variables, $M_u$, $M_i$, $B$, and $T$,  from concentration to population size (\cite{allen2010introduction}), since the changes $(\Delta \overrightarrow{X})_k$ do not need to be integer-valued.
It is assumed that intracellular Mtb bacteria are reproduced with an infected macrophage and released only upon the lysis of the infected macrophage. 
If the death of   one unit of the infected macrophage, denoted by $-1$, is due to the killing of the cell-mediated immunity, $N_1$   units of Mtb bacteria are released, $(\Delta \overrightarrow{X})_4=(0,-1,N_1,0)^{tr}$. 
If it is caused by excessive intracellular bacterium load, $N_2$   units of Mtb bacteria are released, $(\Delta \overrightarrow{X})_5=(0,-1,N_2,0)^{tr}$. 
Furthermore, a loss of extracellular Mtb bacteria can also be caused by phagocytosis of extracellular Mtb bacteria by uninfected macrophages, $(\Delta \overrightarrow{X})_7=(0,0,-1,0)^{tr}$.
All the possible changes in the small-time interval $\Delta t$ include these eleven changes in Table \ref{tab2} and a case of no change, $(\Delta \overrightarrow{X})_{12}=(0,0,0,0)^{tr}$, with the possibility $p_{12}\Delta t= 1-\sum^{11}_{k=1} p_k \Delta t$. Notice that the higher order of $\Delta t$ is neglected.
\begin{table}[h!]
    \centering
    \begin{tabular}{|c l l p{6cm}|}
\hline
$k$ & State Change $(\Delta \overrightarrow{X})_k$ & Probability
$p_i \Delta t + o(\Delta t)$& Description \\\hline
 1& $(1,0,0,0)^{tr}$ & $s_M \Delta t + o(\Delta t)$ & $M_u$ recruitment\\
 2& $(-1,0,0,0)^{tr}$ & $\mu_M M_u \Delta t + o(\Delta t)$ & $M_u$ death\\
 3& $(-1,1,0,0)^{tr}$ & $\beta M_u B \Delta t + o(\Delta t)$ & $M_u$ infected by $B$\\
 4& $(0,-1,N_1,0)^{tr}$ & $bM_i\Delta t + o(\Delta t)$ & loss of $M_i$ and release bacteria\\
 5& $(0,-1,N_2,0)^{tr}$ & $\gamma M_i \frac{T/M_i}{T/M_i+c}\Delta t + o(\Delta t)$ & $M_i$ killed by cell-mediated immunity and release bacteria \\
 6& $(0,0,1,0)^{tr}$ & $\delta B (1-B/K)\Delta t + o(\Delta t)$ &  Mtb bacteria proliferation\\
 7& $(0,0,-1,0)^{tr}$ & $M_u B (\eta+N_3 \beta) \Delta t + o(\Delta t)$ & extracellular loss due to $M_i$ killing and phagotysis of Mtb\\
 8& $(0,0,0,1)^{tr}$ & $s_T\Delta t + o(\Delta t)$ & T-cell recruitment\\
 9& $(0,0,0,1)^{tr}$ & $\frac{c_M M_i T}{e_M T+1}\Delta t + o(\Delta t)$ & T-cell activated by $M_i$\\
 10& $(0,0,0,1)^{tr}$ & $\frac{c_B B T}{e_B T+1}\Delta t + o(\Delta t)$ & T-cell activated by $B$\\
 11& $(0,0,0,-1)^{tr}$ & $\mu_T T\Delta t + o(\Delta t)$ & T-cell death\\
 \hline
\end{tabular}
    \caption{Possible state changes during $\Delta t$ and their associated probabilities.}
    \label{tab2}
\end{table}
The stochastic processes are formed from the underling deterministic ODE model \eqref{eqn1}.

In the second step of formulation, omitting the higher-order terms of $\Delta t$, we compute the first order approximation of the expectation $\mathbb{E}(\Delta \overrightarrow{X})$ and the covariance $\mathbb{E}[(\Delta \overrightarrow{X})(\Delta \overrightarrow{X})^{tr}]$.
The approximated expectation of $\Delta \overrightarrow{X}$ can be expressed as follows
\begin{equation*}
    \mathbb{E}\left(\Delta \overrightarrow{X}\right)= \sum^{11}_{k=1} p_k \Delta t \left(\Delta \overrightarrow{X}\right)_k= \overrightarrow{f}\left(M_u, M_i,B,T\right)\,\Delta t,
\end{equation*}
where
\begin{equation*}
   \overrightarrow{f}= 
    \left( \begin{array}{c} 
    f_1 \left(M_u, M_i,B,T\right)\\
     f_2 \left(M_u, M_i,B,T\right)\\
      f_3 \left(M_u, M_i,B,T\right)\\
       f_4 \left(M_u, M_i,B,T\right)
       \end{array}
   \right)=
   \left( \begin{array}{c}
  s_M - \mu_M M_u - \beta M_u B \\[2.0ex]
 \beta M_u B - b M_i - \gamma M_i \dfrac{T/M_i}{T/M_i + c} \\[2.0ex]
\delta B \left(1-\dfrac{B}{K}\right) + M_i \left( N_1 b + N_2 \gamma
\dfrac{T/M_i}{T/M_i +c} \right ) - M_u B (\eta + N_3 \beta) \\[2.0ex]
s_T + \dfrac{c_M M_i T}{e_M T + 1} + \dfrac{c_B B T}{e_B T + 1} - \mu_T T.
   \end{array}
   \right)
\end{equation*}
is the drift vector, which is the same as the right side of the deterministic ODE model \eqref{eqn1}. Here $p_i$s are coefficients of the $\Delta t$ terms of the probabilities for the eleven events in Table \ref{tab2}.

The covariance matrix is a $4 \times 4$ matrix, 
$\mathbb{E}[(\Delta \overrightarrow{X})(\Delta \overrightarrow{X})^{tr}]-\mathbb{E}(\Delta \overrightarrow{X})[\mathbb{E}(\Delta \overrightarrow{X})]^{tr}$. 
Dropping the $(\Delta t)^2$ order term $\mathbb{E}(\Delta \overrightarrow{X})[\mathbb{E}(\Delta \overrightarrow{X})]^{tr}$ results the  approximation to order $\Delta t$ as follows:
\begin{equation*}
    \mathbb{E}[(\Delta \overrightarrow{X})(\Delta \overrightarrow{X})^{tr}]\approx
    \sum^{11}_{k=1} p_k \Delta t \left(\Delta \overrightarrow{X}\right)_k\left(\Delta \overrightarrow{X}\right)^{tr}_k= \Sigma\left(M_u, M_i,B,T\right)\,\Delta t.
\end{equation*}
Here, the matrix $\Sigma$ takes the following form
\begin{equation*}
    \Sigma=
    \begin{bmatrix}
    p_1+p_2+p_3 & -p_3 &0&0\\
    -p_3 & p_3+p_4+p_5 & -N_1 p_4-N_2 p_5 & 0\\
    0 &-N_1 p_4-N_2 p_5 & N_1^2 p_4+N_2^2 p_5+p_6+p_7 & 0\\
    0&0&0&p_8+p_9+p_{10}+p_{11}
    \end{bmatrix}
    =\left[\Sigma_{i,j}\right]
\end{equation*}
The nonzero elements of $\Sigma$ are as follows:
\begin{equation*}
    \begin{array}{l}
        \Sigma_{11}=B M_u \beta+ M_u  \mu_M+ s_M,    \qquad \Sigma_{12}=-\beta  M_u  B, \\[2.0ex]
       \Sigma_{21}= -\beta  M_u  B,    \Sigma_{22}=\beta  M_u  B+b  M_i+\gamma  M_i  \dfrac{T/M_i}{T/M_i + c},  \\[2.0ex]
       \Sigma_{23}= -b  M_i  N_1-\gamma  M_i \dfrac{T/M_i}{T/M_i + c} N_2,  \\[2.0ex] 
     \Sigma_{32}= -b  M_i  N_1-\gamma  M_i  \dfrac{T/M_i}{T/M_i + c}N_2,  \\[2.0ex]
     \Sigma_{33}=b  M_i  N_1^2+\gamma  M_i  \dfrac{T/M_i}{T/M_i + c} N_2^2+ \delta  B  \left(1-\dfrac{B}{K}\right)+M_u  B  (N_3  \beta+\eta) \\[2.0ex]
    \Sigma_{44}=s_T+c_M  M_i  \dfrac{T}{T e_M+1}+c_B  B  \dfrac{T}{T e_B+1}+\mu_T  T.
    \end{array}
\end{equation*}
Instead of computing the square root of $\Sigma$ to find the diffusion matrix, we calculate a matrix $B\left(M_u, M_i,B,T\right)=[B_{i,j}]_{4\times 11}$, such that $BB^{tr}=\Sigma$. In the matrix $B$, four rows indicate four state variables ($M_u$, $M_i$, $B$, and $T$), eleven columns represent eleven events in Table \ref{tab2}. 
The nonzero elements in the matrix $B$ are as follows:
\begin{equation*}
\begin{array}{l}
     B_{11}=\sqrt{s_M}, \quad
      B_{12}=-\sqrt{\mu_M M_u}, \quad
       B_{13}=-\sqrt{\beta M_u B}, \\[2.0ex]
        B_{23}=\sqrt{\beta M_u B}, \quad
         B_{24}=-\sqrt{b M_i}, \quad
      B_{25}=-\sqrt{\dfrac{\gamma M_i T}{M_i c+T}},\\[2.0ex]
      B_{34}=N_1\sqrt{(b M_i)},\quad
      B_{35}=N_2\sqrt{\dfrac{\gamma M_i T}{M_ic+T}},\quad
      B_{36}=\sqrt{\delta B (1-\dfrac{B}{K})},\quad
      B_{37}=-\sqrt{M_u B (N_3 \beta+\eta)},\\[2.0ex]
       B_{48}=\sqrt{s_T},\quad
       B_{49}=\sqrt{\dfrac{c_M M_i T}{T e_M+1}},\quad
       B_{4,10}=\sqrt{\dfrac{c_B B T}{T e_B+1}},\quad 
       B_{4,11}=-\sqrt{\mu_T T}.
\end{array}
\end{equation*}
The It\^{o} SDE model considering the demographic variations
has the following form:
\begin{equation}
\left \{
\begin{array}{rl}
     dM_u(t)&= f_1(M_u,M_i,B,T) dt + B_{11} dW_1(t)+ B_{12} dW_2(t)+ B_{13} dW_3(t),\\[1.0ex]
     dM_i(t)&= f_2(M_u,M_i,B,T) dt + B_{23} dW_3(t)+ B_{24} dW_4(t)+ B_{25} dW_5(t),\\[1.0ex]
     dB(t)&=f_3(M_u,M_i,B,T) dt + B_{34} dW_4(t)+ B_{35} dW_5(t)+ B_{36} dW_6(t)+ B_{37} dW_7(t),\\[1.0ex]
     dT(t)&= f_4(M_u,M_i,B,T) dt + B_{48} dW_8(t)+ B_{49} dW_9(t)+ B_{4,10} dW_{10}(t)+ B_{4,11} dW_{11}(t),\\[1.0ex]
     \overrightarrow{X}(0)&=\left(M_u(0), M_i(0),B(0),T(0)\right), 
\end{array} \right.
    \label{sde1}
\end{equation}
where the vector  $W_i(t)$, $i=1\, ...\, 11$ are eleven independent Wiener processes, and $f_i$, $i=1\,...\,4$ are four elements in the drift vector $\overrightarrow{f}(M_u,M_i,B,T)$.

\subsection{SDE Model with Demographic and Environmental Variations}\label{sec_sde2}
Host-directed therapies can modulate host immune pathways. The therapies can also affect the corresponding parameters with random manner. We model these environmental effects by including additional random variables in the SDE model \eqref{sde1}.
 
In previous work (\cite{zhang2021investigation}), we identify four parameters, which significantly affect the disease progression.
These four parameters are the bacterial proliferation rate $\delta$, the infected macrophages' loss rate $b$, the cell-mediated immunity rate $\gamma$ and macrophages' killing rate $\eta$. 
We then consider environmental variations on them.
We adopt a mean-reverting process (\cite{allen2007introduction,allen2010introduction}), which was modeled to study drug administrations and was shown to fit well with real data.
Let $\overrightarrow{C}(t)=\left(C_1(t),C_2(t),C_3(t),C_4(t)\right)^{tr}$, where $C_1(t)=\delta(t)$, $C_2(t)=b(t)$, $C_3(t)=\gamma(t)$, and $C_4(t)=\eta(t)$, $C_s=(C_{s1},C_{s2},C_{s3},C_{s4})=(\delta_s,b_s,\gamma_s,\eta_s)$, and $C_0=(C_{01},C_{02},C_{03},C_{04})=(\delta_0,b_0,\gamma_0,\eta_0)$, we assume the vector $\overrightarrow{C}(t)$ is a vector of continuous, mean-reverting, stochastic processes.
Each element of $\overrightarrow{C}(t)$ satisfies an It\^{o} stochastic differential equation as follows:
\begin{equation*}
    C_i(t)=\alpha_i (C_{si}-C_i(t)) dt +\sigma_i C_i(t) d W_{ei}(t), \quad \text{for}\; i=1,2,3,4,
\end{equation*}
where $\alpha_i>0$ is the return rate to the mean concentration $C_{si}$. Here, $C_{si}$ denotes the desired therapy result, which is modulated by drug therapies. The timescale for $(C_{si}-C_i(t))$ returning to one half of its original level is $ln(2)/\alpha_i$. Therefore, a larger return rate value $\alpha_i$ implies a shorter return time to the desired therapy outcome. 
Parameters $\alpha_i$ and $C_{si}$ can be tuned through drug administration strategies. 
Moreover, $\sigma_i$ denotes the variability of the corresponding process, and $W_e(t)=(W_{e1}(t), W_{e2}(t),W_{e3}(t),W_{e4}(t))$ is a vector of four independent Wiener processes.
The asymptotic mean and variance of $C_i$ are derived by \cite{allen2007introduction} and \cite{allen2010introduction} as follows: 
\begin{equation*}
   \lim_{t\to \infty} \mathbb{E}(C_i(t)|C_{0i})=C_{si} \;\text{and}\;
   \lim_{t\to \infty} Var(C_i(t)|C_{0i})= \dfrac{C_{si}^2\sigma_i^2}{2\alpha_i-\sigma_i^2} \;\text{for}\; \alpha_i>\dfrac{\sigma_i^2}{2};\; \infty \;\text{for}\; \alpha_i\leq \dfrac{\sigma_i^2}{2}.
\end{equation*}
Therefore, to avoid a large variability and maintain a steady-state rate of immunity in therapies, the return rate should be sufficiently large, such that $\alpha_i>\frac{\sigma_i^2}{2}$. 
For the simulations in the next section, we assume $\alpha_i=2\sigma_i^2$, which yields $ \lim_{t\to \infty} Var(C_i(t)|C_{0i})= \frac{C_{si}^2}{3}$. Therefore, under this assumption, the limiting mean and variance are independent of the return rate $\alpha_i$, but only associated with the mean rate $C_{si}$.
The system of It\^{o} SDEs with both demographic and environmental variations is written as follows:
\begin{equation}
    \left \{
\begin{array}{rl}
     dM_u(t)&= f_1 dt + B_{11} dW_1(t)+ B_{12} dW_2(t)+ B_{13} dW_3(t),\\[1.0ex]
     dM_i(t)&= f_2 dt + B_{23} dW_3(t)+ B_{24} dW_4(t)+ B_{25} dW_5(t),\\[1.0ex]
     dB(t)&=f_3 dt + B_{34} dW_4(t)+ B_{35} dW_5(t)+ B_{36} dW_6(t)+ B_{37} dW_7(t),\\[1.0ex]
     dT(t)&= f_4 dt + B_{48} dW_8(t)+ B_{49} dW_9(t)+ B_{4,10} dW_{10}(t)+ B_{4,11} dW_{11}(t),\\[1.0ex]
     d\delta(t)&=\alpha_1 (\delta_s-\delta(t))dt+\sigma_1 \delta(t) dW_{e1}(t),\\[1.0ex]
     db(t)&=\alpha_2 (b_s-b(t))dt+\sigma_2 b(t) dW_{e2}(t),\\[1.0ex]
     d\gamma(t)&=\alpha_3 (\gamma_s-\gamma(t))dt+\sigma_3 \gamma(t) dW_{e3}(t),\\[1.0ex]
     d\eta(t)&=\alpha_4 (\eta_s-\eta(t))dt+\sigma_4 \eta(t) dW_{e4}(t),
\end{array} \right.
    \label{sde2}
\end{equation}
where $f_i=f_i(M_u(t),\,M_i(t),\,B(t),\,T(t),\,\delta (t),\,b(t),\,\gamma (t),\,\eta (t))$, and
$B_{i,j}=B_{i,j}(M_u(t),\,M_i(t),\,B(t),$ $T(t),\,\delta (t),\,b(t),\,\gamma (t),\,\eta (t))$, for $i=1...4$ and
$j=1..11$.
Here, $\left(M_u(0), M_i(0),B(0),T(0),\delta_0,b_0,\gamma_0,\eta_0\right)$ denotes
 initial conditions, and $W_i(t)$ and $W_{ej}(t)$, $i=1,...,11$ and $j=1,...,4$ denote independent Wiener processes.

In the following section, numerical methods are used to approximate the solutions of the proposed SDE models.
In our simulations, we apply a straightforward Euler-Maruyama method, and set zero as an absorbing state.

\begin{figure}[h!]
    \centering
    \begin{subfigure}{.7\textwidth}
        \centering
    \includegraphics[width=1\textwidth]{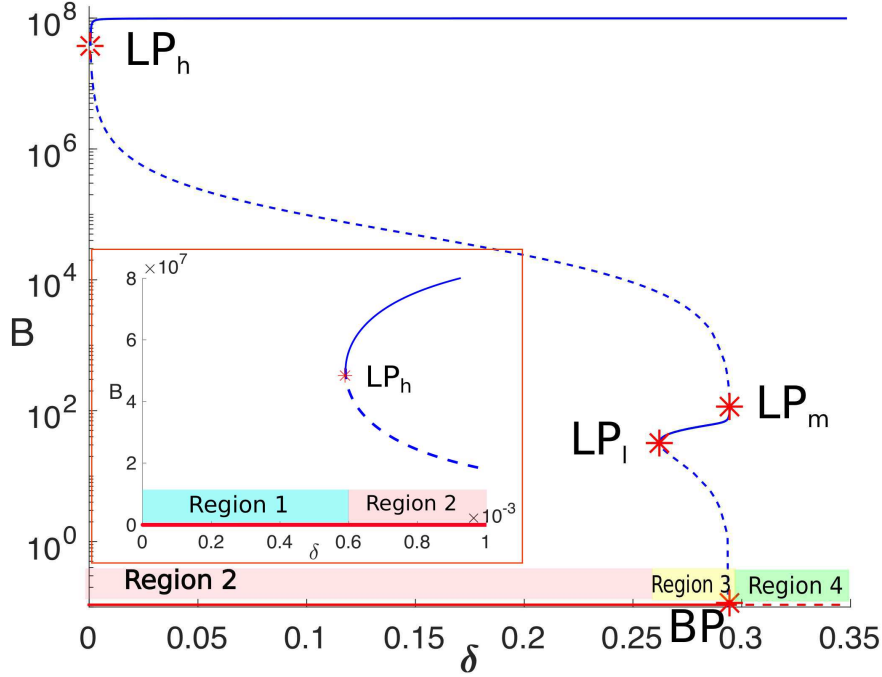}
        \caption{Bifurcation diagram of $B$ vs $\delta$.}
    \end{subfigure}
        \begin{subfigure}{.48\textwidth}
        \centering
    \includegraphics[width=1\textwidth]{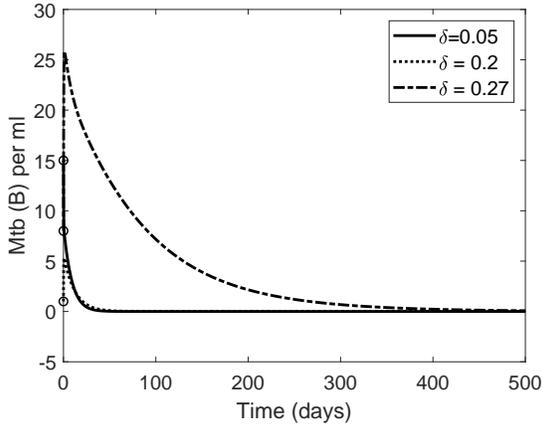}
        \caption{Simulated disease clearance.}
    \end{subfigure}
            \begin{subfigure}{.48\textwidth}
        \centering
    \includegraphics[width=1\textwidth]{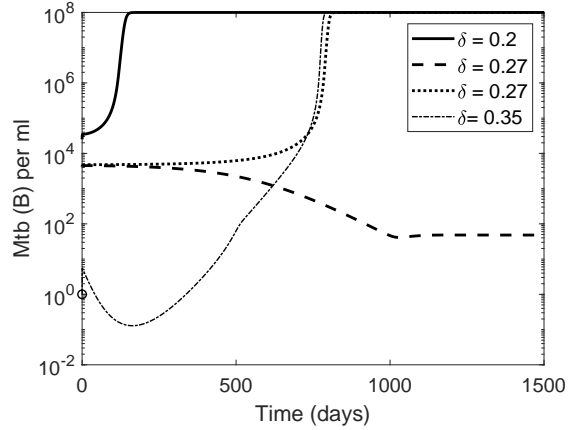}
        \caption{Simulated LTBI and active disease.}
    \end{subfigure}
    \caption{
    Bacterial concentrations under various bacterial proliferation rate. (a) Bifurcation diagram of the Mtb concentration ($B$ per ml) vs the Mtb proliferation rate ($\delta$ per day).
    Region 1 is demonstrated in a closeup with $\delta \in (0,\, 1 \times 10^{-3})$. 
    Sub-figures (b) and (c) demonstrate simulated clearance, LTBI, and active disease. 
   Parameter values are chosen as $\delta=0.05$ in Region 1, $\delta=0.2$ in Region 2,  $\delta=0.27$ in Region 3, and $\delta=0.35$ in Region 4. 
   Notice that simulated trajectories with $\delta=0.2$ and $\delta=0.27$ demonstrate multiple stable steady states. We circle the point $B(0)$ for four trajectories in (b) and (c) and take the other parameter values from Table \ref{tab1}.}
    \label{fig_bif_sim_B}
\end{figure}

\section{Mtb-host Dynamics of the ODE Model and SDE Model with demographic variation}\label{var_bac}
Mtb pathogens are extremely adaptive within a host. 
Their survival strategies involve multiple underlying mechanisms. 
Among these mechanisms, we focus on the bacterial proliferation rate $\delta$, which is the parameter most affected by the first-line antibiotic drugs.
We apply numerical bifurcation analyses to delimit the parameter space into four regions. 
Then we investigate disease clearance, latent infection, and active disease on the ODE model \eqref{eqn1} in the four identified parameter regions. 
Furthermore,  we study the stochastic fluctuation in cellular levels in each parameter region through approximated stationary distributions from simulations of the SDE model \eqref{sde1}. 

\subsection{Multiple Disease Outcomes Derived from Bifurcation Analysis}
Bacterial concentration reflects the local bacterial load, 
which is an important variable associated with Mtb virulence (\cite{smith2003mycobacterium}). The relation between bacterial concentration  and its proliferation rate $\delta$ is demonstrated
through a 1-dimensional bifurcation diagram in Figure \ref{fig_bif_sim_B} (a). The trivial steady state denoting disease clearance, $(\Bar{M}_{u0},\,\Bar{M}_{i0},\,\Bar{B}_0,\,\Bar{T}_0) = (s_M/\mu_M, \,0,\,0,\,s_T/\mu_T)$, exists for all positive values of $\delta$.
But it loses its stability when $\delta$ increases and crosses
a branching point (BP) bifurcation (transcritical bifurcation)
at $\delta_{BP}=0.2956$ /day, where
\begin{equation}
(\Bar{M}_{u0},\,\Bar{M}_{i0},\,\Bar{B}_0,\,\Bar{T}_0) =  (5\times 10^{5},\,0,\,0,\,20) \quad \text{cells/ml.}
    \label{triv_eqili}
\end{equation}
Without infection, a reference value for activated antigen-specific CD4+ T cells expressing CD25 is reported as $7-52$ cells/ml (\cite{bisset2004reference}), and a reference value for resting and activated macrophages is in the order of $10^5$ cells/ml (\cite{antony1993recruitment}  and \cite{schwander1998enhanced}).
The branching point bifurcation generates a disease established steady state $(\Bar{M}_u,\,\Bar{M}_i,\,\Bar{B},\,\Bar{T})$, which denotes the case that intracellular Mtb escape macrophages killing and become extracellular.
The bacterial concentration can stabilize at a lower level around $10^1-10^2$ cells/ml, if $\delta$ takes a value between two saddle node (LP) bifurcations LP$_l$ and LP$_m$.
The bacterial concentration can also reach a high-level range $10^7-10^8$ cells/ml, if $\delta$ takes a value between the saddle-node bifurcation LP$_h$ and the transcritical bifurcation BP.   A summary is as follows:
\begin{itemize}
    \item Bistability happens between LP$_h$ and LP$_l$. 
Taking $\delta=0.2\;\in [\delta_h,\,\delta_l]$ for example, clearance occurs if the initial infection contains a low-level of Mtb bacteria, illustrated by the dotted curve (with $M_u(0)= 1 \times 10^6$, $M_i(0)= 1$, $B(0)=1$, and $T(0)= 40$) in Figure \ref{fig_bif_sim_B} (b). 
Active disease happens in the case of a high-level of bacteria in the initial infection, represented by the solid curve (with $M_u(0)= 3 \times 10^6$, $M_i(0)= 2 \times 10^3$, $B(0)=2.5 \times 10^4$, and  $T(0)= 3.7 \times 10^6$) in Figure \ref{fig_bif_sim_B} (c).
    \item Triple stable steady states occur between LP$_l$ and LP$_m$. 
Taking $\delta=0.27\;\in [\delta_l,\,\delta_m]$ for example, disease clearance, LTBI, and active disease are demonstrated by the dash-dot curve (with $M_u(0)= 6\times 10^5$, $M_i(0)= 1$, $B(0)=8$,  and $T(0)= 90$) in Figure \ref{fig_bif_sim_B} (b), dashed curve (with $M_u(0)= 4.5\times 10^6$, $M_i(0)= 270$, $B(0)=4.2\times 10^3$,  and $T(0)= 7 \times 10^6$) in Figure \ref{fig_bif_sim_B} (c), and dotted curve (with $M_u(0)= 4.5\times 10^6$, $M_i(0)= 270$, $B(0)=4.8\times 10^3$,  and $T(0)= 7 \times 10^6$) in Figure \ref{fig_bif_sim_B} (c).
\item
The infection will certainly reach to a high-level if $\delta >\delta_m$, since $\delta_{BP} < \delta_m$.
Taking $\delta=0.35$ for example, a simulated trajectory starting from a low Mtb level and ending at high Mtb level, which is plotted in the dash-dot curve (with $M_u(0)= 1 \times 10^6$, $M_i(0)= 1$, $B(0)=1$,  and $T(0)= 40$) in Figure \ref{fig_bif_sim_B} (c).
\item Disease clearance will certainly be achieved if $\delta$ $<$ LP$_h$,
Taking $\delta=0.05$, a simulated trajectory showing Mtb bacteria die out is the solid curve (with $M_u(0)= 1 \times 10^6$, $M_i(0)= 1$, $B(0)=15$, and  $T(0)= 40$)  in Figure \ref{fig_bif_sim_B} (b).
\end{itemize}

The corresponding bifurcation points are LP$_l$: $(\Bar{M}_u,\,\Bar{M}_i,\,\Bar{B},\,\Bar{T},\,\delta_l)=(4.99672\times 10^5,\,2.366,\,32.78,\,40.16,\,0.2621)$, LP$_m$: $(\Bar{M}_u,\,\Bar{M}_i,\,\Bar{B},\,\Bar{T},\,\delta_m)=(4.98848\times 10^5,\,7.171,\,115.4,\,7751,\,0.2945)$, and LP$_h$: $(\Bar{M}_u,\,\Bar{M}_i,\,\Bar{B},\,\Bar{T},\,\delta_h)=(499.7,\,3102,\,$
$4.997\times 10^7,\,0.7572\times 10^{10},\,0.00059)$.
Hopf bifurcation does not occur for the parameter values in Table \ref{tab1}.
Moreover, saddle node bifurcation indicates the creation and destruction of equilibriums. Transcritical bifurcation, in this case, is the intersection of disease-free and infected equilibriums. These two bifurcations are sufficient to predict various disease outcomes, which include clearance, latent infection, and active disease.  We, therefore, do not consider higher codimension bifurcations in this project.
Instead, we divide the parameter range of the bacterial proliferation $\delta$ into four regions shown in Figure \ref{fig_bif_sim_B} (a), i.e. Region 1: $[0,\,\delta_h]$, Region 2: $[\delta_h,\,\delta_l]$, Region 3: $[\delta_l,\,\delta_m]$, and Region 4: $[\delta_m,\,+\infty]$. 
In Region 1, only stable disease-free equilibrium exists. The ODE model predicts that the infection will die out eventually. 
In Regions 2 and 3, the ODE model contains multiple stable equilibriums. Depending on the initial invading bacterial level, multiple disease outcomes may occur.
In Region 4, the ODE model has one stable infected equilibrium. It implies that the infection will grow to a high-level eventually.

\begin{figure}
    \centering
    \begin{subfigure}{0.45\textwidth}
    \includegraphics[width=1\textwidth]{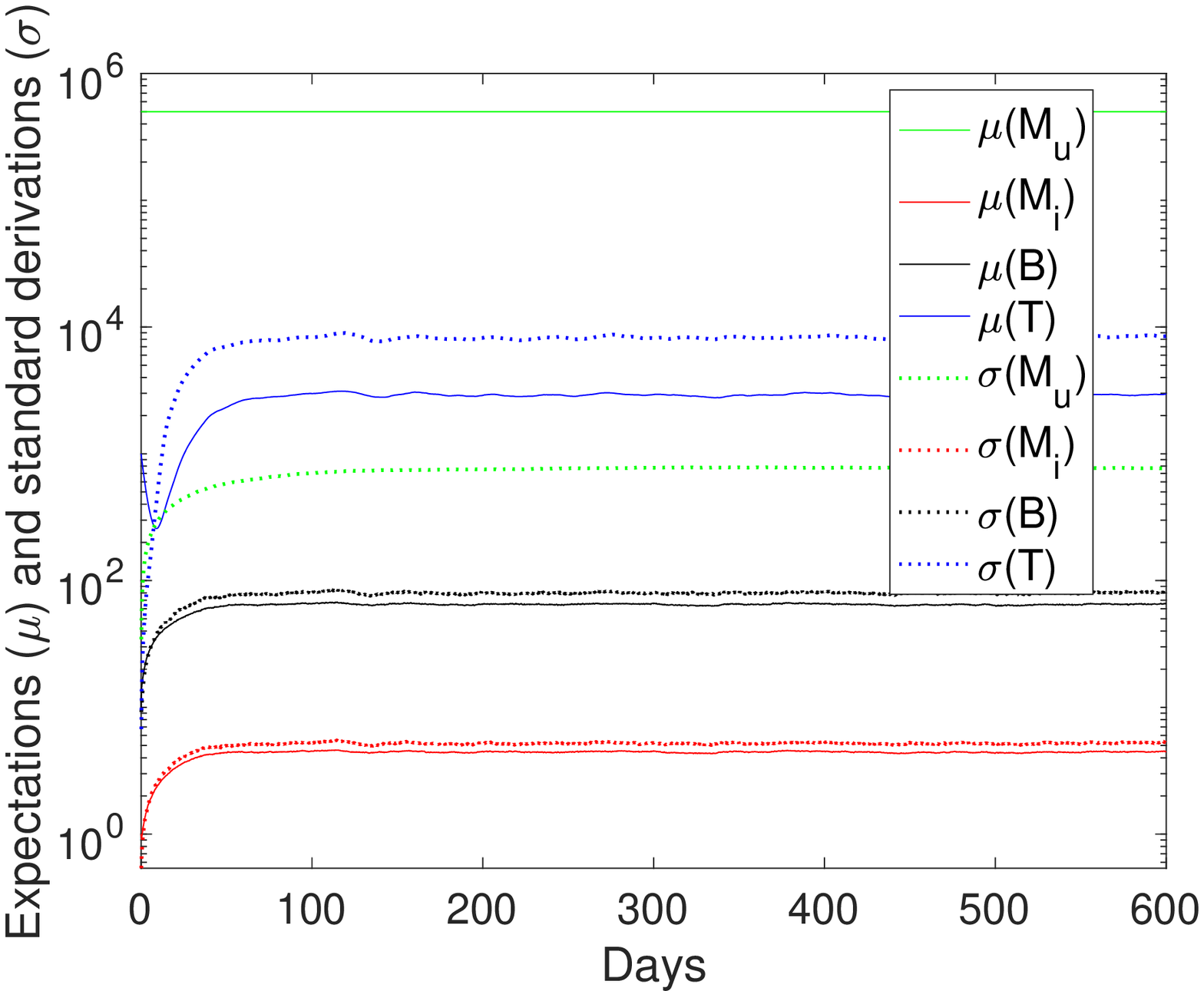}
    \caption{Mean and Std over time.}
        \label{fig_sd1a}
    \end{subfigure}
     \begin{subfigure}{0.42\textwidth}
    \includegraphics[width=1\textwidth]{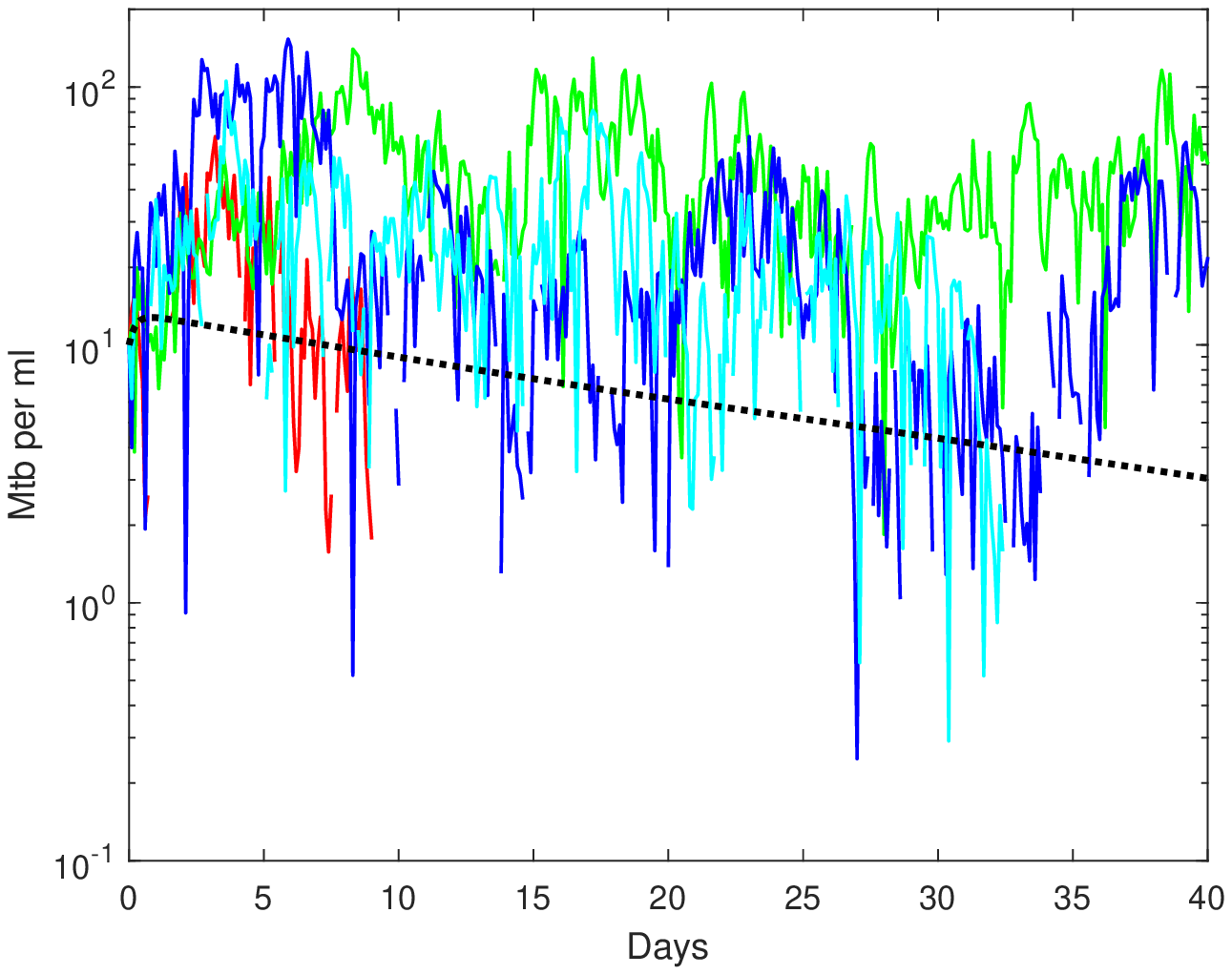}
    \caption{Colored sample paths, black ODE soln.}
         \label{fig_sd1b}
    \end{subfigure}\\
     \begin{subfigure}{0.45\textwidth}
        \includegraphics[width=1\textwidth]{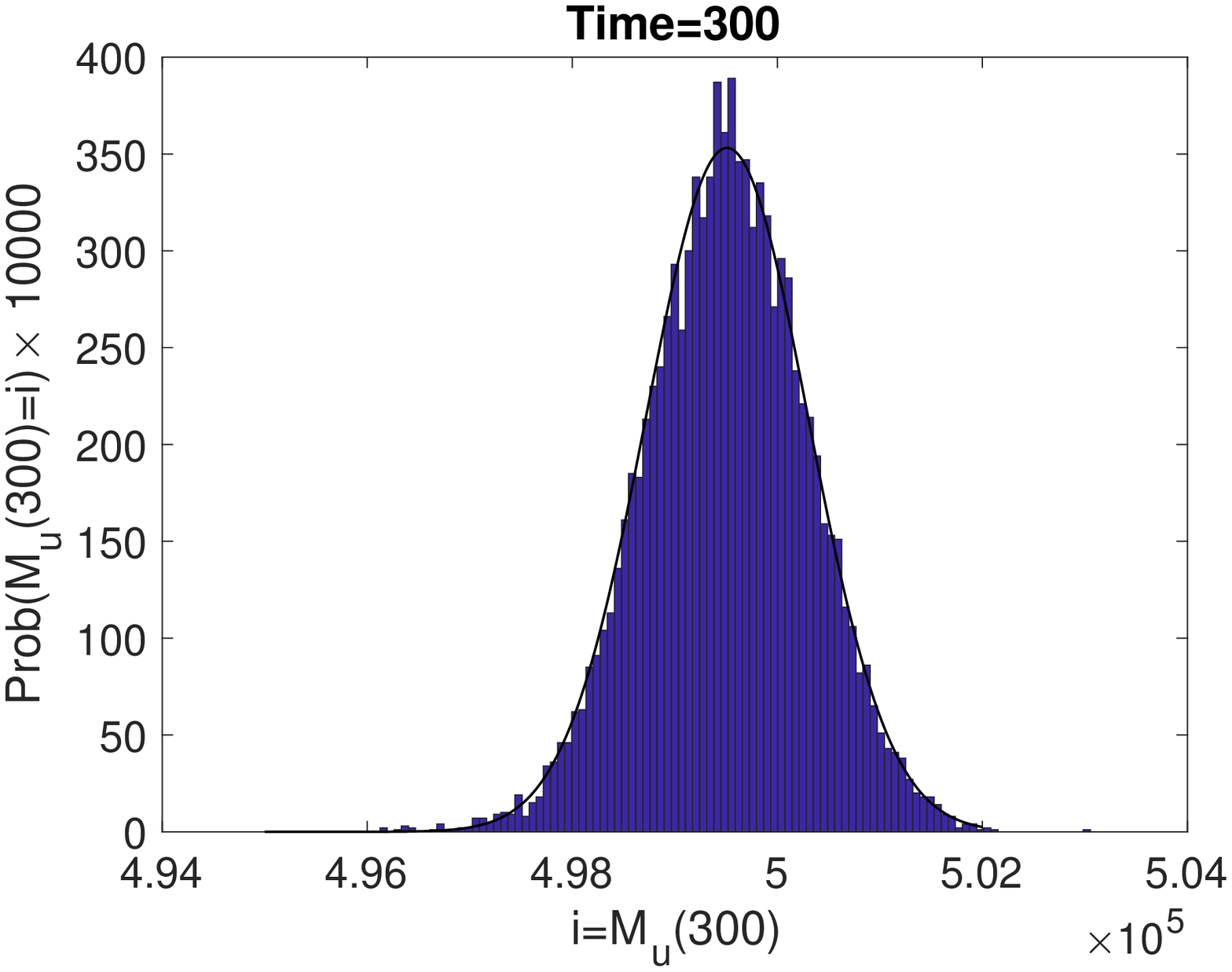}
        \caption{$\mu=4.9951\times 10^{5}$, $\sigma=790.8575$.}
        \label{fig_sd1c}
    \end{subfigure}
 \begin{subfigure}{0.45\textwidth}
    \includegraphics[width=1\textwidth]{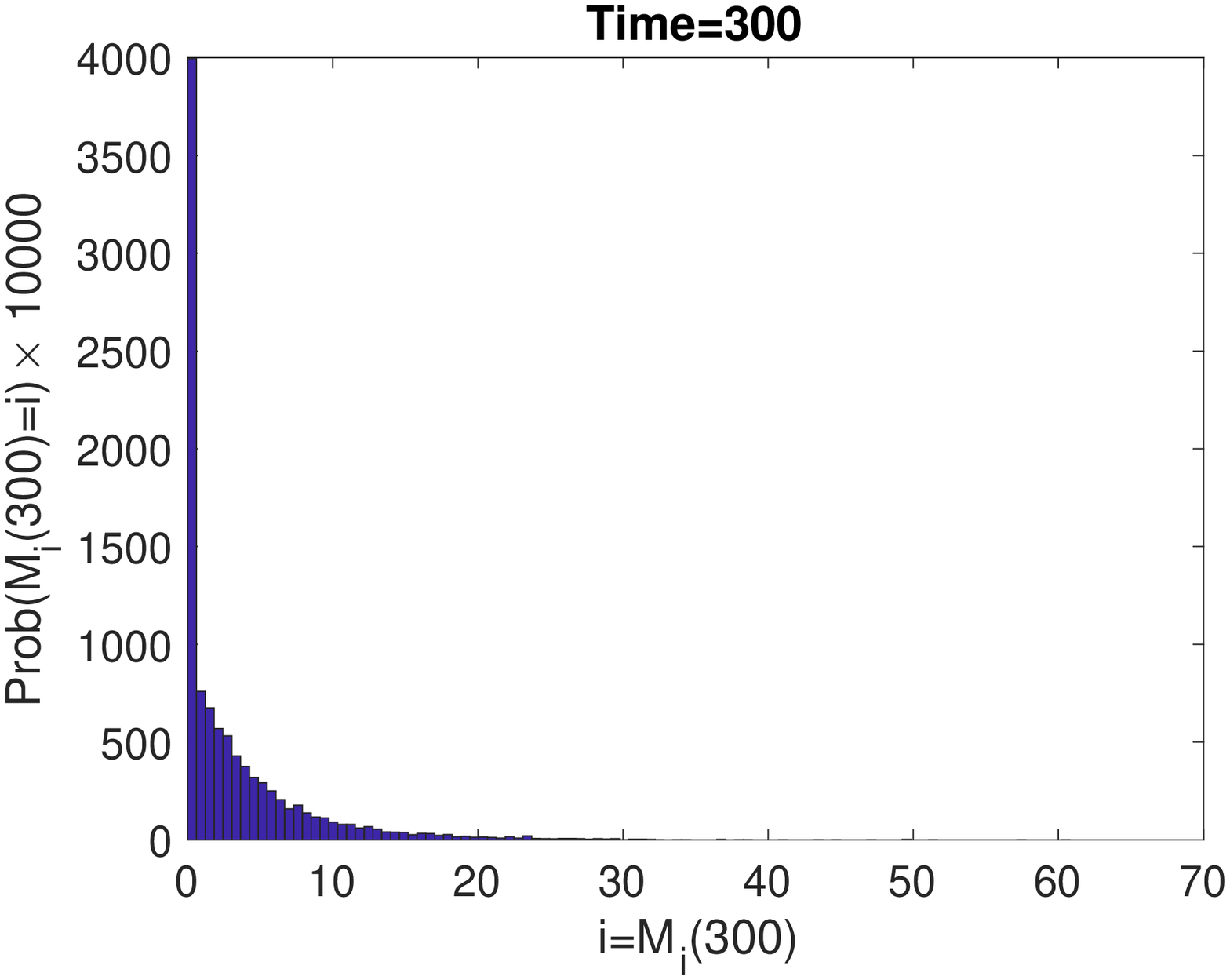}
             \label{fig_sd1d}
      \caption{$\mu=3.2950$, $\sigma=4.9366$.}
    \end{subfigure}
 \begin{subfigure}{0.45\textwidth}
        \includegraphics[width=1\textwidth]{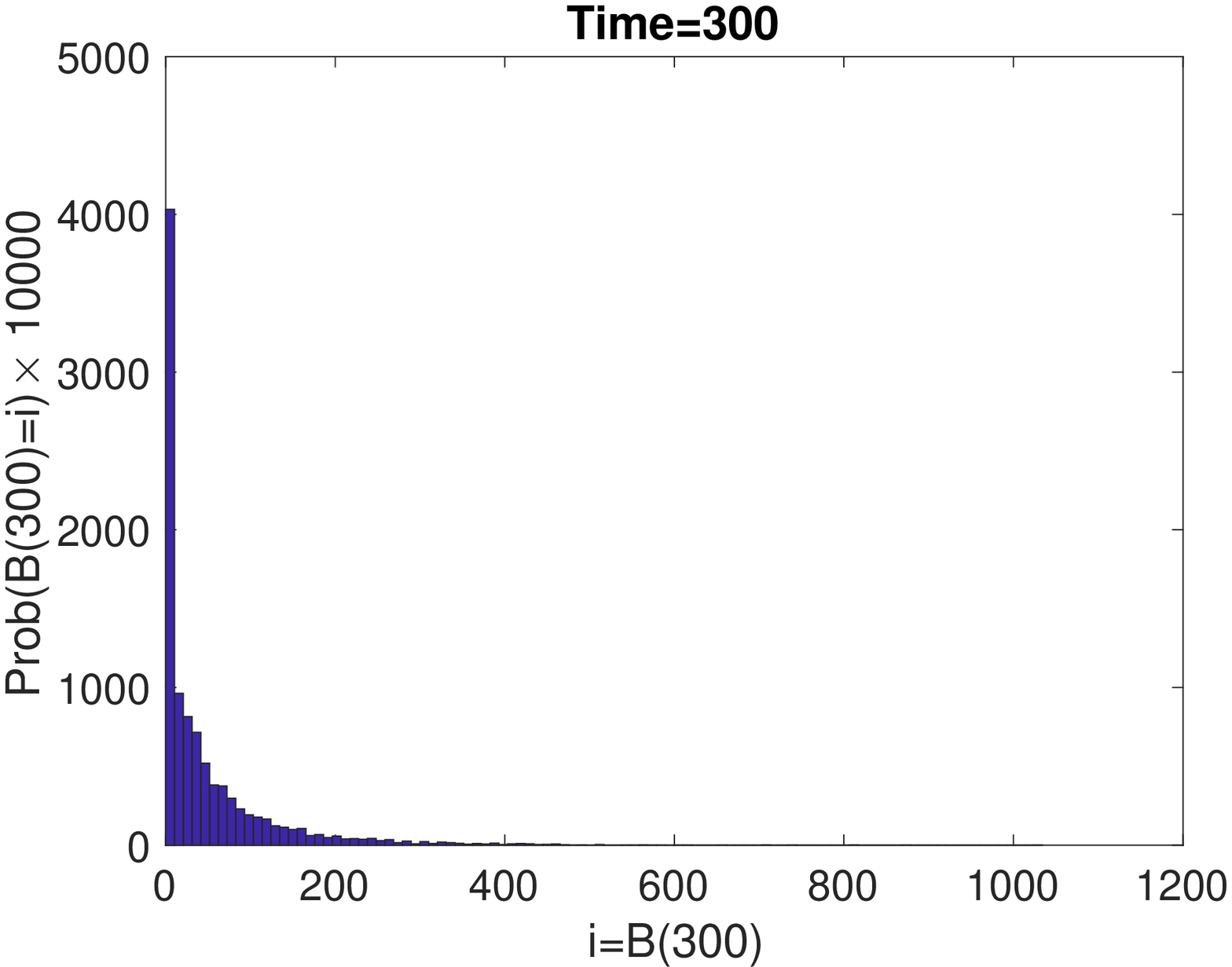}
               \label{fig_sd1e}
        \caption{$\mu=48.6576$, $\sigma=76.0503$.}
    \end{subfigure}
 \begin{subfigure}{0.45\textwidth}
    \includegraphics[width=1\textwidth]{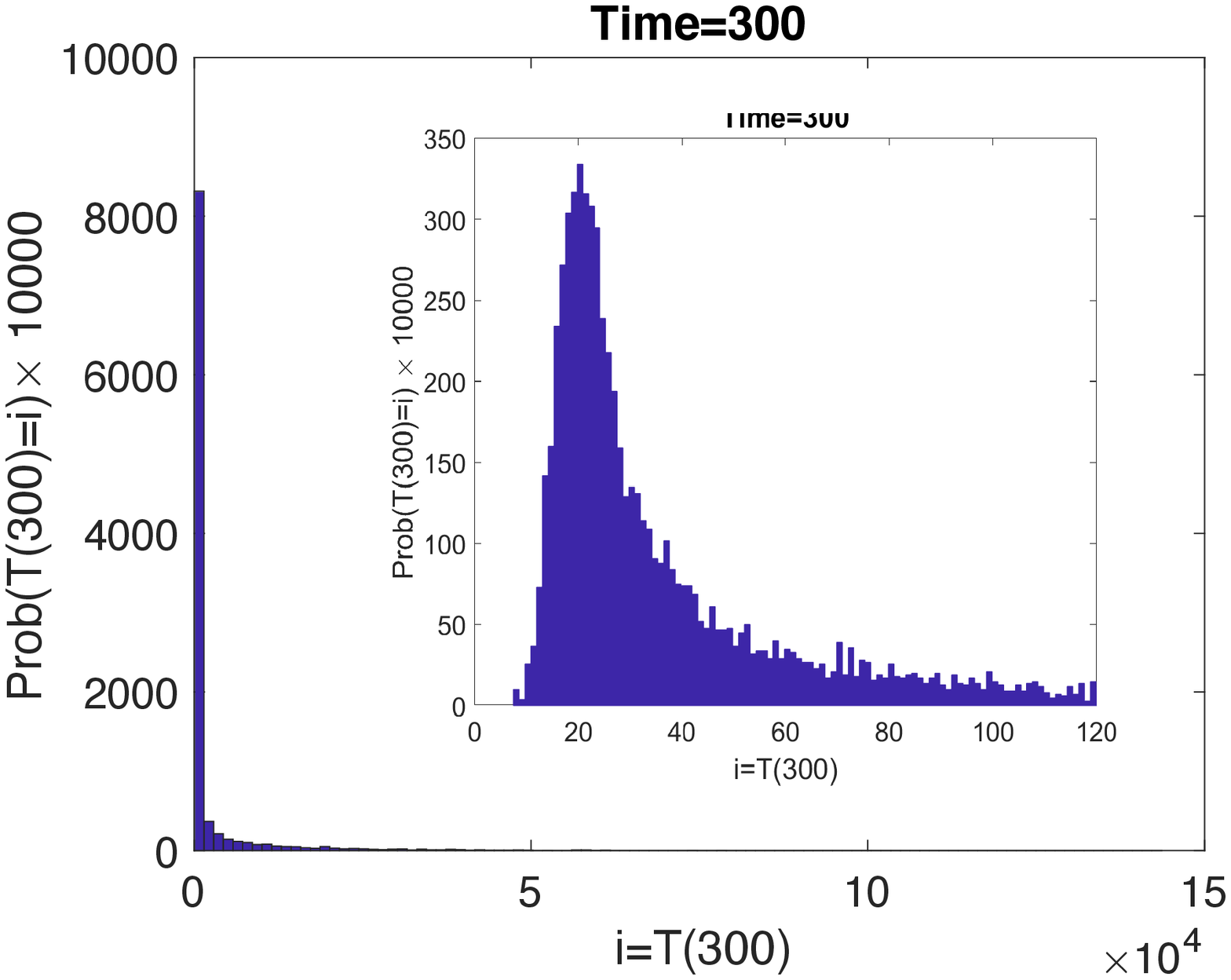}
        \caption{$\mu=2.2029\times 10^{3}$, $\sigma=7.4485\times 10^3$.}
       \label{fig_sd1f}
    \end{subfigure}
        \caption{Disease clearance or latent TB infection case at $\delta=0.2$. Mean: $\mu$, standard deviation: $\sigma$. (a): Expectations and standard deviations of 10,000 sample paths of the SDE model \eqref{sde1}. (b) Four stochastic realizations are in colored curves.   Since y-axis is in log scale, the red and cyan curves terminate when they hit the absorption state $B=0$. The corresponding ODE solution is in the dotted black curve. (c)-(f): Approximate stationary distributions for $M_i$, $M_u$, $B$, and $T$ at $t=300$.  $100$ bins are used for all four histograms. Initial conditions are $M_u(0)=4.99 \times 10^5$, $M_i(0)=1$, $B(0)=10$, and $T(0)=1000$. Note that the other parameter values are taken from Table \ref{tab1}.}
    \label{fig_sde1}
\end{figure}

\begin{figure}
    \centering
    \begin{subfigure}{0.45\textwidth}
    \includegraphics[width=1\textwidth]{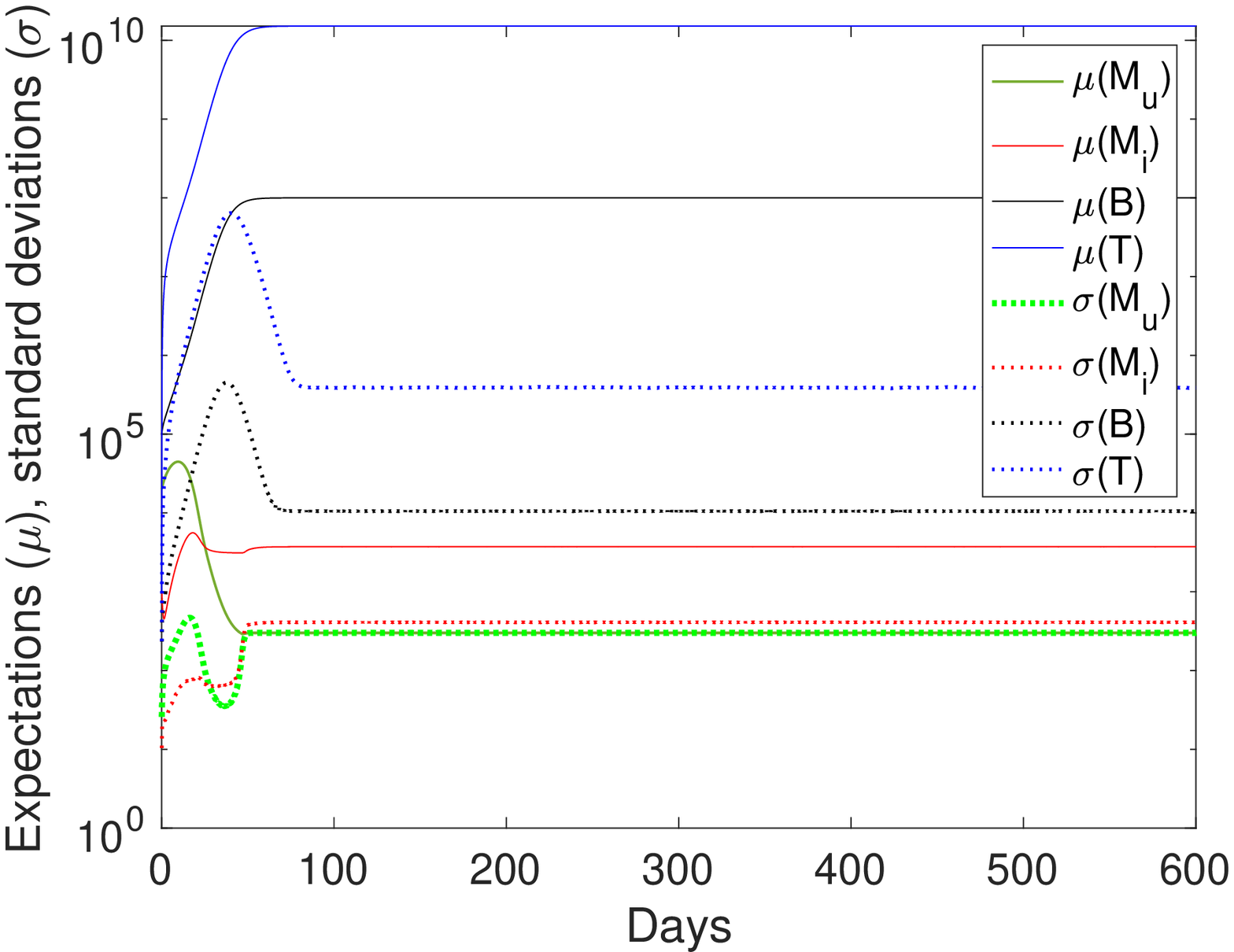}
    \caption{Mean and Std over time.}
        \label{fig_sd2a}
    \end{subfigure}
     \begin{subfigure}{0.45\textwidth}
    \includegraphics[width=1\textwidth]{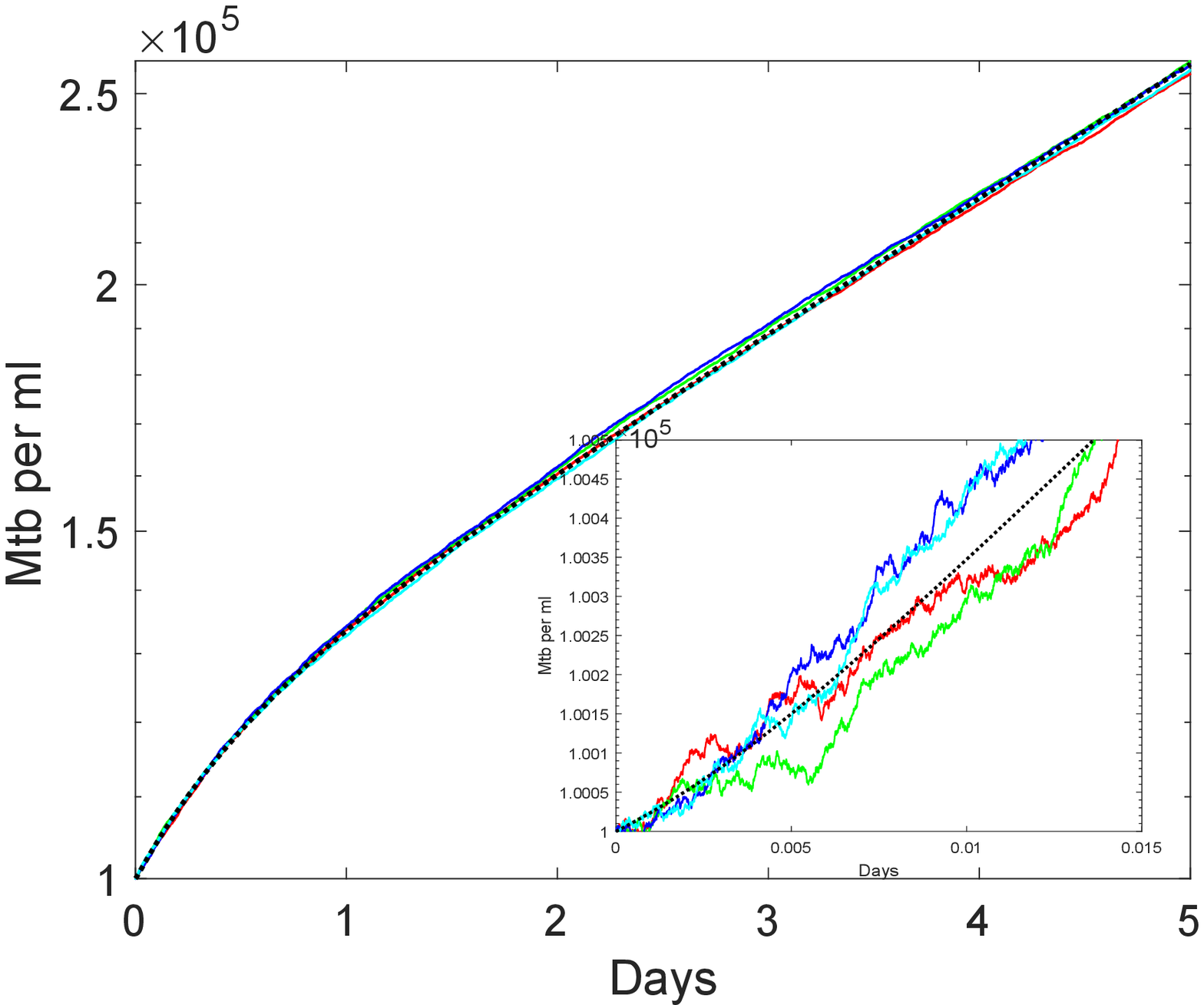}
    \caption{Colored sample paths, black ODE soln.}
         \label{fig_sd2b}
    \end{subfigure}\\
        \vspace{1cm}
     \begin{subfigure}{0.45\textwidth}
        \includegraphics[width=1\textwidth]{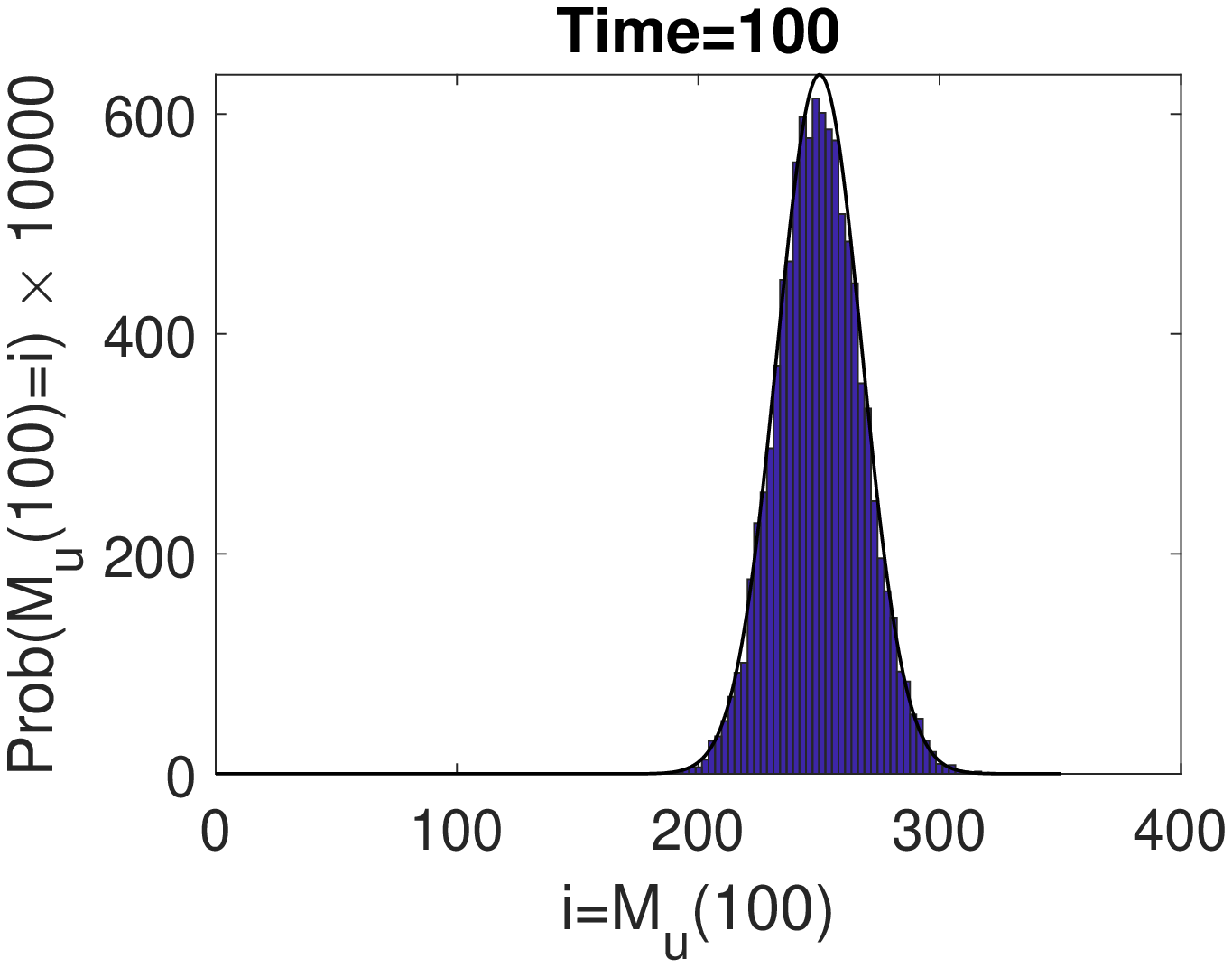}
        \caption{$\mu=250.2504$, $\sigma=17.5698$.}
        \label{fig_sd2c}
    \end{subfigure}
 \begin{subfigure}{0.45\textwidth}
    \includegraphics[width=1\textwidth]{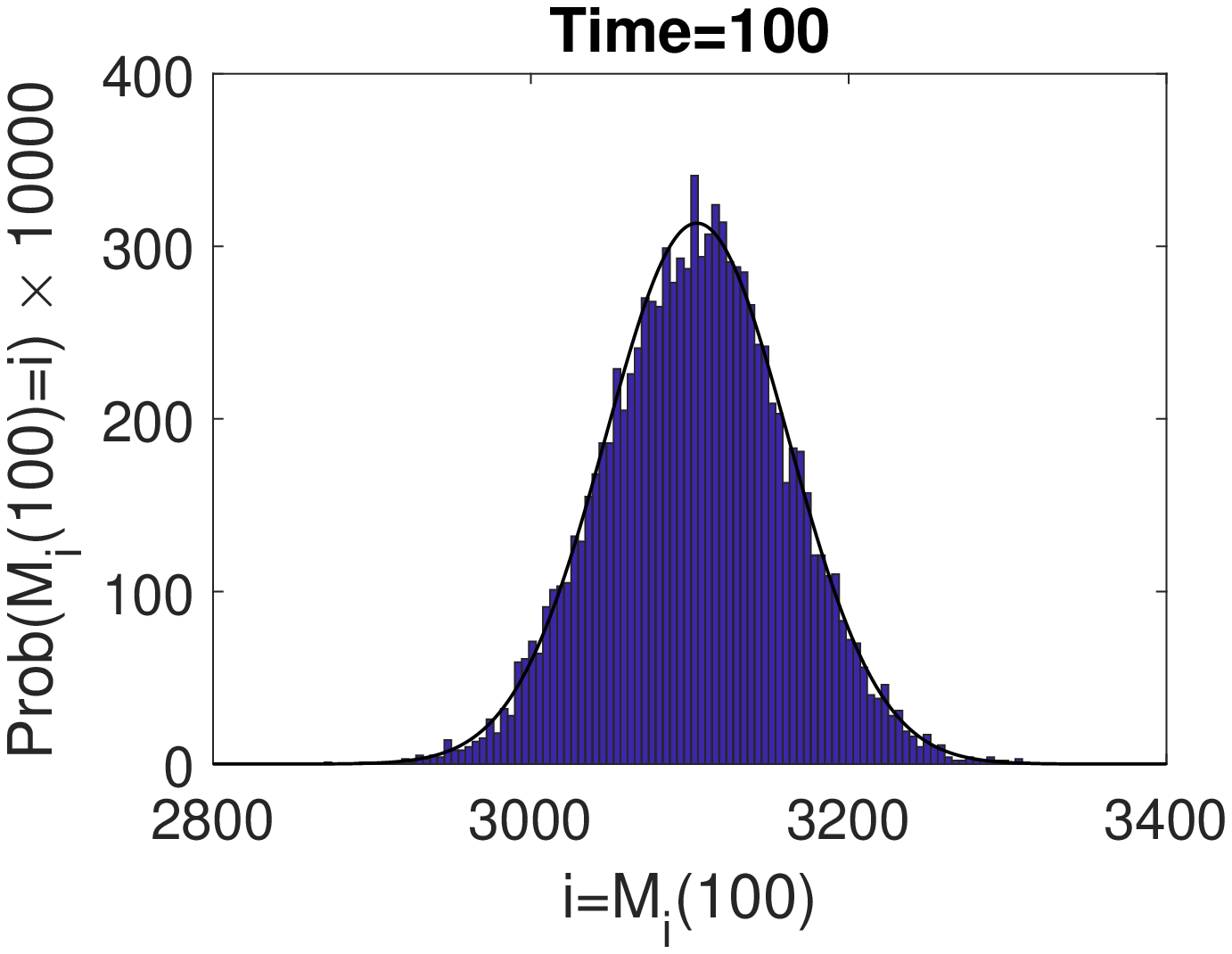}
             \label{fig_sd2d}
      \caption{$\mu=3.1045 \times 10^3$, $\sigma=57.2839$.}
    \end{subfigure}\\
        \vspace{1cm}
 \begin{subfigure}{0.45\textwidth}
        \includegraphics[width=1\textwidth]{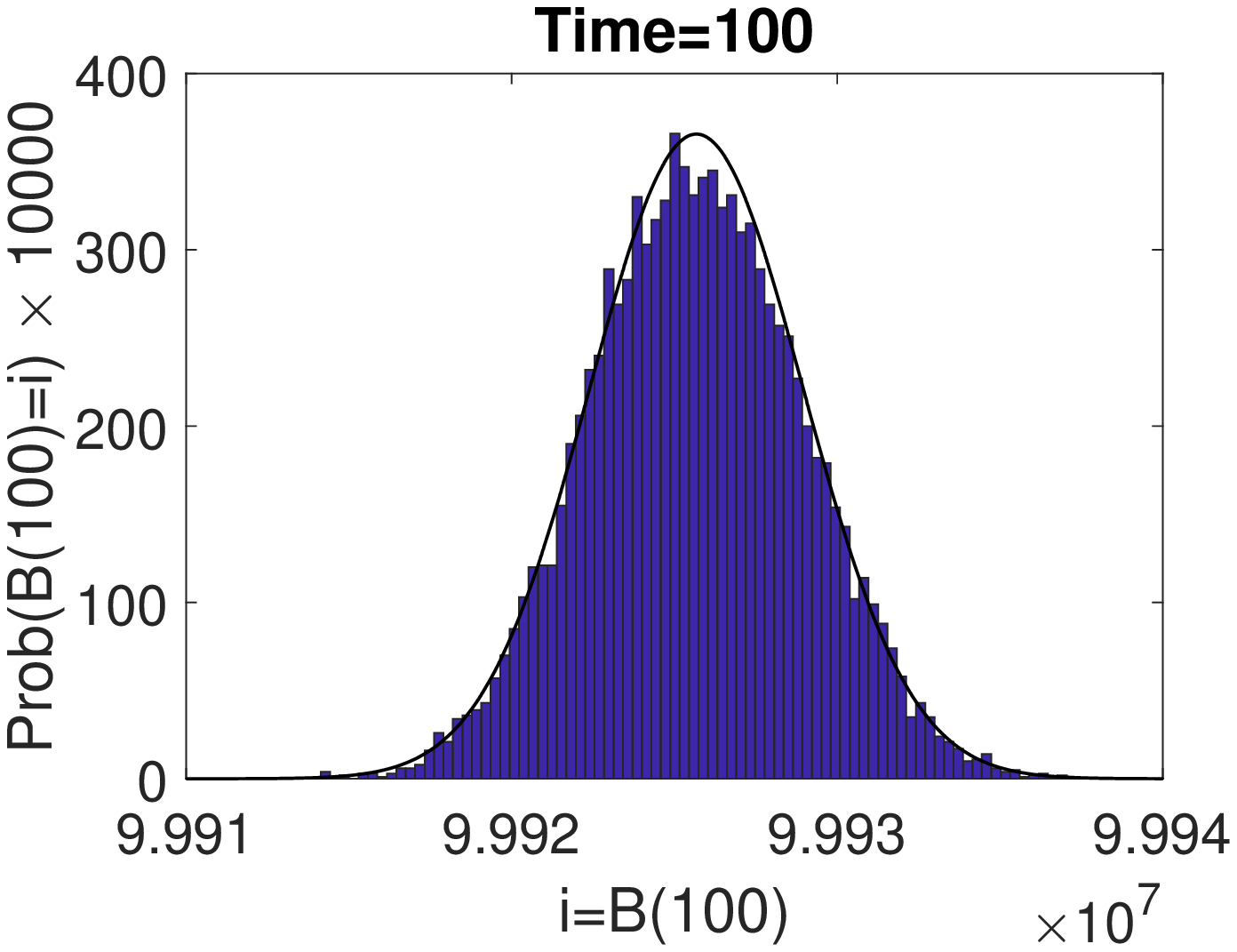}
               \label{fig_sd2e}
        \caption{$\mu=9.9926 \times 10^7$, $\sigma=3.2734 \times 10^3$.}
    \end{subfigure}
 \begin{subfigure}{0.45\textwidth}
    \includegraphics[width=1\textwidth]{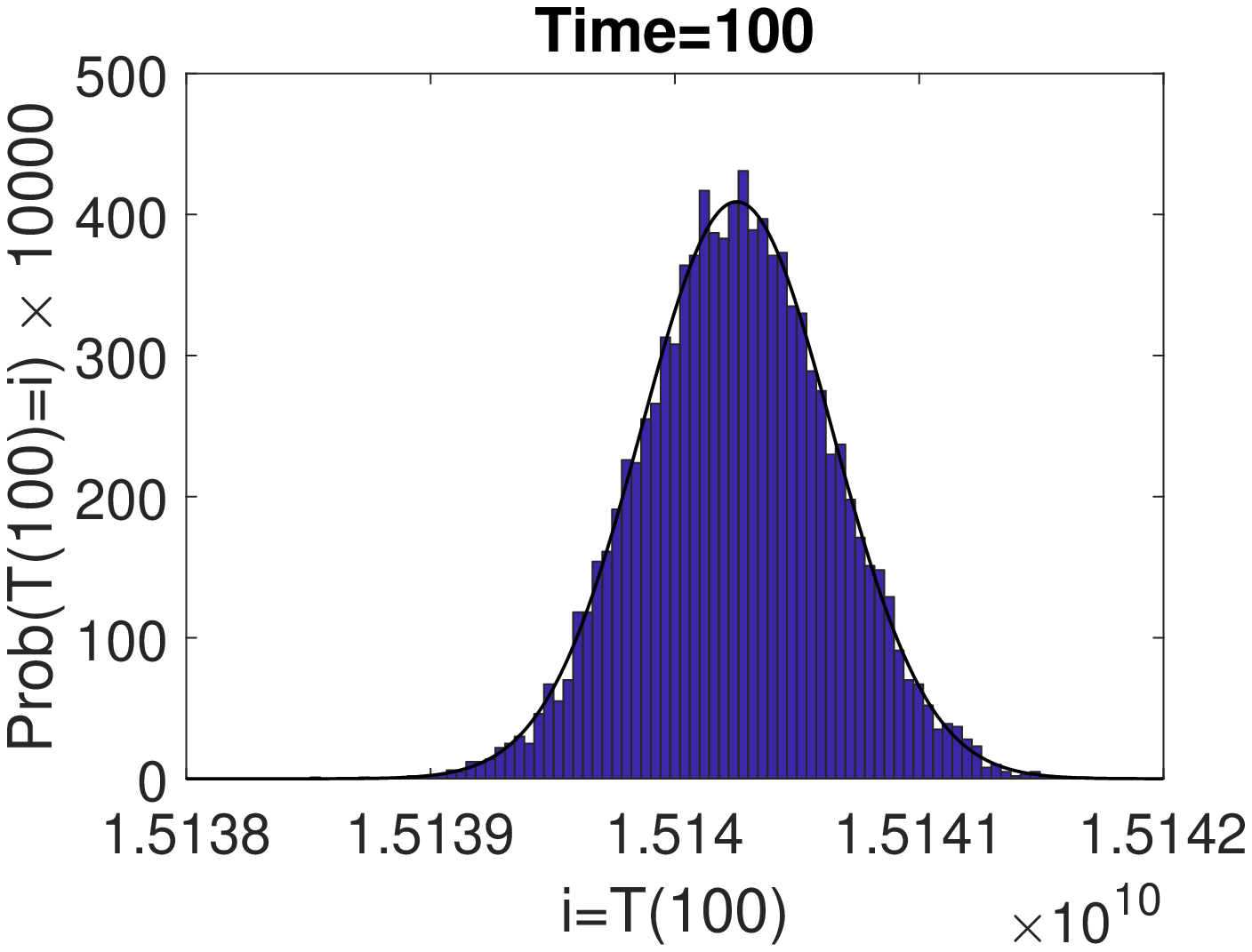}
        \caption{$\mu=1.5140\times 10^{10}$, $\sigma=3.9010\times 10^5$.}
       \label{fig_sd2f}
    \end{subfigure}
    \caption{Active disease case at $\delta=0.2$. Mean: $\mu$, standard deviation: $\sigma$. (a): Expectations and standard deviations of 10,000 sample paths of the SDE model \eqref{sde1}.   (b) Four stochastic realizations in colored curves vary around the corresponding ODE solution in the dotted black curve. Disease clearance doesn't occur due to the large number of bacterial load. (c)-(f): Approximate stationary distributions for $M_i$, $M_u$, $B$, and $T$ at $t=100$.  $100$ bins are used for all four histograms. Initial conditions are $M_u(0)= 2 \times 10^4$, $M_i(0)=10^3$, $B(0)=10^5$, and $T(0)=1000$. Note that the other parameter values are taken from Table \ref{tab1}.}
    \label{fig_sde2}
\end{figure}

\begin{figure}
    \centering
    \begin{subfigure}{0.45\textwidth}
    \includegraphics[width=1\textwidth]{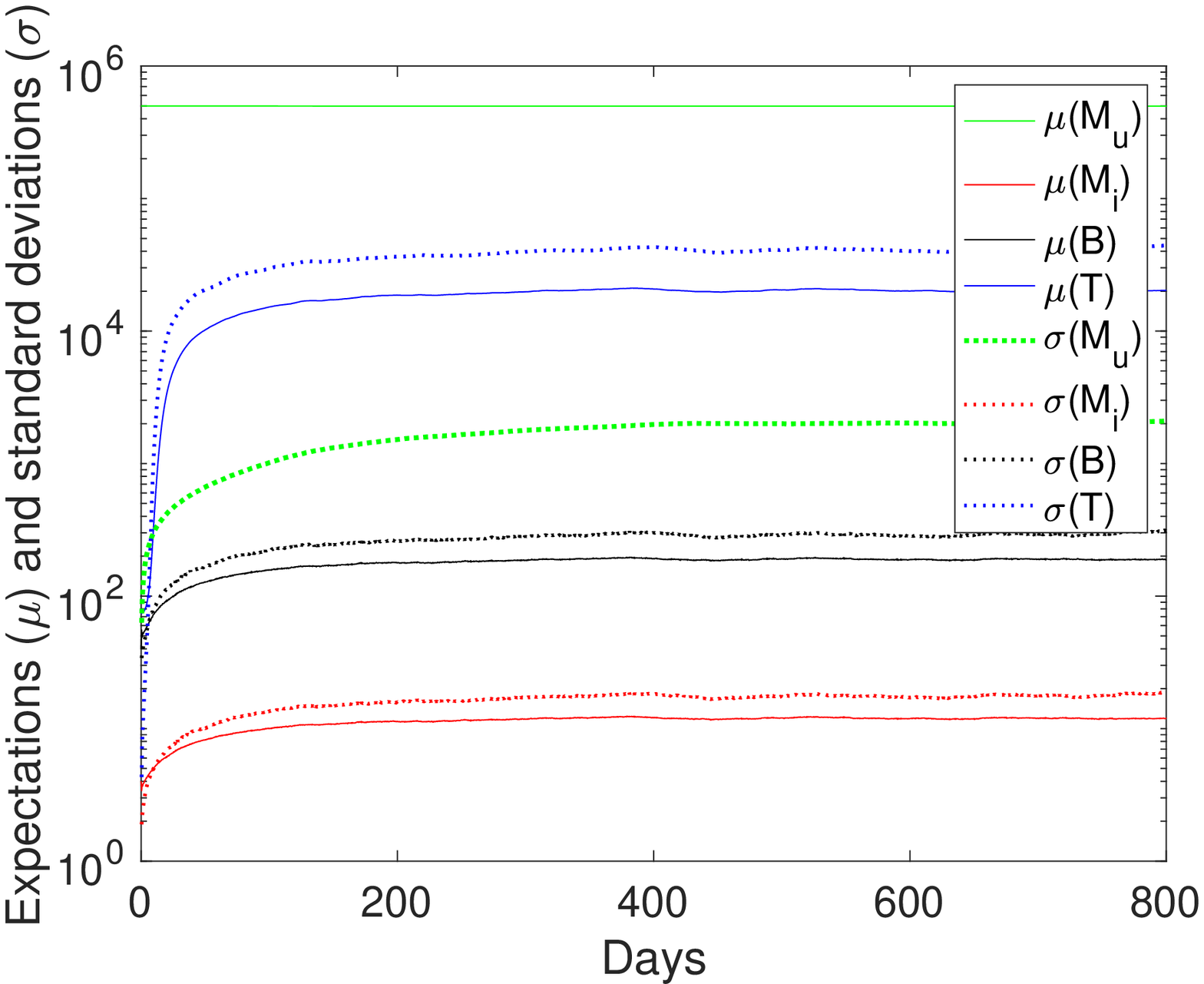}
    \caption{Mean and Std over time.}
        \label{fig_sd3a}
    \end{subfigure}
     \begin{subfigure}{0.45\textwidth}
    \includegraphics[width=1\textwidth]{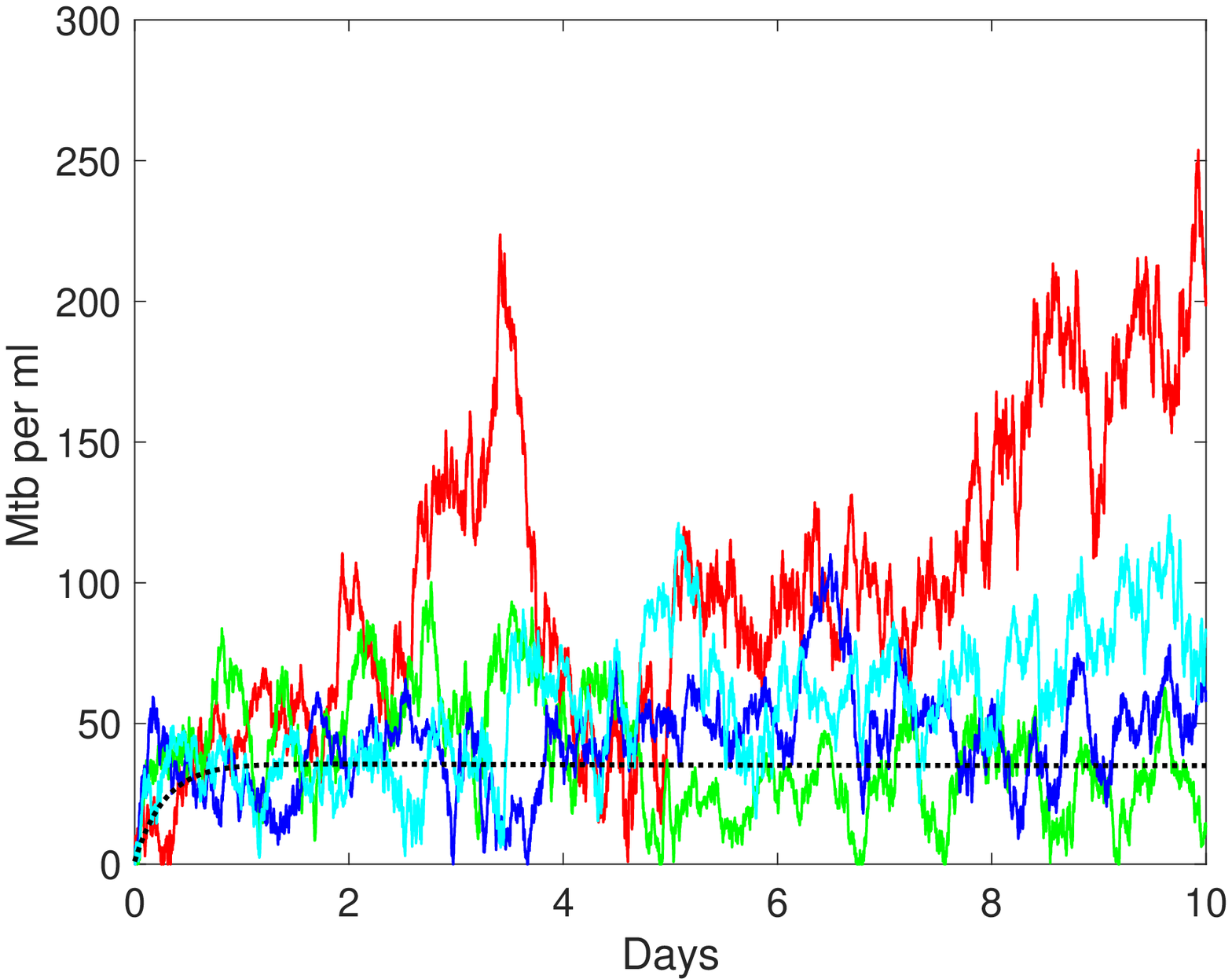}
    \caption{Colored sample paths, black ODE soln.}
         \label{fig_sd3b}
    \end{subfigure}\\
     \begin{subfigure}{0.45\textwidth}
        \includegraphics[width=1\textwidth]{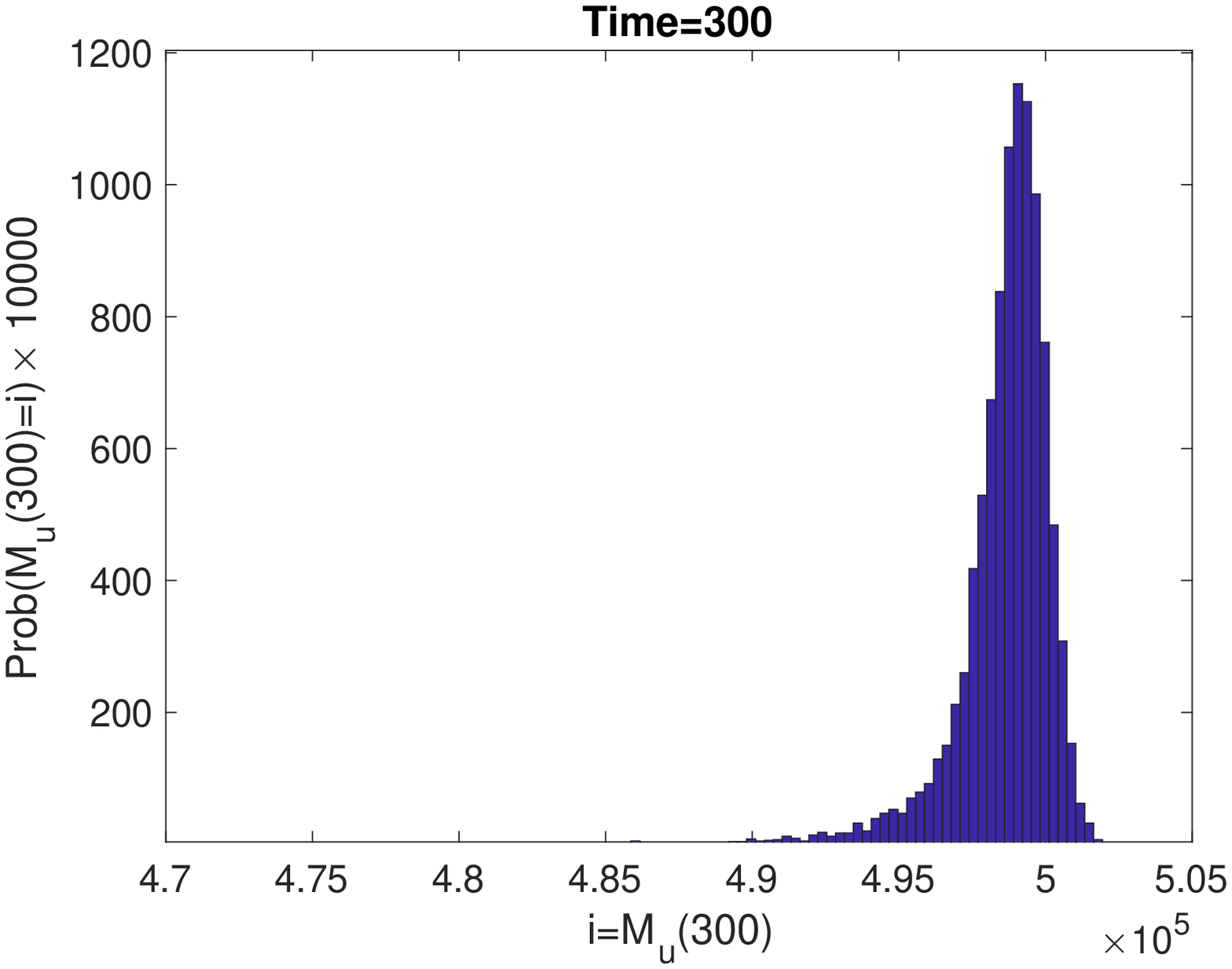}
        \caption{$\mu=4.9860\times 10^{5}$, $\sigma=1.7721\times 10^3$.}
        \label{fig_sd3c}
    \end{subfigure}
 \begin{subfigure}{0.45\textwidth}
    \includegraphics[width=1\textwidth]{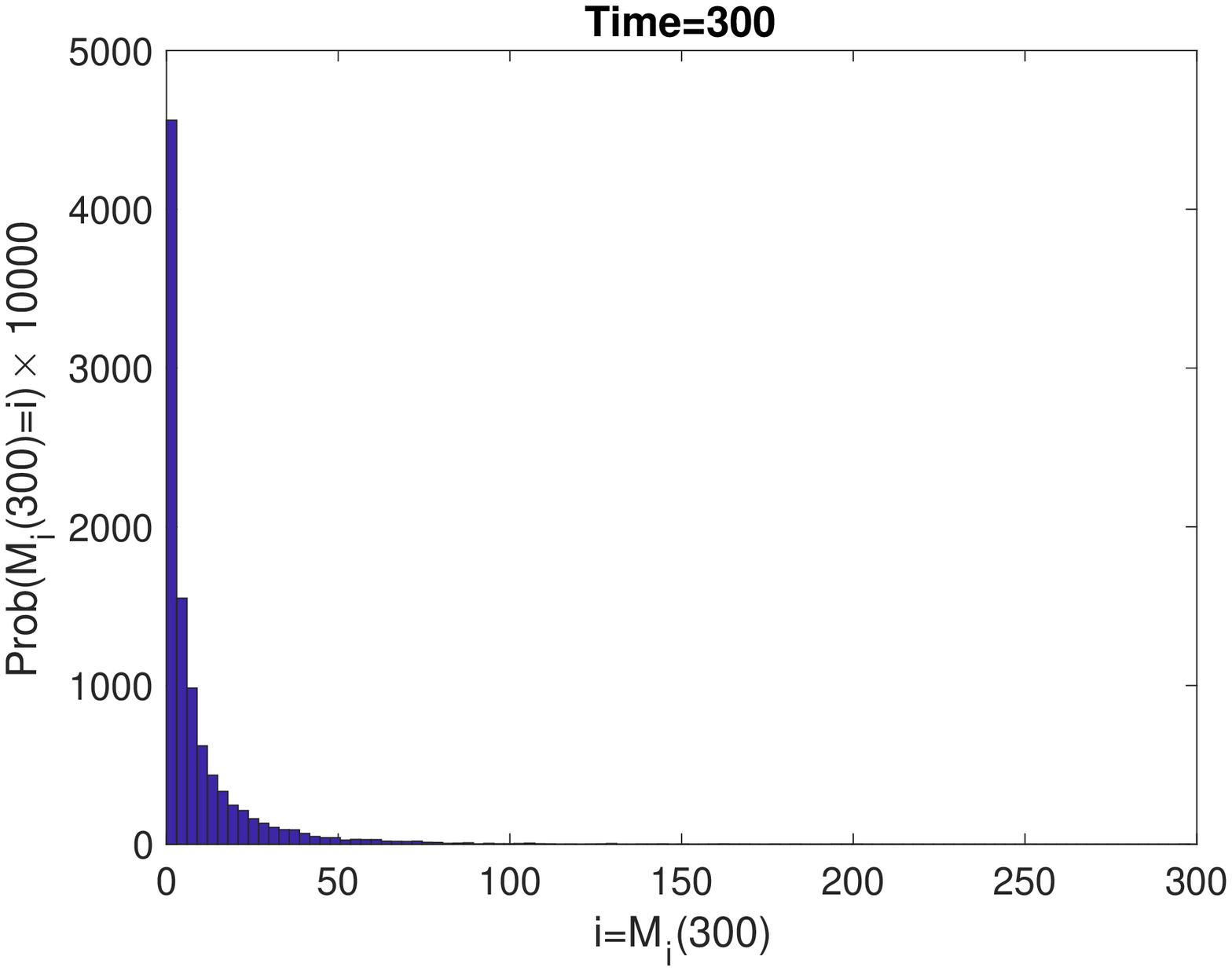}
             \label{fig_sd3d}
      \caption{$\mu=9.2332$, $\sigma=15.8717$.}
    \end{subfigure}\\
 \begin{subfigure}{0.45\textwidth}
        \includegraphics[width=1\textwidth]{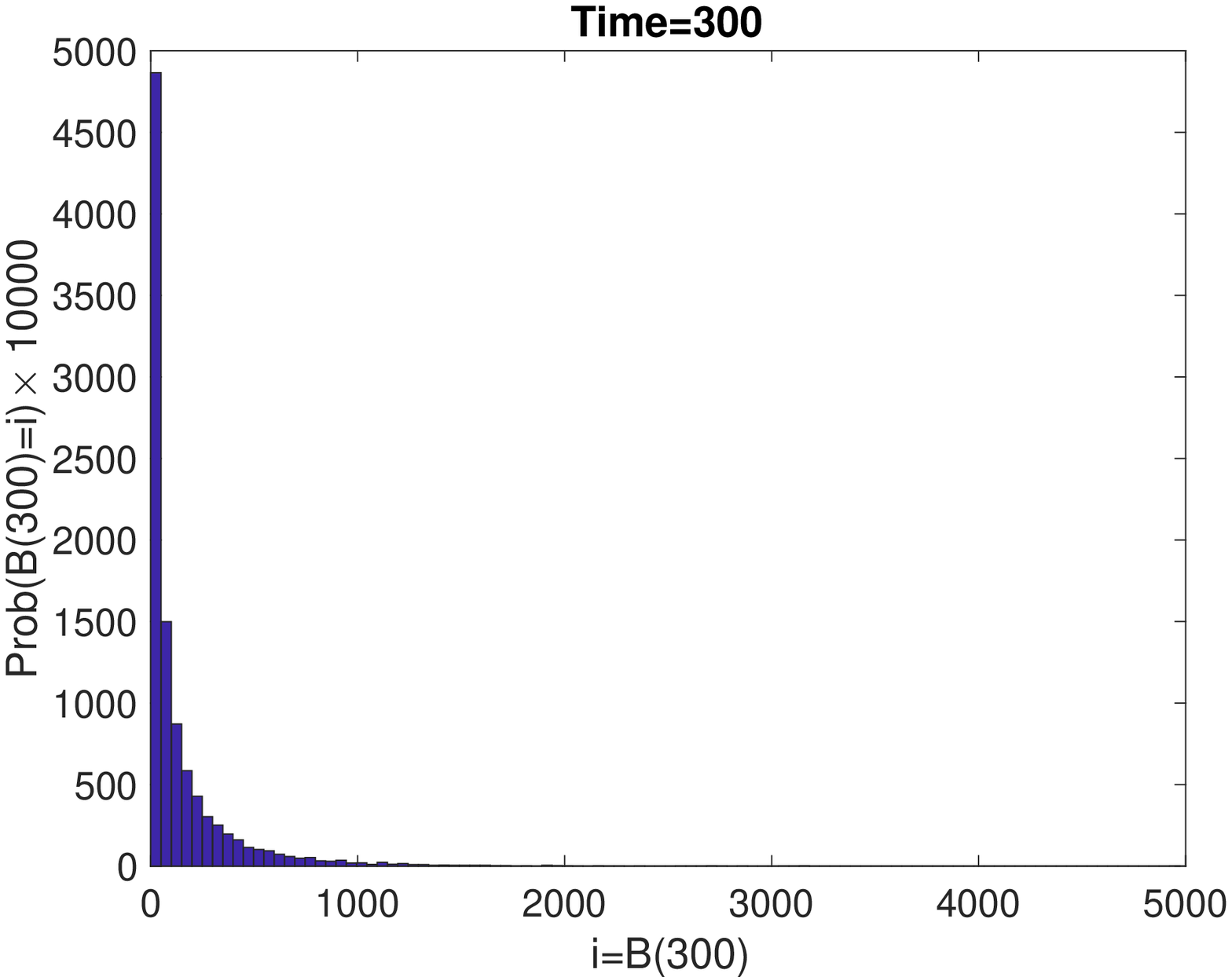}
               \label{fig_sd3e}
        \caption{$\mu=146.4230$, $\sigma=257.8087$.}
    \end{subfigure}
 \begin{subfigure}{0.45\textwidth}
    \includegraphics[width=1\textwidth]{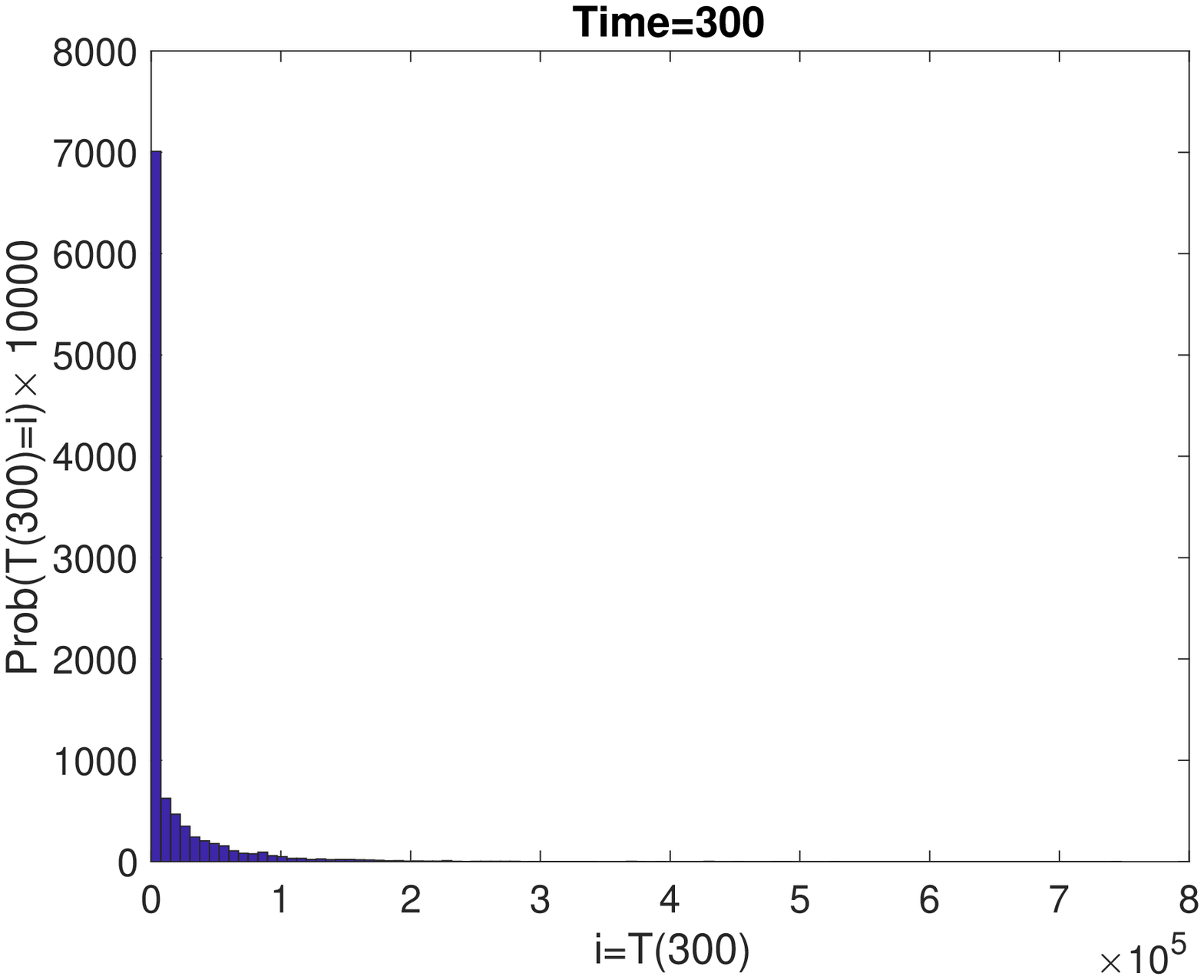}
        \caption{$\mu=1.5253\times 10^4$, $\sigma=3.6356\times 10^4$.}
       \label{fig_sd3f}
    \end{subfigure}
    \caption{Latent infected TB case at $\delta=0.27$. Mean: $\mu$, standard deviation: $\sigma$. (a): Expectations and standard deviations of 10,000 sample paths of the SDE model \eqref{sde1}.  (b) Four stochastic realizations are in colored curves. The corresponding ODE solution is in the dotted black curve. (c)-(f): Approximate stationary distributions for $M_i$, $M_u$, $B$, and $T$ at $t=300$.  $100$ bins are used for all four histograms. Initial conditions are $M_u(0)=4.99 \times 10^5$, $M_i(0)=4$, $B(0)=4$, and $T(0)=75$. Note that the other parameter values are taken from Table \ref{tab1}.}
    \label{fig_sde3}
\end{figure}

\begin{figure}
    \centering
    \begin{subfigure}{0.45\textwidth}
    \includegraphics[width=1\textwidth]{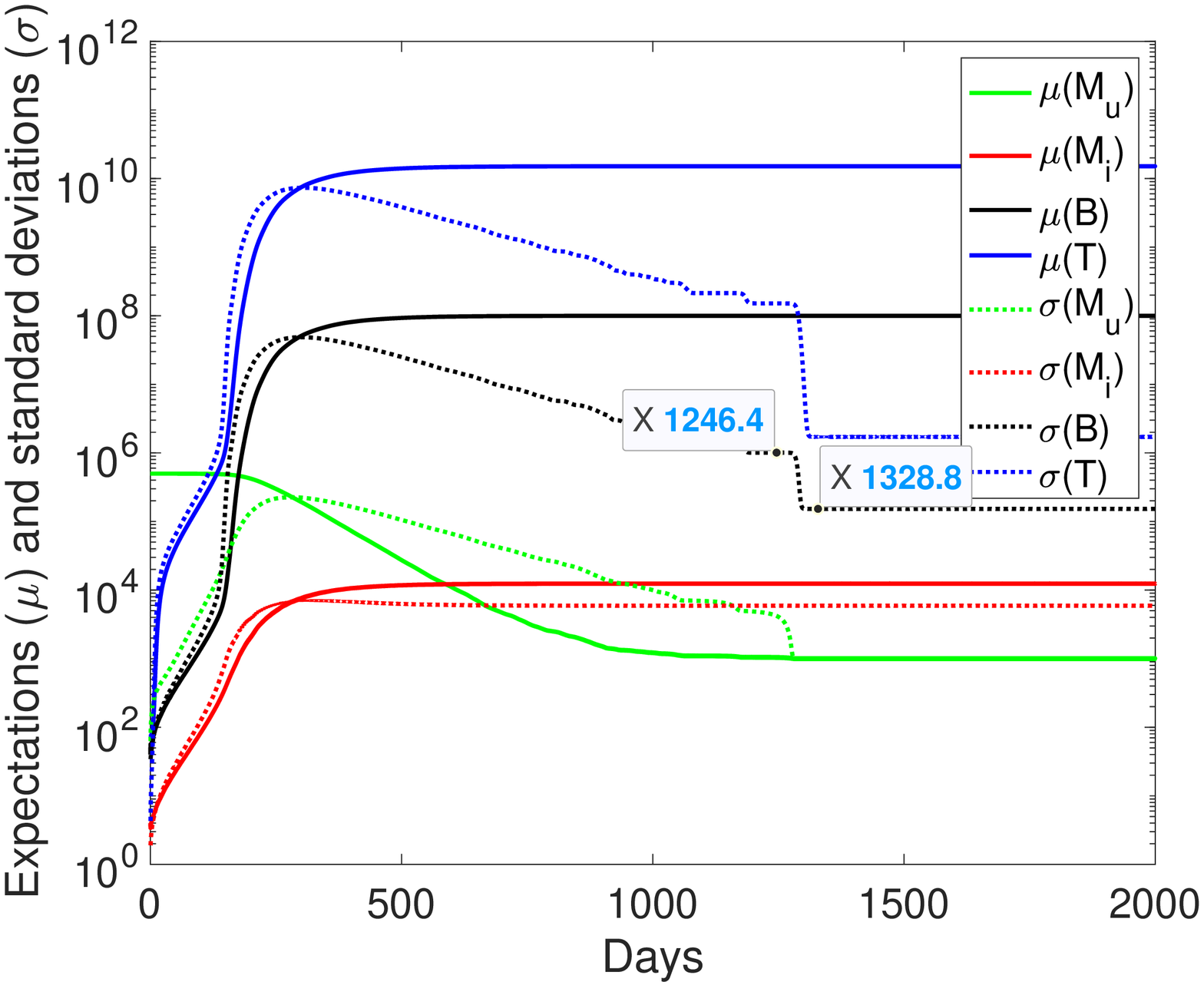}
    \caption{Mean and Std over time.}

    \end{subfigure}
     \begin{subfigure}{0.45\textwidth}
    \includegraphics[width=1\textwidth]{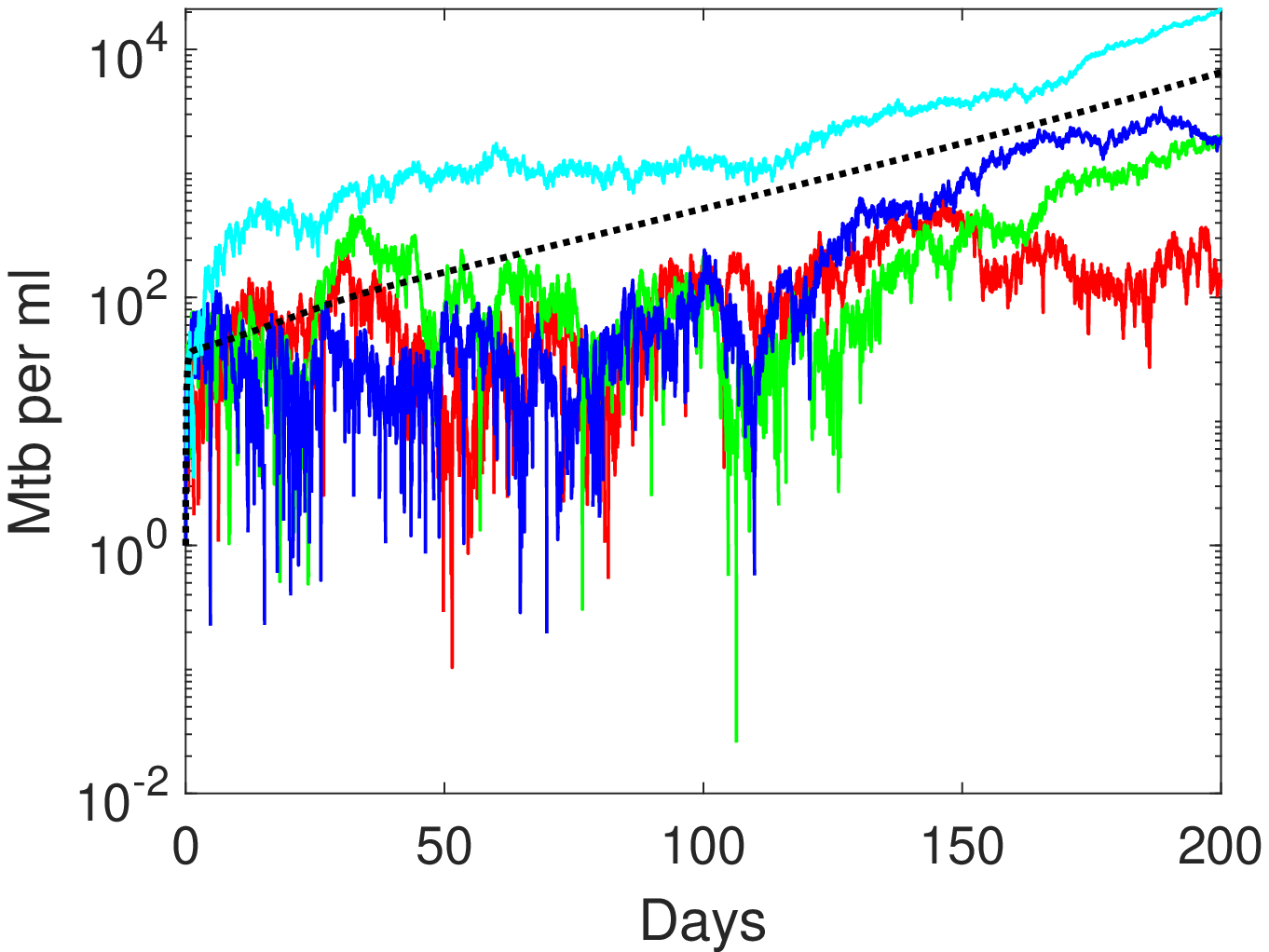}
    \caption{Colored sample paths, black ODE soln.}
    \end{subfigure}\\
     \begin{subfigure}{0.45\textwidth}
        \includegraphics[width=1\textwidth]{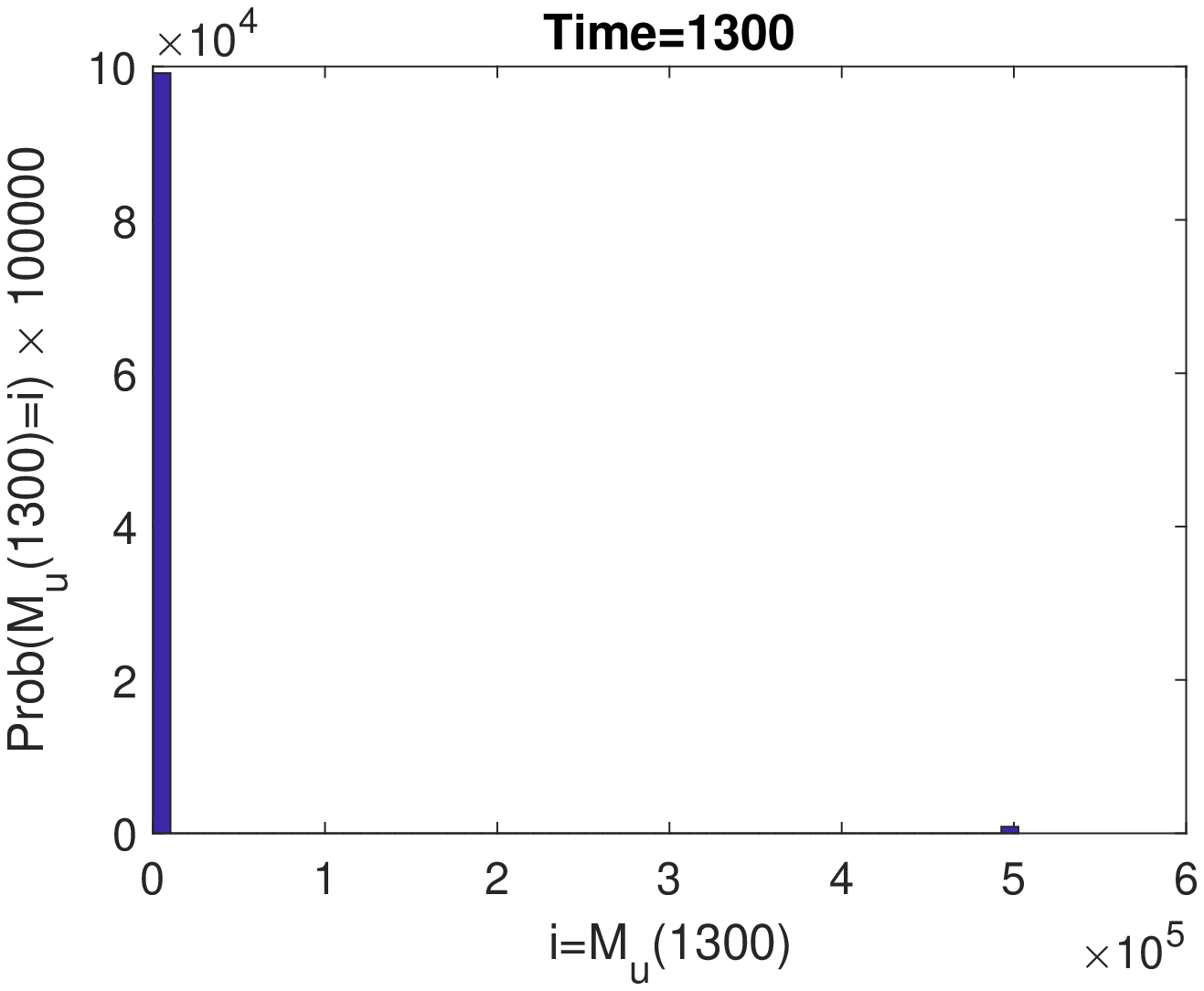}
        \caption{$\mu=4.4955\times 10^{3}$, $\sigma=4.5856\times 10^4$.}
    \end{subfigure}
 \begin{subfigure}{0.45\textwidth}
    \includegraphics[width=1\textwidth]{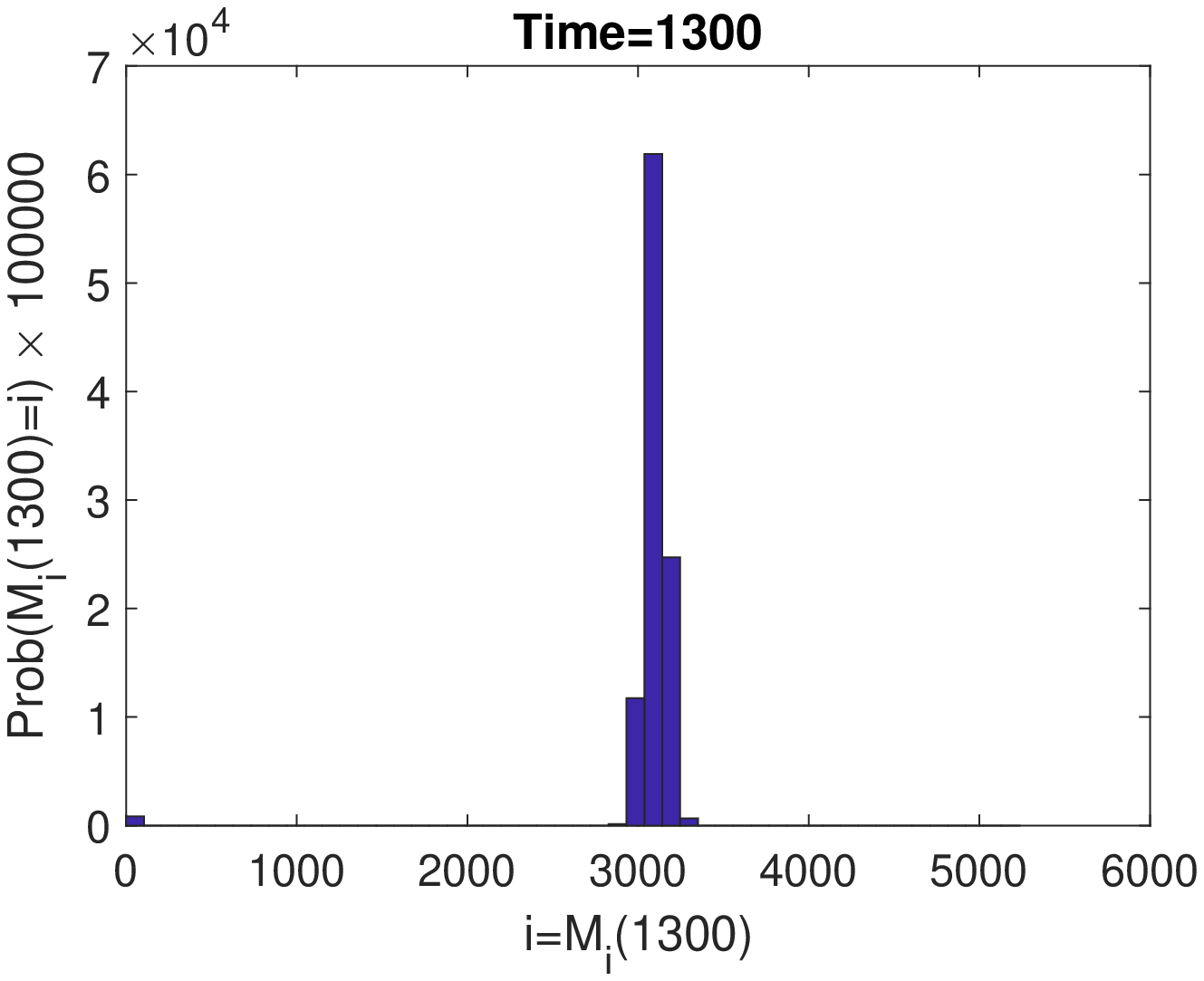}
      \caption{$\mu=3.0777\times 10^{3}$, $\sigma=290.3366$.}
    \end{subfigure}\\
 \begin{subfigure}{0.45\textwidth}
        \includegraphics[width=1\textwidth]{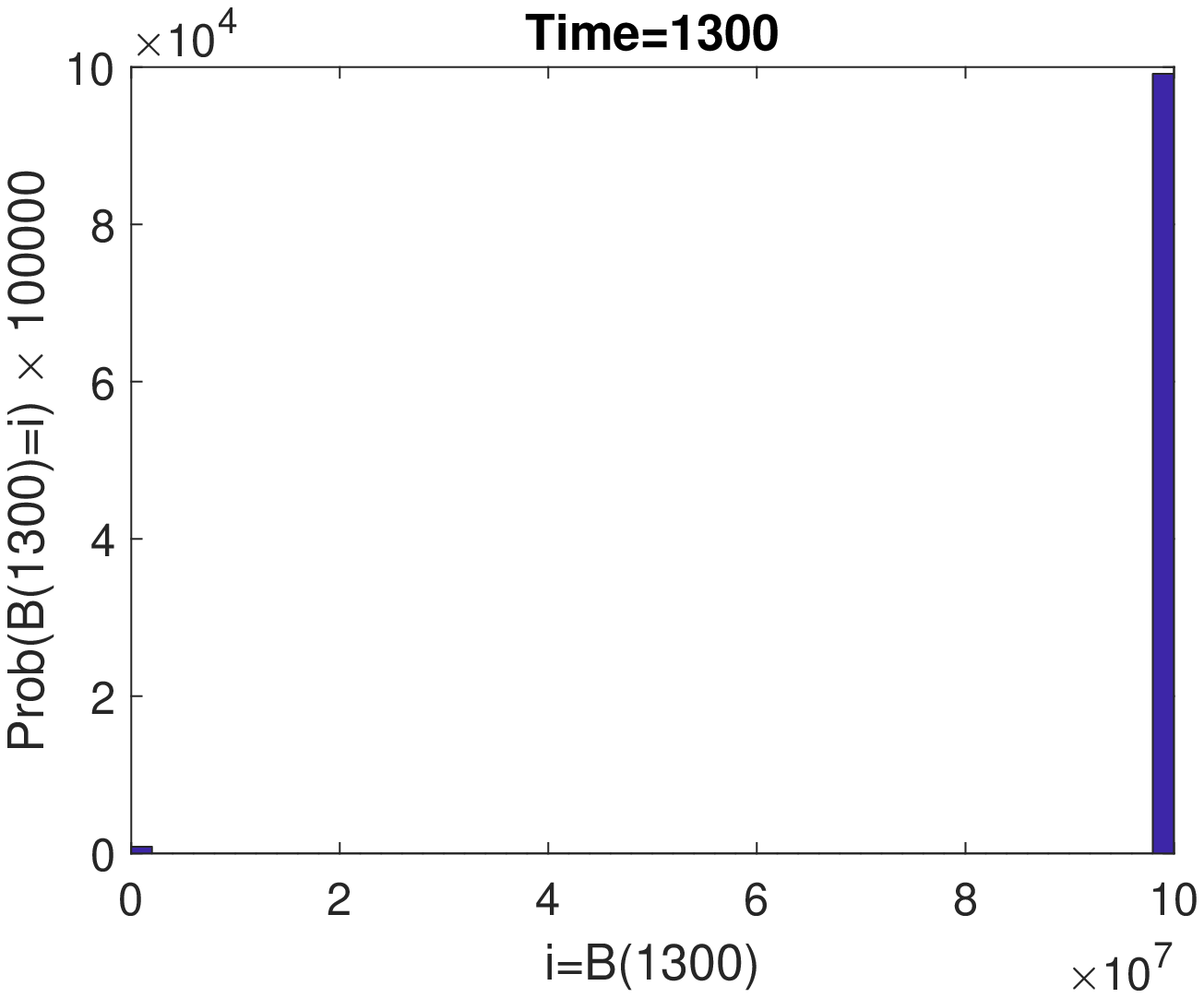}
        \caption{$\mu=9.9108\times 10^7 $, $\sigma=9.1779\times 10^6$.}
    \end{subfigure}
 \begin{subfigure}{0.45\textwidth}
    \includegraphics[width=1\textwidth]{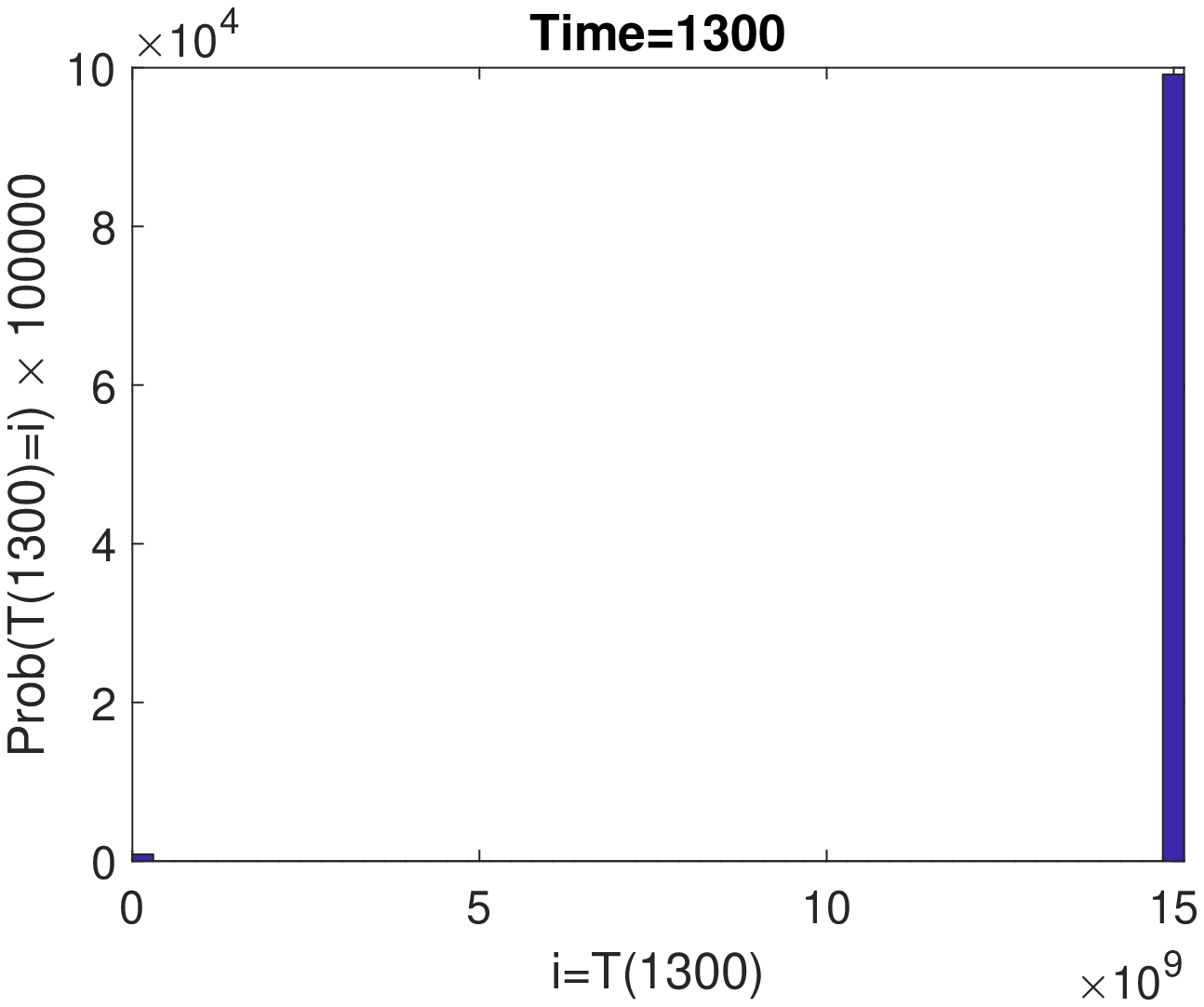}
        \caption{$\mu=1.5016\times 10^{10}$, $\sigma=1.3908\times 10^9$.}
    \end{subfigure}
    \caption{Active disease case at $\delta=0.35$. Mean: $\mu$, standard deviation: $\sigma$. (a): Expectations and standard deviations of 10,000 sample paths of the SDE model \eqref{sde1}.  (b) Four stochastic realizations are in colored curves. The corresponding ODE solution is in the dotted black curve. (c)-(f): Approximate stationary distributions for $M_i$, $M_u$, $B$, and $T$ at $t=1300$. $100$ bins are used for all four histograms. Initial conditions are $M_u(0)=4.99 \times 10^5$, $M_i(0)=4$, $B(0)=4$, and $T(0)=75$. Note that the other parameter values are taken from Table \ref{tab1}.}
    \label{fig_sde5}
\end{figure}

\subsection{Various Disease Outcomes with Demographic Variations in Different Parameter Regions}
  In Region 2, the disease stays only if a very high dose  of bacterium is inhaled, which is shown in Figure \ref{fig_bif_sim_B} (b) and (c) for $\delta=0.2 \in [\delta_h,\,\delta_l]$. 
The separatrix of  different disease outcomes is the unstable equilibrium shown in the dashed blue curves in Figure \ref{fig_bif_sim_B} (a). These dashed blue curves show Mtb concentration increases with the decrease of the $\delta$ value. This indicates that a reduction of the bacterial proliferation rate facilitates the disease progression.
For comparison, simulations for the SDE model \eqref{sde1} with demographic variations demonstrate disease clearance and active disease, as well, which are shown in Figures \ref{fig_sde1} and \ref{fig_sde2}. The initial values also determine different disease outcomes. Here, we set the parameter value $\delta=0.2 \in [\delta_h,\,\delta_l]$ and fix the other parameter values as in Table \ref{tab1}.
Based on 10,000 stochastic realizations of model \ref{sde1},  expectations and standard deviations of $M_u$, $M_i$, $B$, and $T$ stabilizes by roughly $t=300$ days.
This indicates that solutions of the SDEs \eqref{sde1} approach an approximated distribution by time $t=300$ days.
The corresponding simulations in Figure \ref{fig_sd1a} start from $M_u(0) = 5\times 10^5$, $M_i(0) = 1$, $B(0) = 10$, $T(0) =1\times 10^3$.
Four stochastic realizations along with one ODE solution are plotted in Figure \ref{fig_sd1b}.

At time $t=300$ days,  histograms of approximated stationary distributions for $M_u$, $M_i$, $B$, and $T$ are plotted in Figures \ref{fig_sd1c}-\ref{fig_sd1f}.
The histogram of the $M_u$ population in Figure \ref{fig_sd1c} fits well with the corresponding normal distribution in the black curve. Moreover, $\mathbb{E}(M_u)=4.9951 \times 10^5 \approx \Bar{M_{u0}}$ in \eqref{triv_eqili} implies that the macrophage concentration is at the uninfected level.
The histogram of $M_i$ and $B$ are heavily skewed to the right with a peak at zero, which indicates the disease clearance.
The histogram of $T$ peaks around $T\approx 20=\Bar{T}_0$ in \eqref{triv_eqili} with a long tail. This indicates the T-cell concentration peaks at the uninfected level, with a large variation.

For another set of initial conditions, $M_u(0) = 2\times 10^4$, $M_i(0) = 1\times 10^3$, $B(0) = 1\times 10^5$, $T(0) =1\times 10^3$, expectations and standard deviations of $M_u$, $M_i$, $B$, and $T$, based on 10,000 stochastic realizations, stabilize roughly by $t=100$ days,  shown in Figure \ref{fig_sd2a}.
All four stochastic realizations grow along with the ODE solution, shown in Figure \ref{fig_sd2b}.
To end at time $t=100$ days, histograms of approximated stationary distributions for $M_u$, $M_i$, $B$, and $T$ fit well with the corresponding normal distributions shown in Figures \ref{fig_sd2c}-\ref{fig_sd2f}.
Moreover, $\mu(M_u)\approx \bar{M_u}(\delta=0.2)$, $\mu(M_i)\approx \bar{M_i}(\delta=0.2)$, $\mu(B)\approx \bar{B}(\delta=0.2)$,  and $\mu(T)\approx \bar{T}(\delta=0.2)$. This indicates an active disease case.

In Region 3, the disease outcome depends on the inhaled dose of bacteria.
The disease clearance is shown in the dash-dot curve in  Figure \ref{fig_bif_sim_B} (b).
The disease progression to LTBI or active TB is shown in the dashed and dotted curves in  Figure \ref{fig_bif_sim_B} (c).
Simulation results for SDE model \eqref{sde1} in Region 3 show similar predictions for the disease clearance and active disease cases in Region 2. We, therefore, omit these two cases in Region 3 and focus on the
additional state with low bacterial level generated from the saddle-node bifurcation LP$_1$.
We take parameter value $\delta=0.27 \in [\delta_l,\,\delta_m]$ in Region 3. In the ODE model \eqref{eqn1}, the stable low bacterial level equilibrium denoting LTBI  is $(\bar{M}_u,\bar{M}_i,\bar{B},\bar{T})=(499519.6475,3.3532,48.0814,74.9548)$.
Starting at $M_u(0) = 4.99\times 10^5$, $M_i(0) = 4$, $B(0) = 4$, $T(0) = 75$, based on 10,000 stochastic realizations, expectations and standard deviations of $M_u$, $M_i$, $B$, and $T$ stabilize roughly by $t=300$ days, as shown in Figure \ref{fig_sd3a}.
Four sample paths and one corresponding ODE solution are graphed in Figure \ref{fig_sd3b}.
At time $t=100$ days, approximated stationary distributions for $M_u$, $M_i$, $B$, and $T$ have medians as
Median$(M_u)=4.9893\times10^5$, Median$(M_i)=3.6589$, Median$(B)=52.7568$, and Median$(T)=417.5100$. 
Even though the histograms for $M_i$, $B$, and $T$ have skewed shapes, the median values for $M_u$, $M_i$, and $B$ are close to the stable latent infected equilibrium $(\bar{M}_u,\bar{M}_i,\bar{B},\bar{T})=(499519.6475,3.3532,48.0814,74.9548)$.
Notice that the Median$(T)=417.5100$ is much larger than $\bar{T}$. 
This indicates that the demographic variations significantly affect   T-cell population.

In Region 4, the disease will progress to active TB shown in the dash-dot curve in Figure \ref{fig_bif_sim_B} (b) by the ODE model prediction. The nontrivial equilibrium of the ODE model \eqref{eqn1} at $\delta=0.35$ is $(\bar{M}_u,\bar{M}_i,\bar{B},\bar{T})=(249.9805,3104.0391, 9.9957 \times 10^7, 1.5145 \times10^{10})$. Starting at $M_u(0) = 4.99\times 10^5$, $M_i(0) = 4$, $B(0) = 4$, $T(0) = 75$, based on 100,000 stochastic realizations, expectations and standard deviations of $M_u$, $M_i$, $B$, and $T$ stabilize  roughly by $t=1300$ days, as shown in Figure \ref{fig_sde5} (a). 
Four sample paths blow up along with the ODE solution shown in Figure \ref{fig_sde5} (b).
At time $t=1300$ days, the means of approximated stationary distributions for  $M_i$, $B$, and $T$ are close to  $\bar{M}_i$, $\bar{B}$, and $\bar{T}$.
The mean of $M_u$ is far away from $\bar{M}_u$, but the standard deviation of the approximated stationary distribution is large as well. 
It implies that small exposure to Mtb bacterial is more likely to develop to active disease.
\begin{figure}{r}
  \begin{center}
  (a) \qquad\qquad\qquad\qquad\qquad\qquad\qquad\qquad(b)\\
    \includegraphics[width=0.48\textwidth]{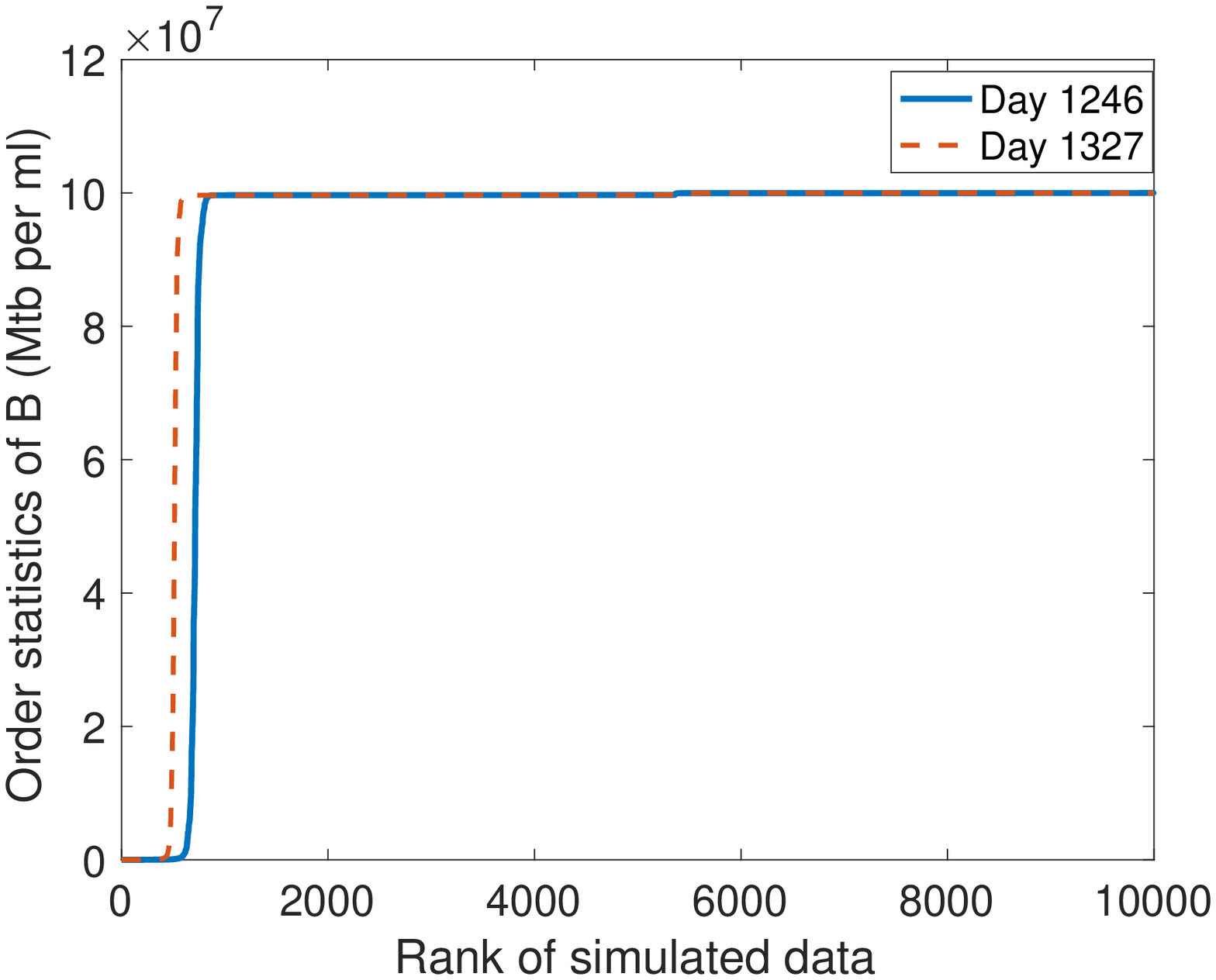}
     \includegraphics[width=0.48\textwidth]{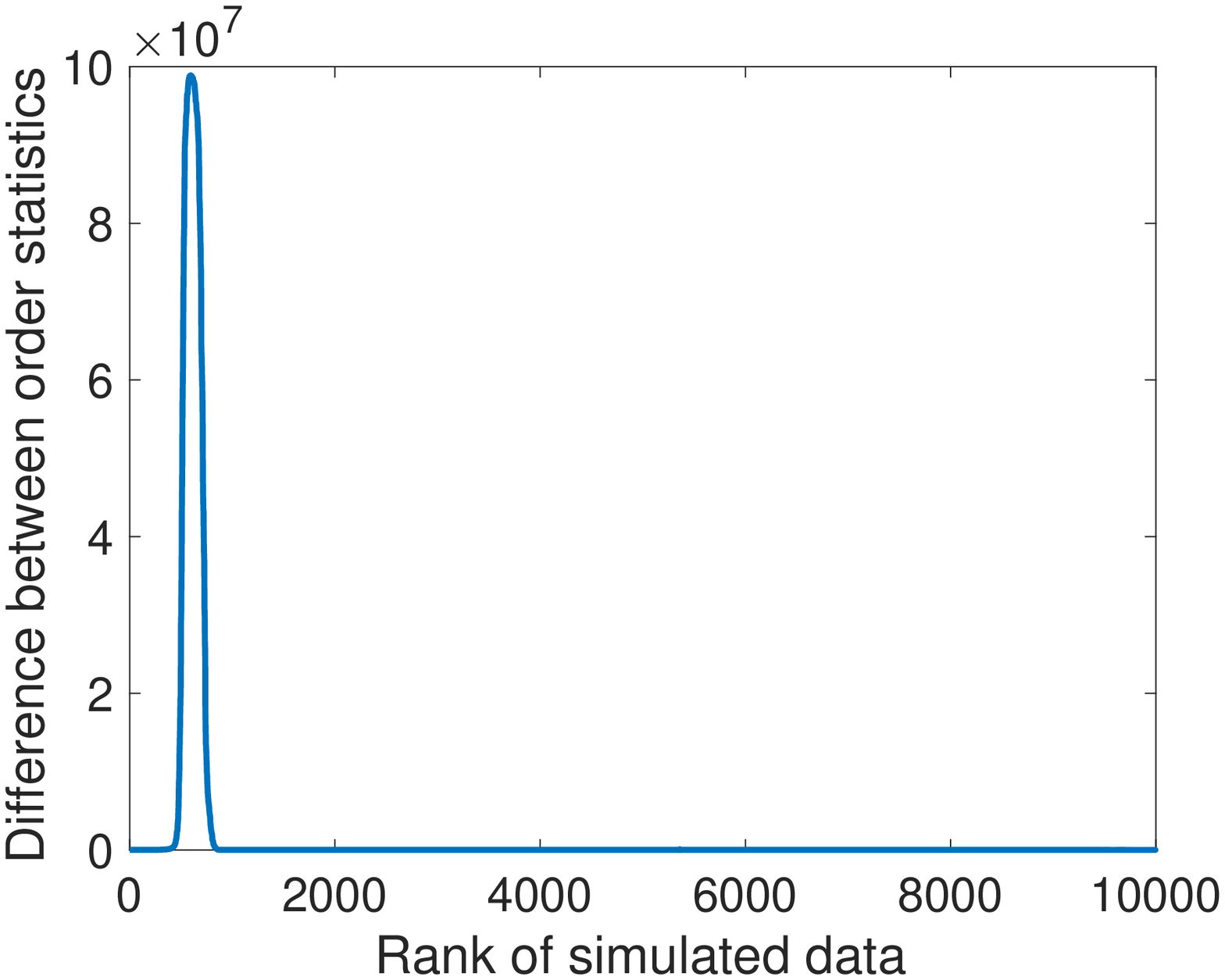}
  \end{center}
  \caption{The horizontal axis denotes the rank of the observations (Mtb concentration $B$) at day $ 1246$ and day $ 1327$. The vertical axis in (a) represents the value of the observation ($B$ value) at each rank. The vertical axis in (b) represents the difference in $B$ value between day $ 1246$ and day $ 1327$.}
  \label{fig_order}
\end{figure}

The decreasing parts of the standard deviation curves in Figure \ref{fig_sde5} (a) indicate that more cell populations converge to their corresponding mean values. 
Moreover, there are sharp drops at the time roughly between $t\in[1246,\,1327]$ in the standard deviation curves of the same figure.
The order statistic Figure \ref{fig_order} suggests that the Mtb concentration, for example, jumps sharply from a low level to a much higher level.
This  causes the sharp drop in the standard deviation during the time interval $[1246,\,1327]$. 
In other words, there are more large values in  day $1327$ than in  day $1246$. Therefore, the standard deviation is smaller on the later day than the earlier day.
It implies that the development of active disease is more likely to happen after  day $250$ and before  day $1328$.
Interestingly, there are some extremely short bars in the histograms in Figures \ref{fig_sde5} (c) - (f) that show that there  still exists a possibility for disease clearance.

Here are some notes about the SDE simulations in Figures \ref{fig_sde1}, \ref{fig_sde2}, \ref{fig_sde3}, and \ref{fig_sde5}:
\begin{enumerate}
    \item  For each figure, we simulate SDE model \eqref{sde1}   $10,000$ times. With the obtained $10,000$ sample paths, we calculate their expectations and standard deviations at each time step, then demonstrate the trends of convergence in  time-series plots shown in the corresponding subfigures (a). Next, we choose simulation times at which the corresponding expectations and standard deviations converge to their stabilized values, that is  $t=300$, $t=100$, $t=300$, $t=1300$ for Figures \ref{fig_sde1}, \ref{fig_sde2}, \ref{fig_sde3}, and \ref{fig_sde5}, respectively.
  \item The bin size for the histogram is chosen as the square root of the number of data points, which is $10,000$ here. The data points are SDE solutions at the simulation end time. Therefore, we choose $100$ as the bin size for all four SDE simulations in Figures \ref{fig_sde1}, \ref{fig_sde2}, \ref{fig_sde3}, and \ref{fig_sde5}. 
  \item The baseline parameter values are shown in Table \ref{tab1}. The parameter ranges are provided when variations of parameters are needed.
Note that both the baseline values and the ranges of parameters are from existing modeling papers. 
  \item In Region 2, the ODE model \eqref{eqn1} contains two stable steady states. 
We demonstrate disease clearance and active disease under demographic variations in Figures \ref{fig_sde1} and \ref{fig_sde2}. 
To generate Figures \ref{fig_sde1} and \ref{fig_sde2}, we take the same set of parameter values, but different initial conditions (see the  figure captions). 
In Region 3, the ODE model \eqref{eqn1} has three stable steady states, which represent disease clearance, LTBI, and active disease.
SDE simulations for disease clearance and active disease show the similar behavior as in Figures \ref{fig_sde1} and \ref{fig_sde2}. We, therefore, omit the graphs for these two cases, but only demonstrate the LTBI case in Figure \ref{fig_sde3}. 
In Region 4, the ODE model \eqref{eqn1} has only one steady state, which represents active disease. Interestingly, the SDE simulation predicts a possibility for the occurrence of disease clearance (see the tiny bars in histograms in Figure \ref{fig_sde5}). 
Moreover, the curves representing standard derivations in Figure \ref{fig_sde5} have decreasing parts with sharp drops.
  These drops imply that the development of active disease is more likely to happen from the end of the first year to the third year after the initial exposure.
\end{enumerate}

\section{Influences of Bacterial and Host Factors on Disease Progression Speed and Combination Therapy Outcomes}
  
In this section, we consider the impact of both the bacterial factor (that is the infected macrophage proliferation rate $\delta$) and host factors (including the infected macrophage loss rate $b$,  the cell-mediated immunity rate $\gamma$, and macrophages' bacteria killing rate $\eta$) in evaluating the pathogen-directed (antibiotic) therapy with adjunctive host-directed therapies. 
Previous work (\cite{zhang2021investigation}) sugests that these four parameters have the most statistically significant impacts on model behaviors.
We first study the disease dynamics influenced  by both the bacterial proliferation rate modulated by pathogen-directed therapy and  the host immune components ($b$, $\gamma$, and $\eta$) altered by host-directed therapy. Our results suggest that host immune responses can both promote and hinder disease progression. Different disease outcomes depend on the relation between the numbers of Mtb engulfed by and released by macrophages.
Considering the stochastic variations from both the cell populations and the adherence to therapies, we investigate the SDE model \eqref{sde2} with both demographic and environmental variations to provide insight into the pathogen- and host-directed combination therapy developments. 
Next, we derive a quadratic relationship between the Mtb spreading speed and the bacterial proliferation rate $\delta$. Then, we further investigate the disease progression speed under the influences of both pathogen-directed therapy (such as antibiotics affecting $\delta$) and host-directed therapy (such as vitamin D affecting $b$, $\gamma$, and $\eta$).
Our results demonstrate the advantages of the combination therapy.
\begin{figure}[h!]
    \centering
    \begin{subfigure}{.3\textwidth}
        \centering
        \includegraphics[width=1\textwidth]{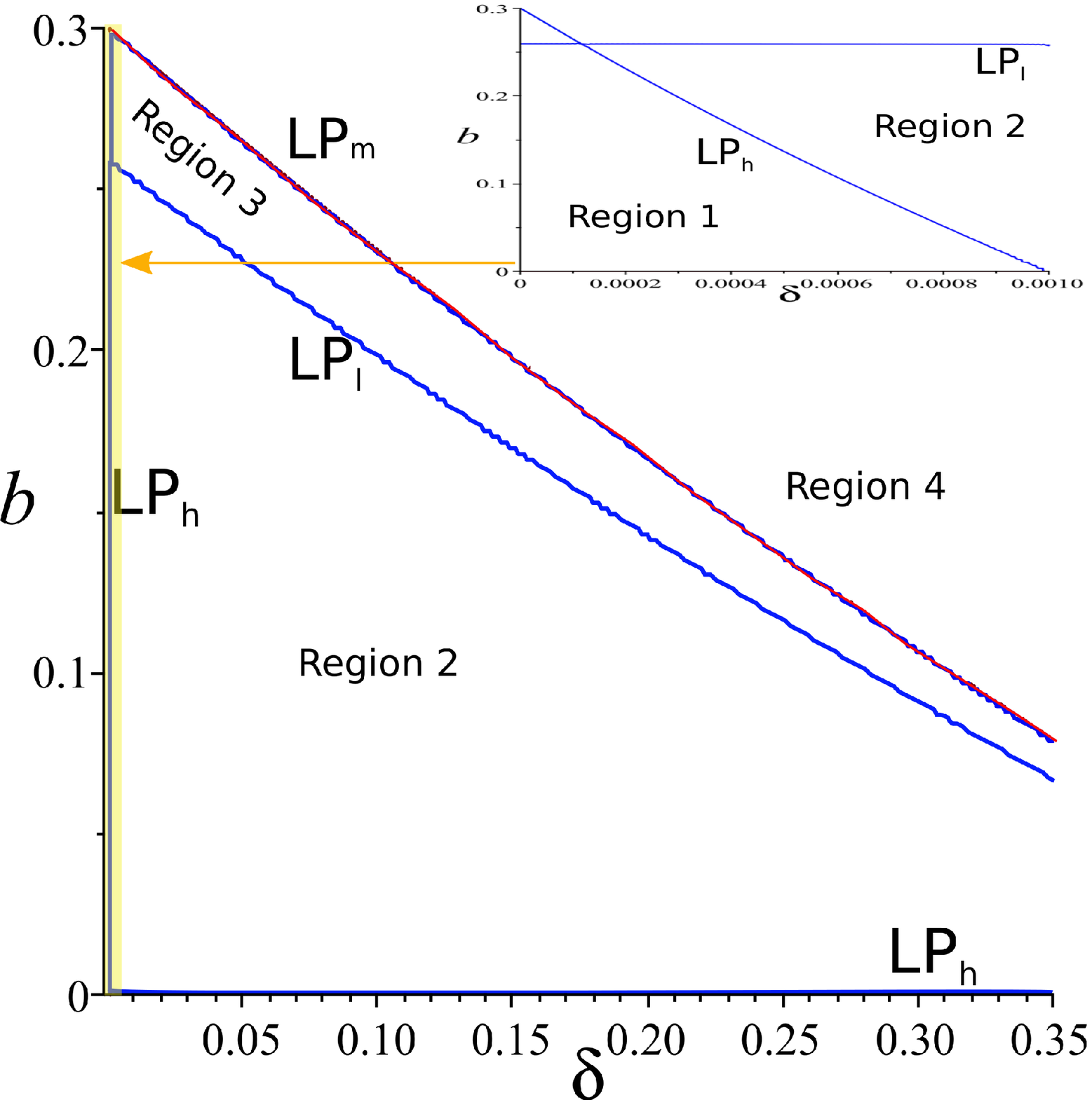}
            \caption{ $\delta$ vs $b$}
    \end{subfigure}
        \begin{subfigure}{.3\textwidth}
            \centering
     \includegraphics[width=.9\textwidth]{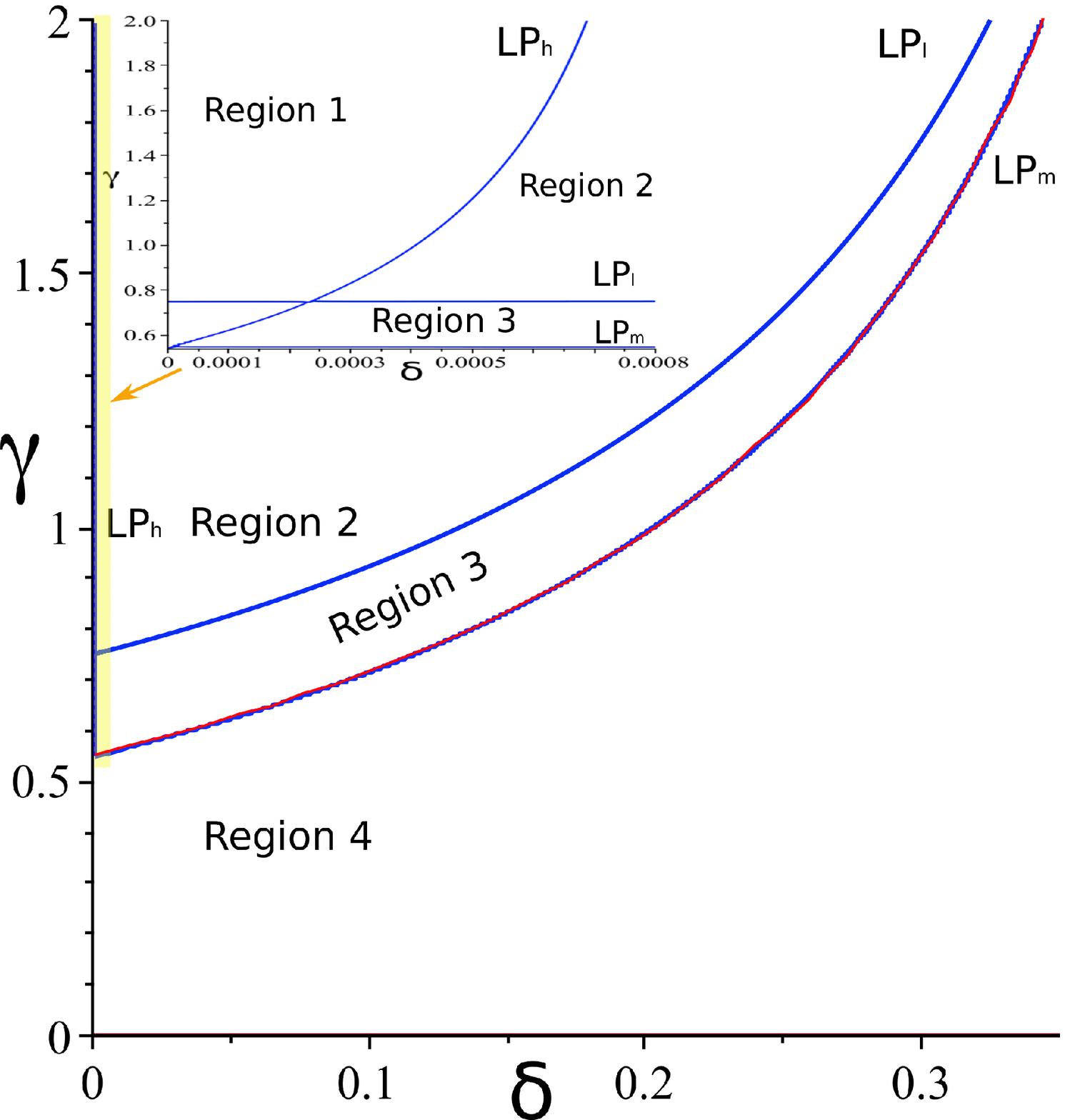}
         \caption{$\delta$ vs $\gamma$}
    \end{subfigure}
        \begin{subfigure}{.3\textwidth}
            \centering
     \includegraphics[width=.9\textwidth]{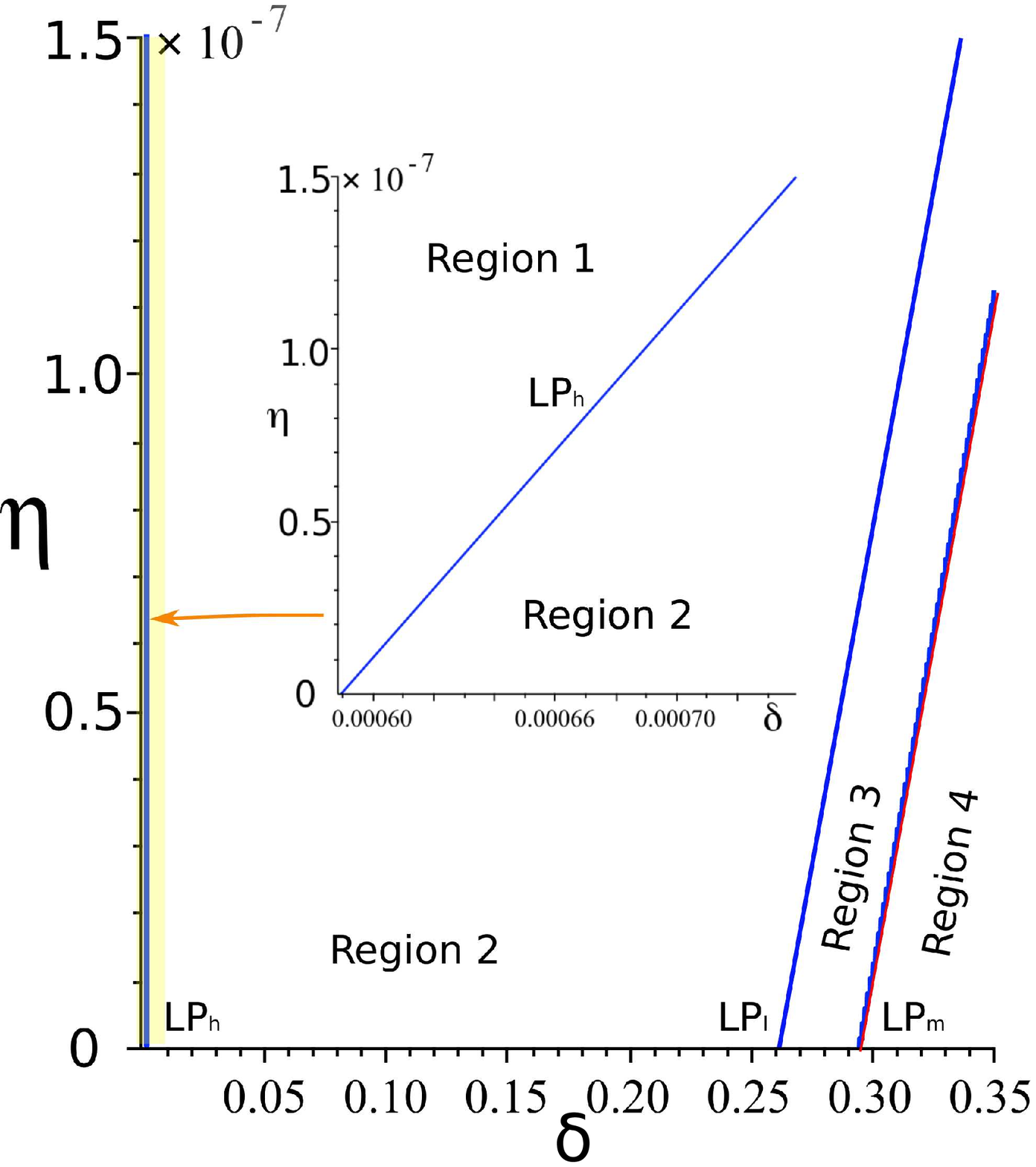}
    \caption{$\delta$ vs $\eta$}
    \end{subfigure}
    
        \begin{subfigure}{.32\textwidth}
        \centering
        \includegraphics[width=1\textwidth]{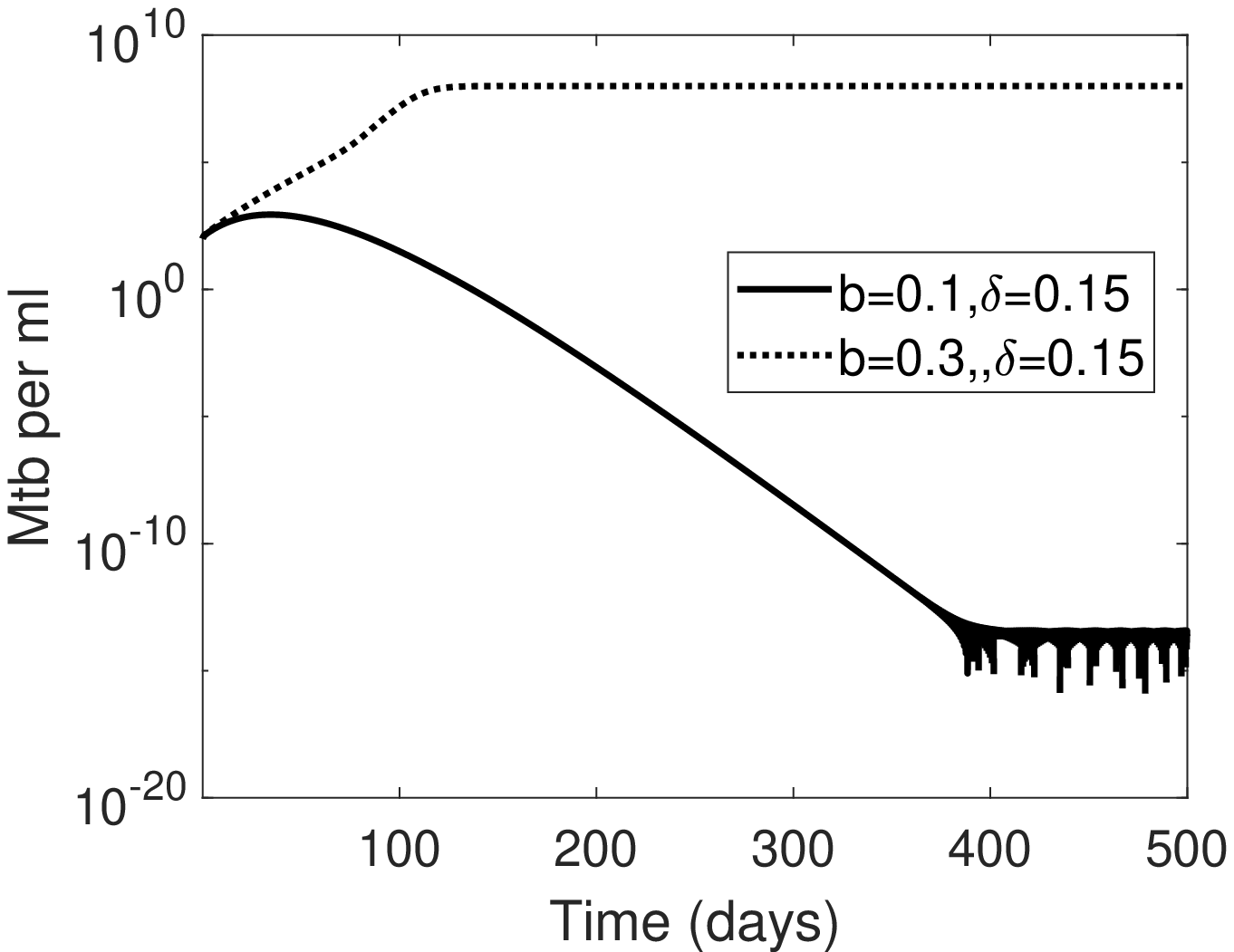}
            \caption{Boost the death of $M_i$}
    \end{subfigure}
        \begin{subfigure}{.32\textwidth}
            \centering
     \includegraphics[width=1\textwidth]{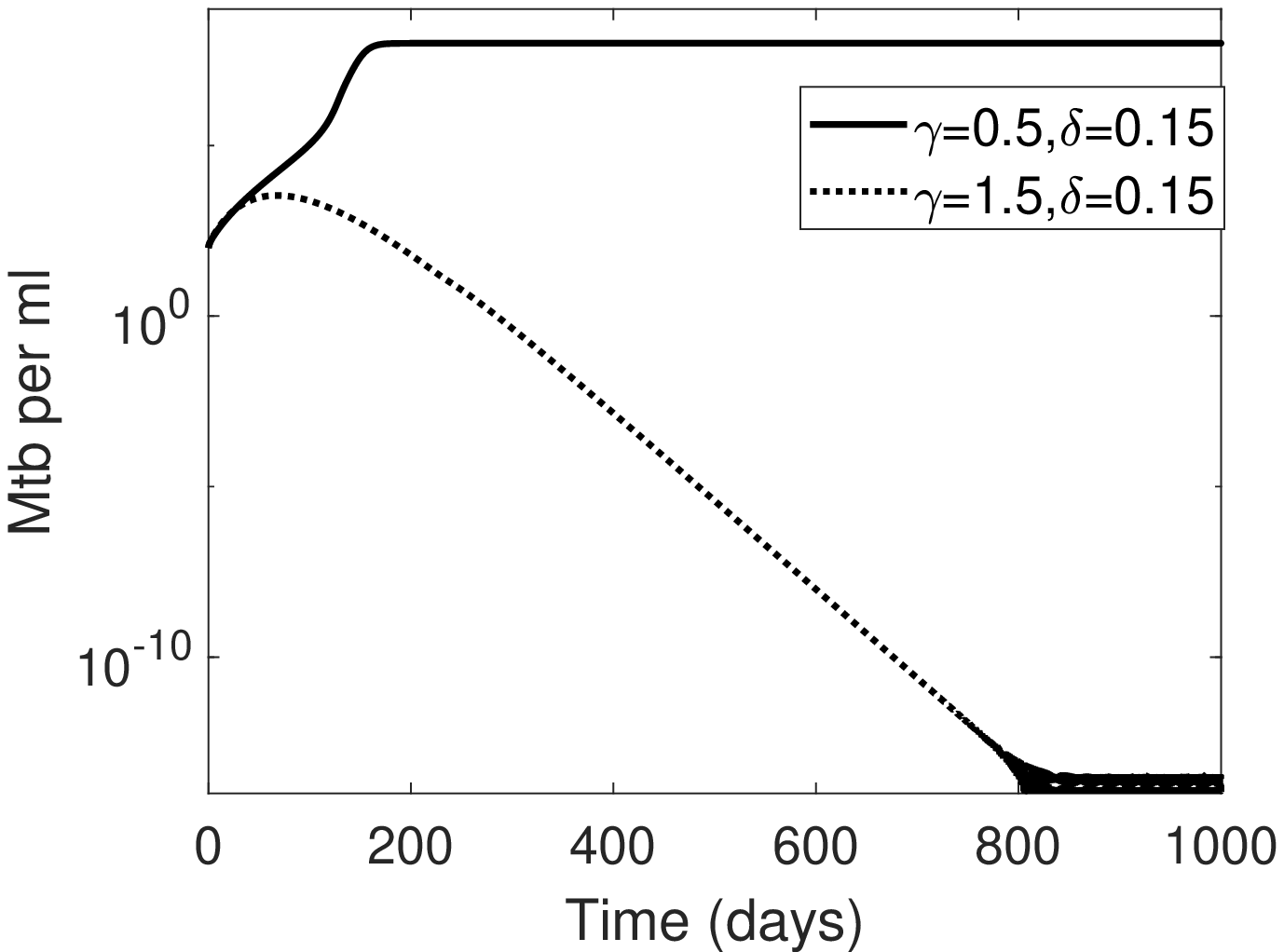}
         \caption{Boost the cell-mediated immunity}
    \end{subfigure}
        \begin{subfigure}{.32\textwidth}
            \centering
     \includegraphics[width=1\textwidth]{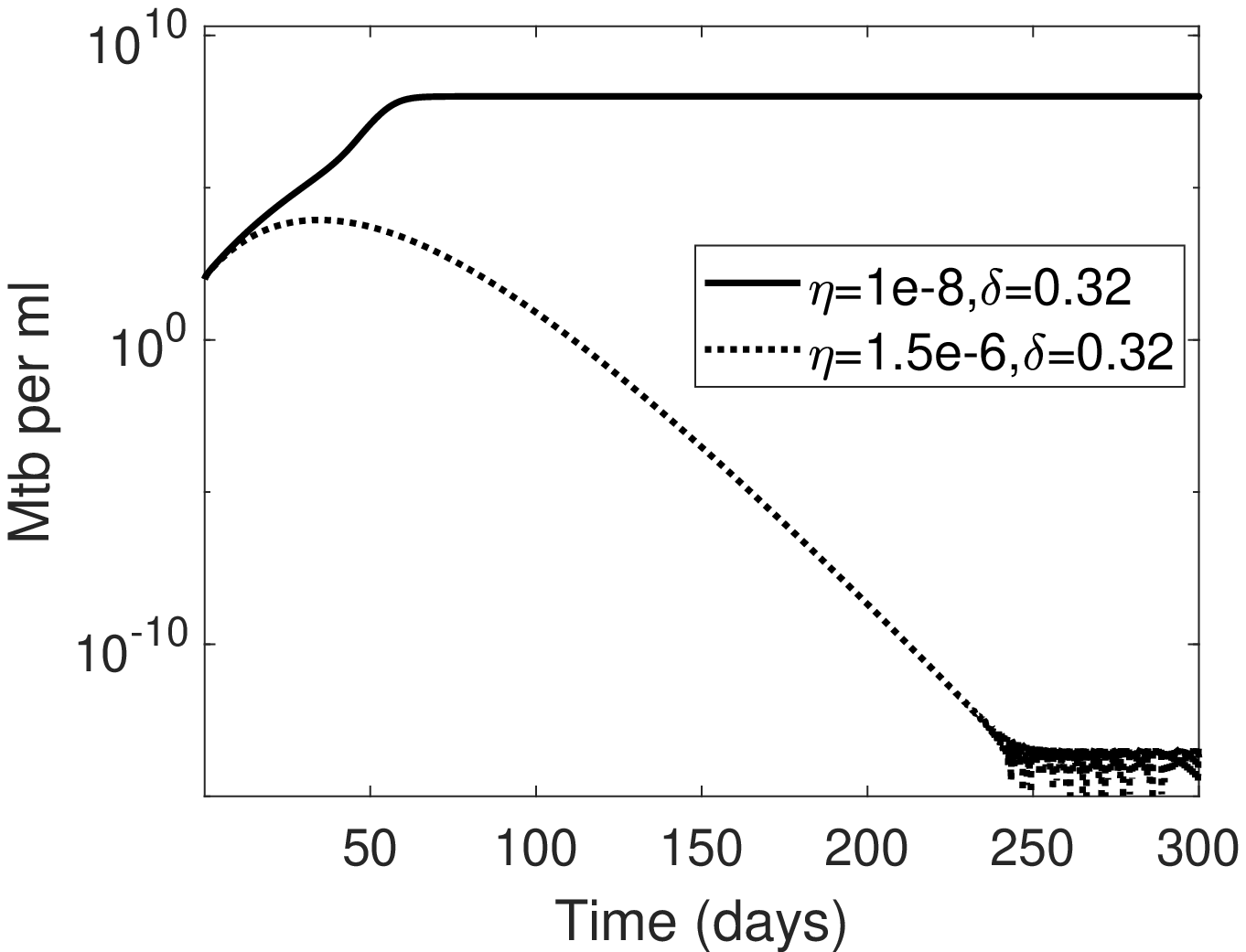}
    \caption{Boosted intracellular Mtb death}
    \end{subfigure}
    \caption{ 
    The influence of host immune response to bactericidal drug therapy.
    (a), (b), and (c) are two-parameter bifurcation diagrams.
    Blue and red curves denote saddle node (LP) and transcritical (BP) bifurcations, respectively.
    Numerical simulation (d) corresponds bifurcation diagram (a); (e) corresponds (b); and (f) corresponds (c). (d) shows that an increase in the death of infected macrophages ($b$) favors Mtb infection. While (e) and (f) shows that boosting the cell-mediated immunity ($\gamma$) and the intracellular Mtb death rate ($\eta$) benefit the host defense. Note that we take the other parameter values from Table \ref{tab1}   and initial values as $M_u(0)= 10$, $M_i(0)= 1$, $B(0)=100$, and $T(0)= 1000$.}
    \label{fig_bif_del_2d}
\end{figure}
\subsection{Bacterial Concentration Affected by Both the Bacterial Proliferation Rate and Acquired Defects}
The analysis results from the ODE model \eqref{eqn1} in Figure \ref{fig_bif_sim_B} indicate the possibility of disease clearance if the bacterial proliferation rate $\delta$ falls in Region 2, the potential to progress to LTBI or active TB if  $\delta$ is in Region 3, and a  progression to active TB if  $\delta$ is in Region 4 by the ODE model prediction.
There also exist $\delta$ values, which guarantee disease clearance. 
However, this $\delta$ value window is slim and demonstrated in yellow strips in Figure \ref{fig_bif_del_2d}.
We denote this  parameter region for disease clearance as Region 1. The four parameter regions are 
delimited by saddle node (LP) and transcritical (BP) bifurcations.
Two parameter bifurcation analyses were carried out via Matcont (\cite{dhooge2003matcont}). The corresponding bifurcation diagrams are plotted in Figure \ref{fig_bif_del_2d} and illustrate
the four regions associated with bacterial proliferation rate $\delta$ and three parameters denoting host immune responses (i.e. $b$, $\gamma$, and $\eta$).
It is shown that an increase in the loss rate of infected macrophages $b$ and a decrease of both the cell-mediated immunity rate $\gamma$ and the bacteria killing rate by macrophages $\eta$ impede disease clearance (in Region 1 and 2) and facilitate the disease development (in Region 3 and 4).
Note that parameter relations $N_1>N_3$ and $N_2>N_3$ are taken for the bifurcation analysis under the assumption that macrophages are infected with virulent Mtb and undergo necrosis (\cite{behar2011apoptosis}). 
The bifurcation parameter ranges are provided by \cite{du2017simple}:
$\delta \in (0,\,0.35)$, $b\,\in (0.05,\,0.5)$, $\gamma \,\in (0.1,\,2)$, and $\eta \in (1.25 \times 10^{-9},\,1.25 \times 10^{-7})$.
Other than the bifurcation parameters, the other parameter values are taken from Table \ref{tab1}.

\begin{figure}[t]
    \centering
    \begin{subfigure}{0.24\textwidth}
    \includegraphics[width=1\textwidth]{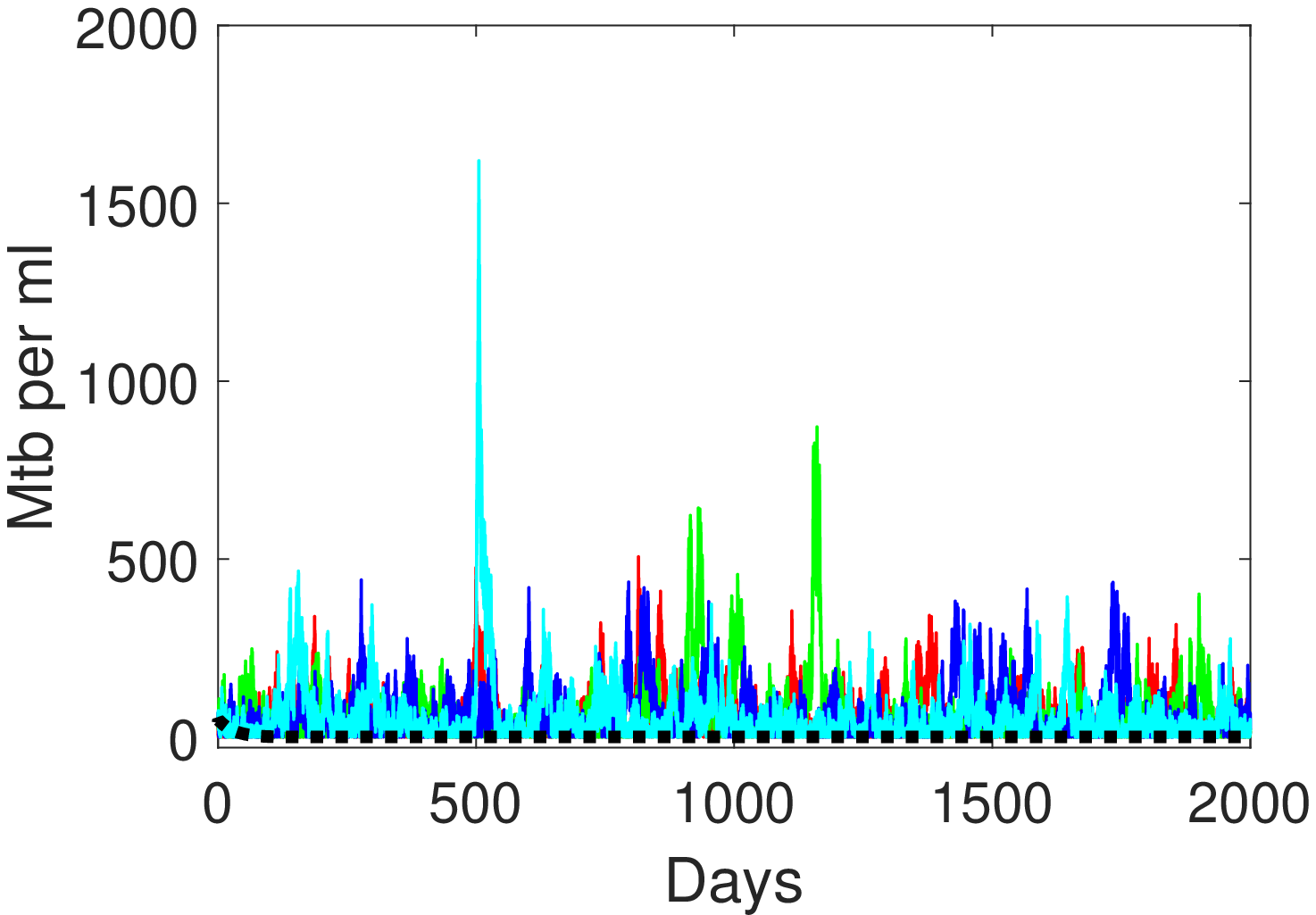}
    \caption{RGN 2, fast return}
    \end{subfigure}
     \begin{subfigure}{0.24\textwidth}
    \includegraphics[width=1\textwidth]{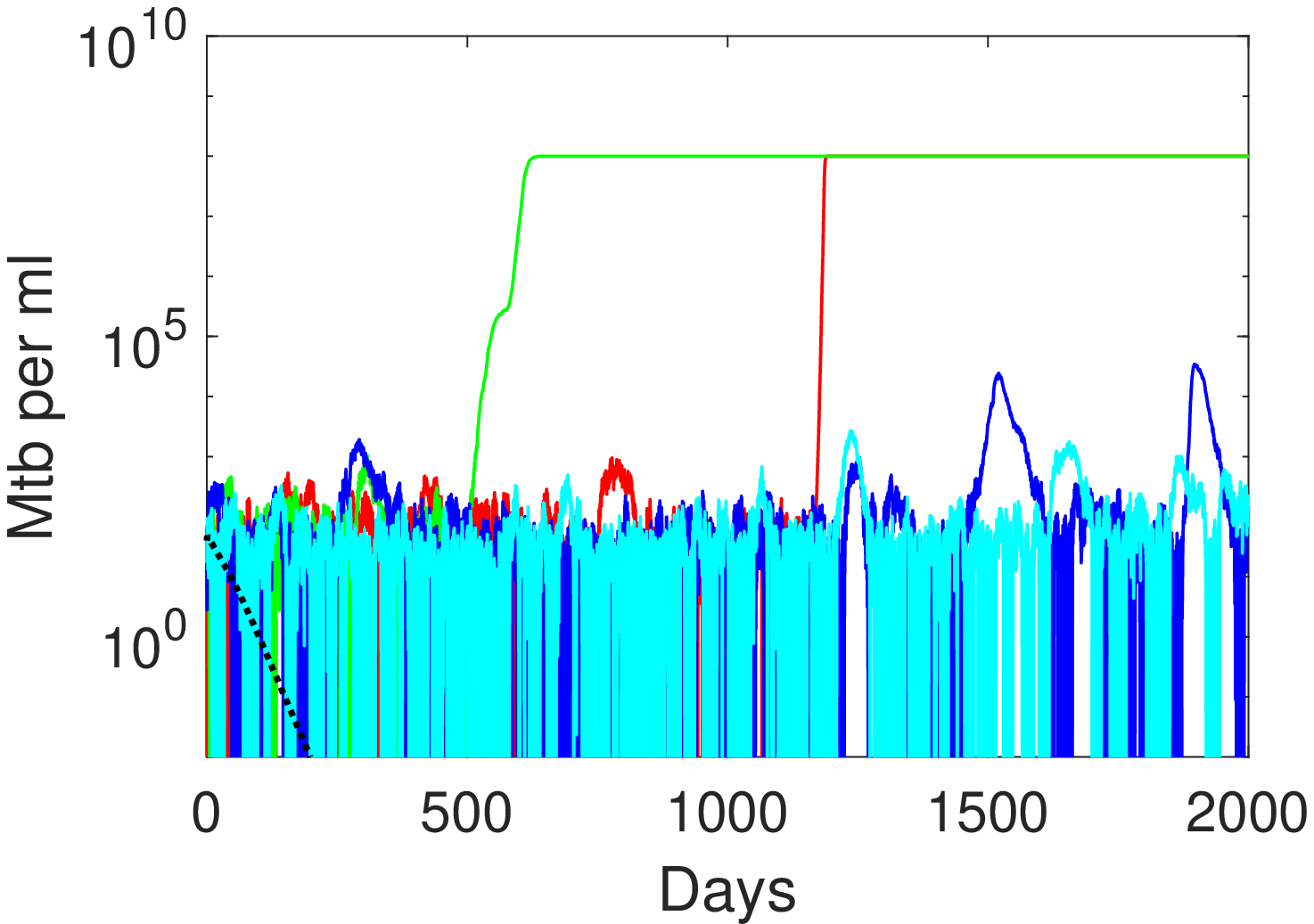}
    \caption{RGN 2, slow return}
    \end{subfigure}
        \begin{subfigure}{0.24\textwidth}
    \includegraphics[width=1\textwidth]{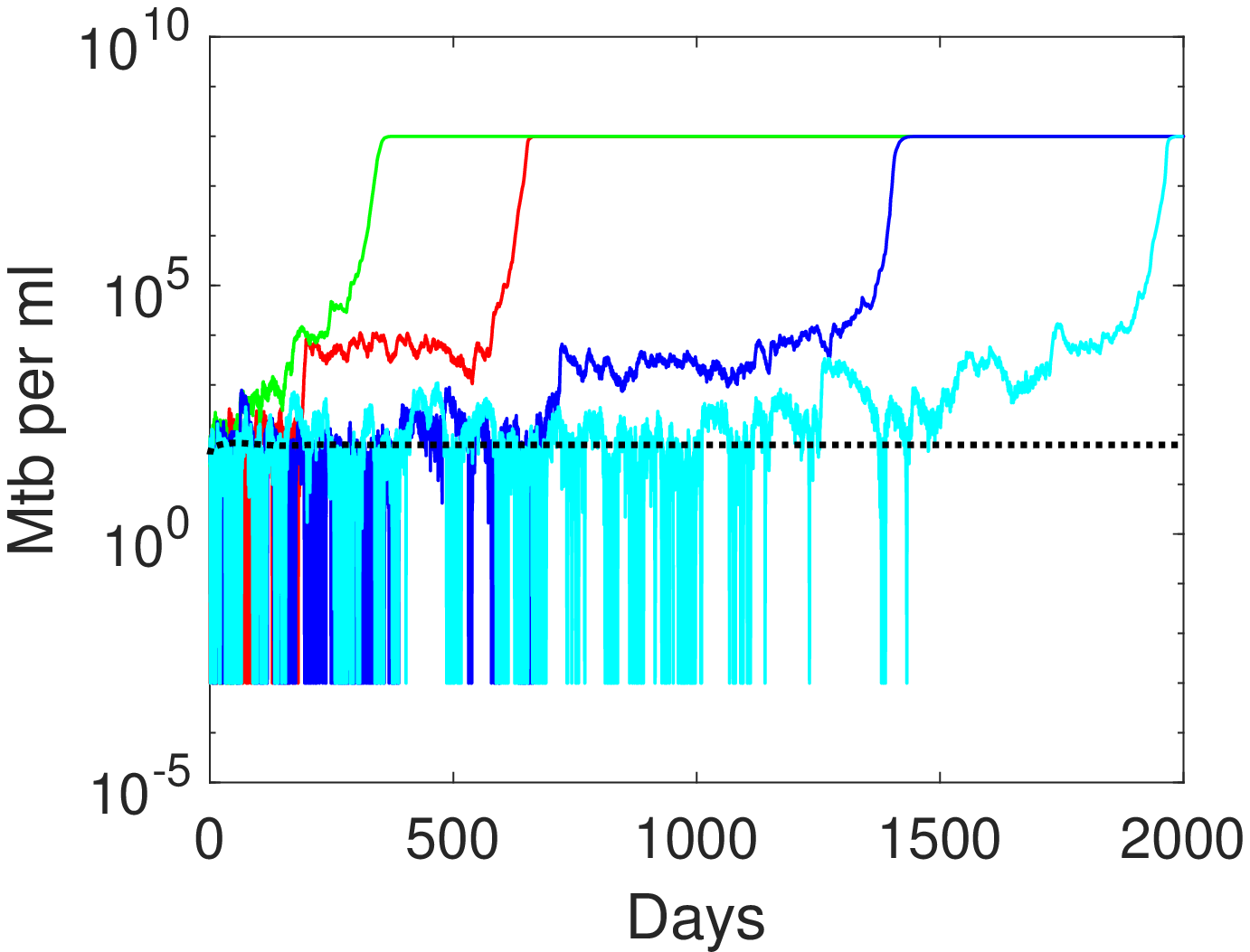}
    \caption{RGN 3, fast return}
    \end{subfigure}
     \begin{subfigure}{0.24\textwidth}
    \includegraphics[width=1\textwidth]{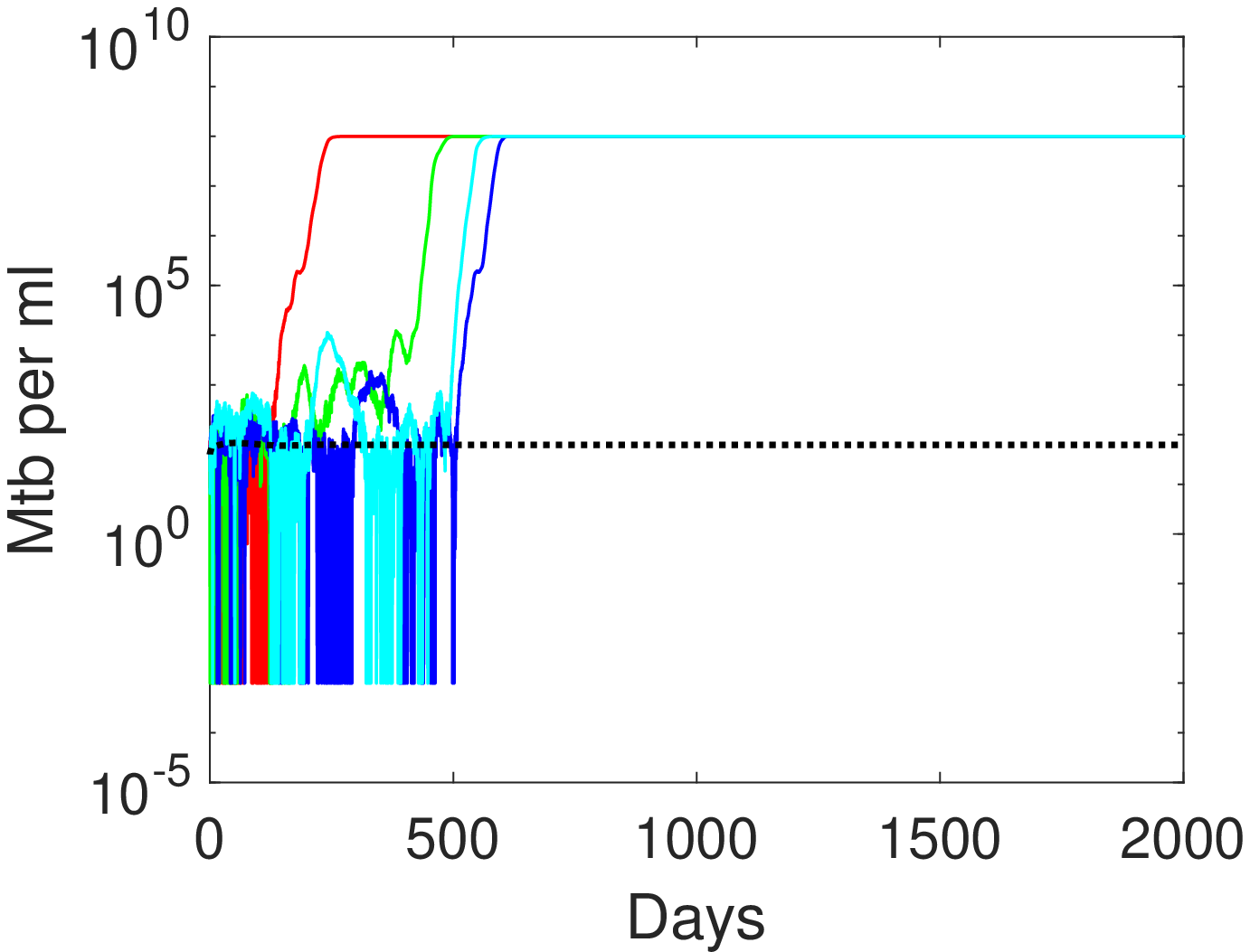}
    \caption{RGN 3, slow return}
    \end{subfigure}\\
        \begin{subfigure}{0.24\textwidth}
    \includegraphics[width=1\textwidth]{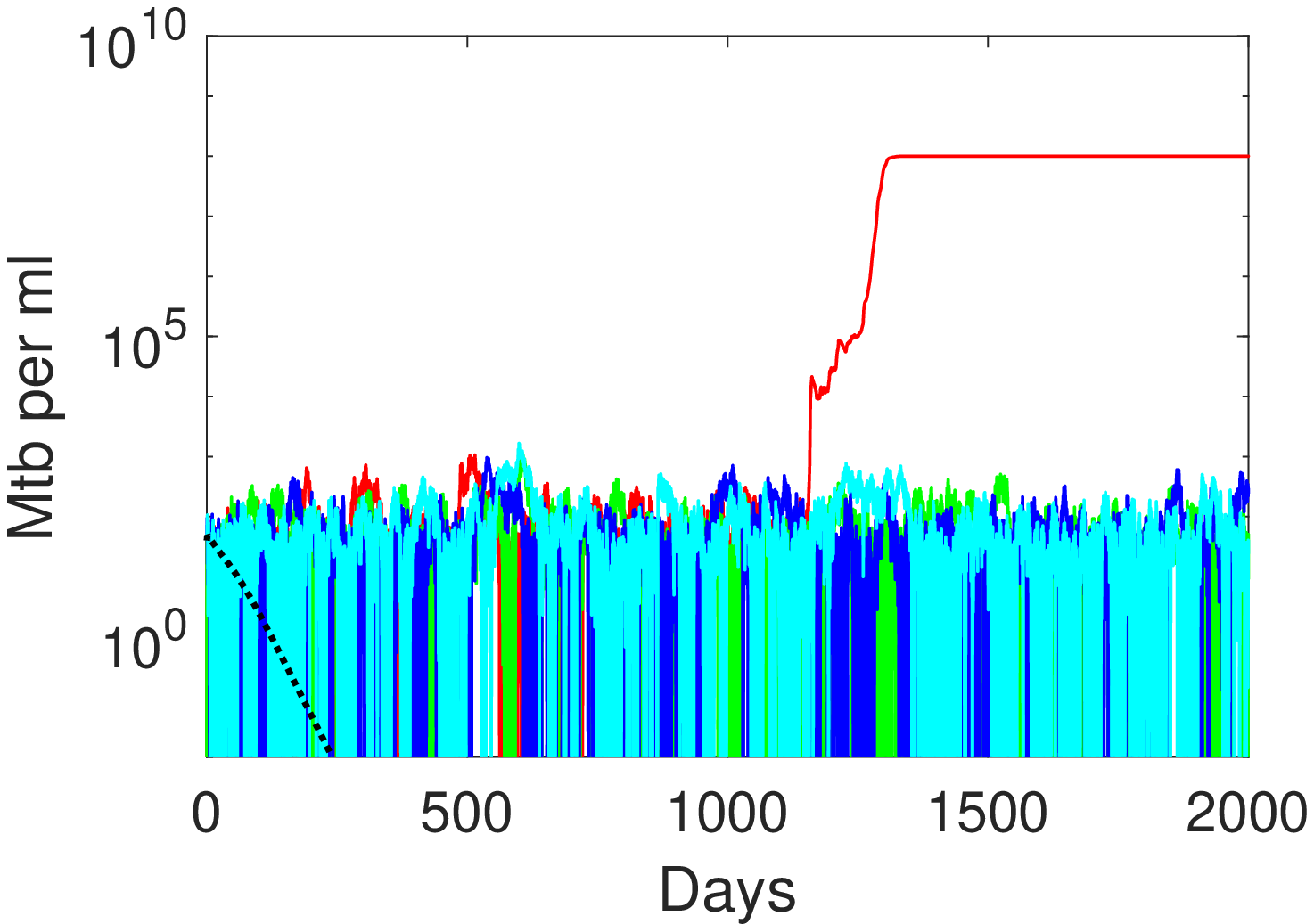}
    \caption{RGN 2, fast return}
    \end{subfigure}
     \begin{subfigure}{0.24\textwidth}
    \includegraphics[width=1\textwidth]{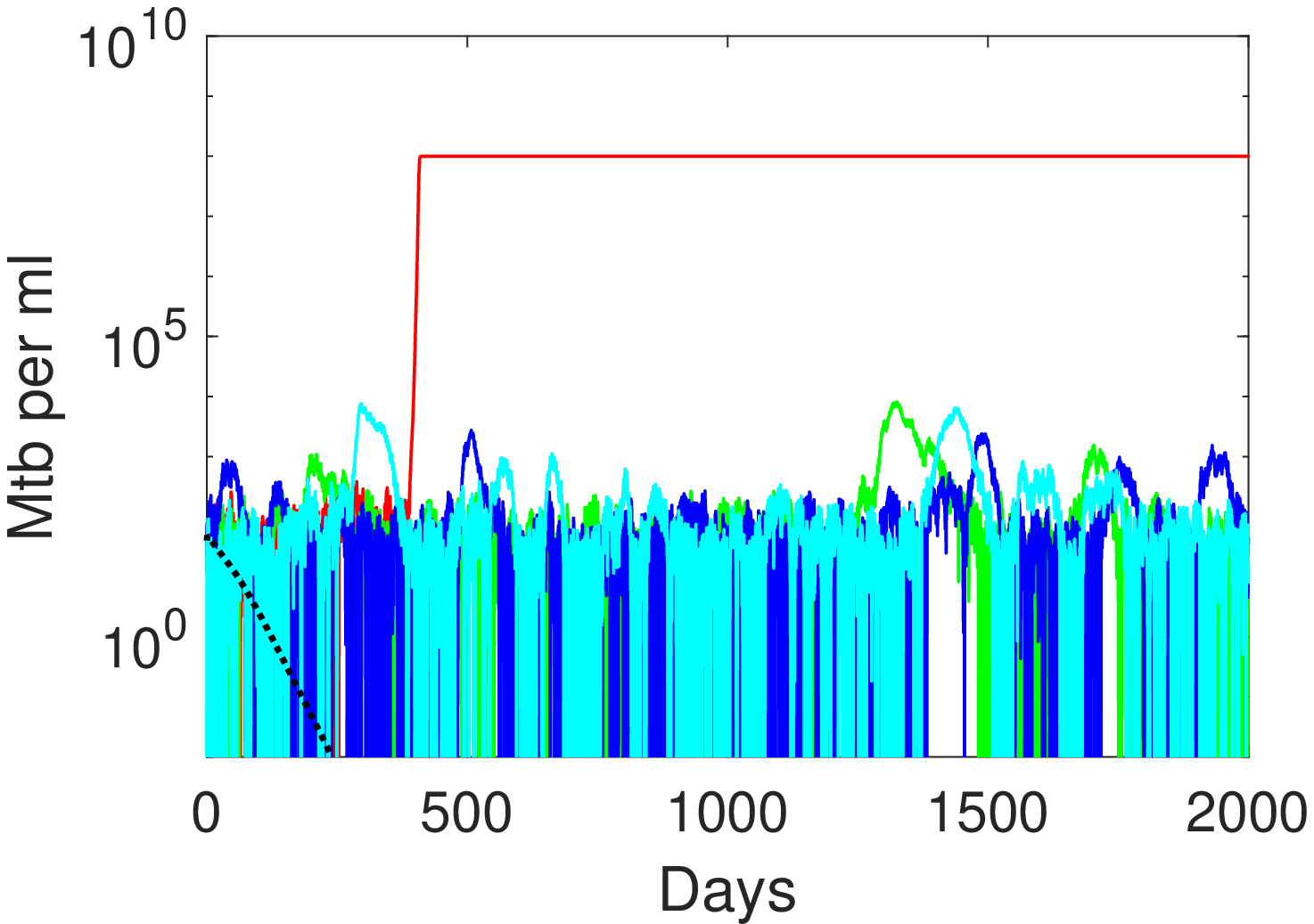}
    \caption{RGN 2, slow return}
    \end{subfigure}
    \begin{subfigure}{0.24\textwidth}
    \includegraphics[width=1\textwidth]{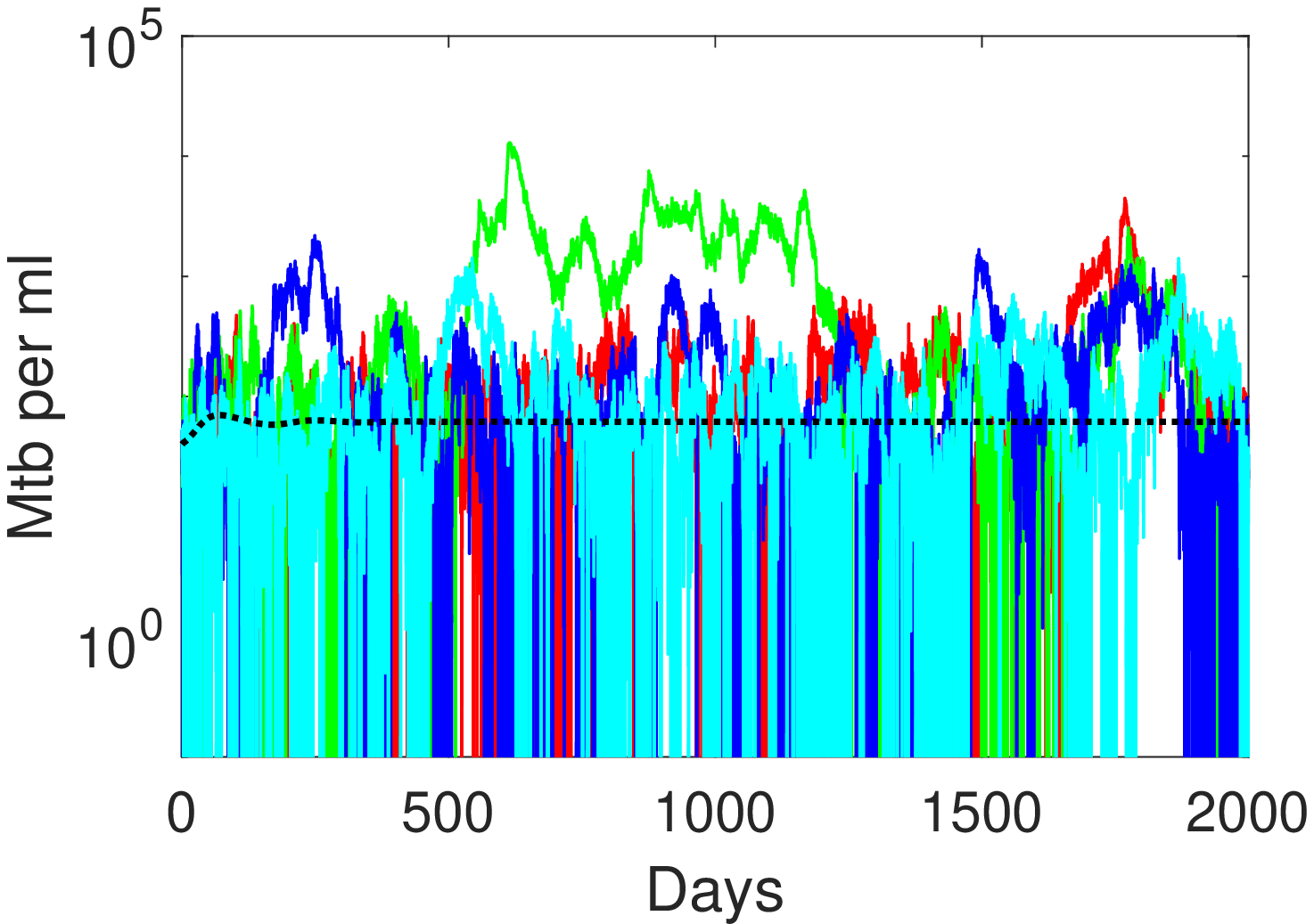}
    \caption{RGN 3, fast return}
    \end{subfigure}
     \begin{subfigure}{0.24\textwidth}
    \includegraphics[width=1\textwidth]{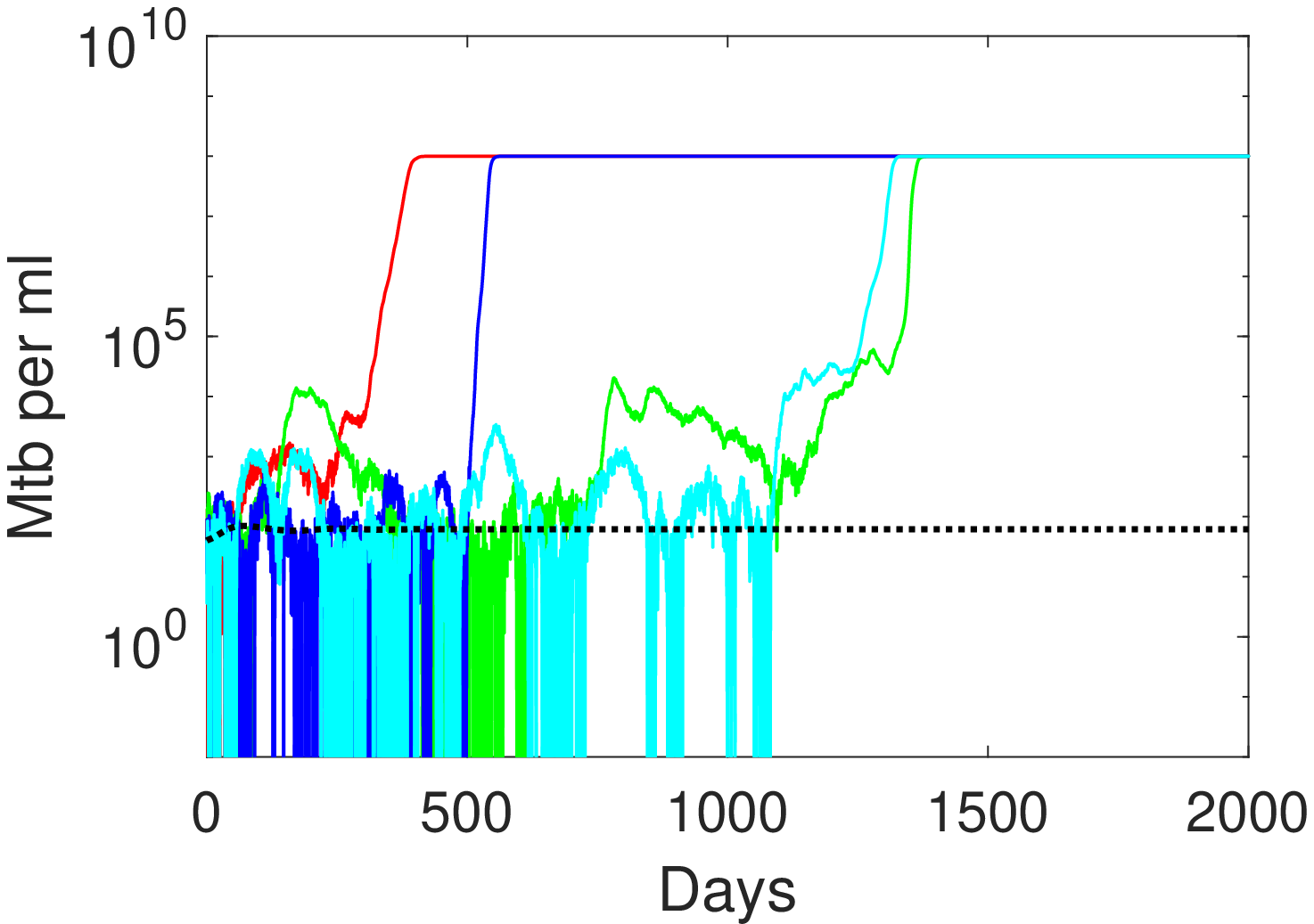}
    \caption{RGN 3, slow return}
    \end{subfigure}\\
        \begin{subfigure}{0.24\textwidth}
    \includegraphics[width=1\textwidth]{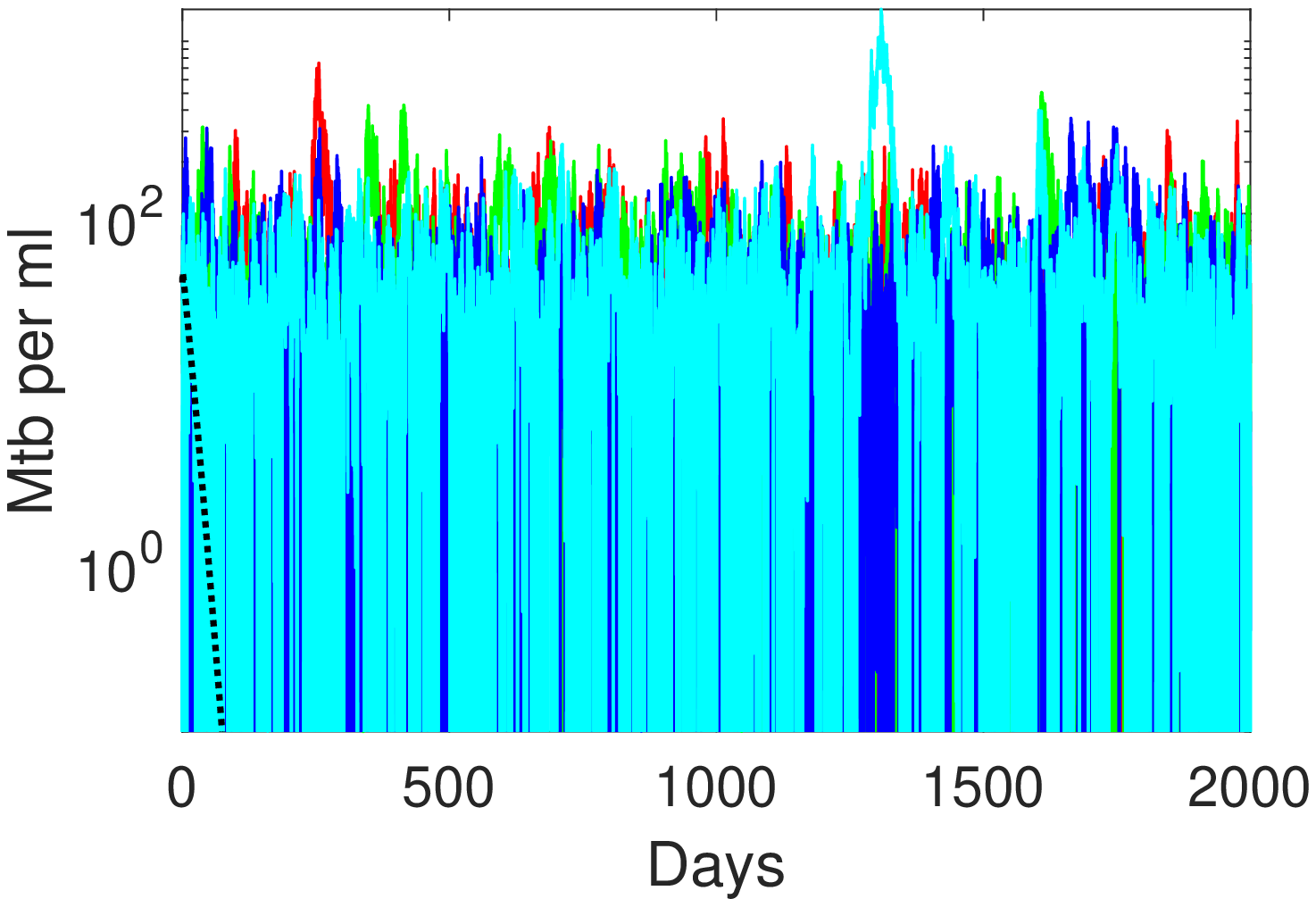}
    \caption{RGN 2, fast return}
    \end{subfigure}
     \begin{subfigure}{0.24\textwidth}
    \includegraphics[width=1\textwidth]{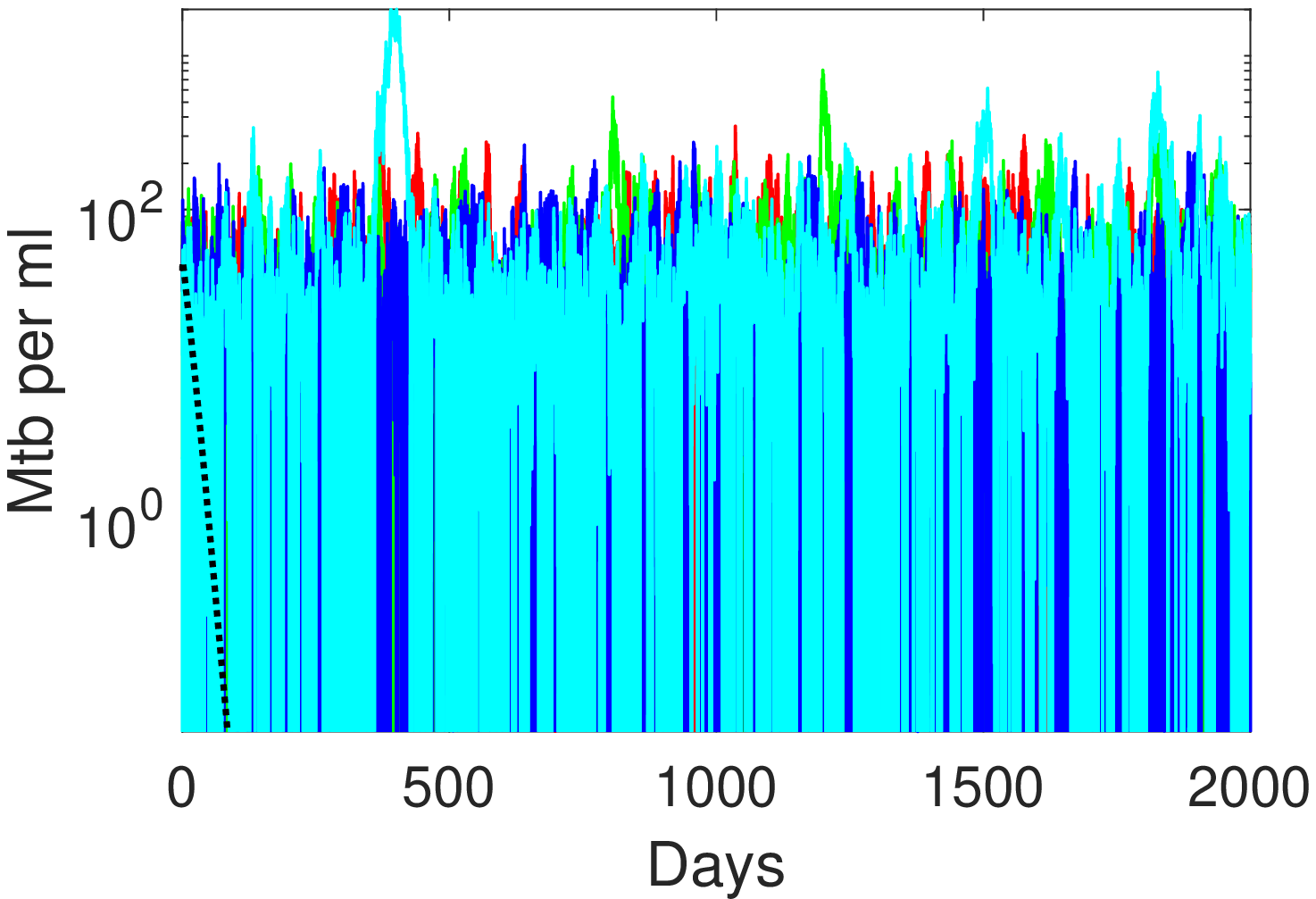}
    \caption{RGN 2, slow return}
    \end{subfigure}
        \begin{subfigure}{0.24\textwidth}
    \includegraphics[width=1\textwidth]{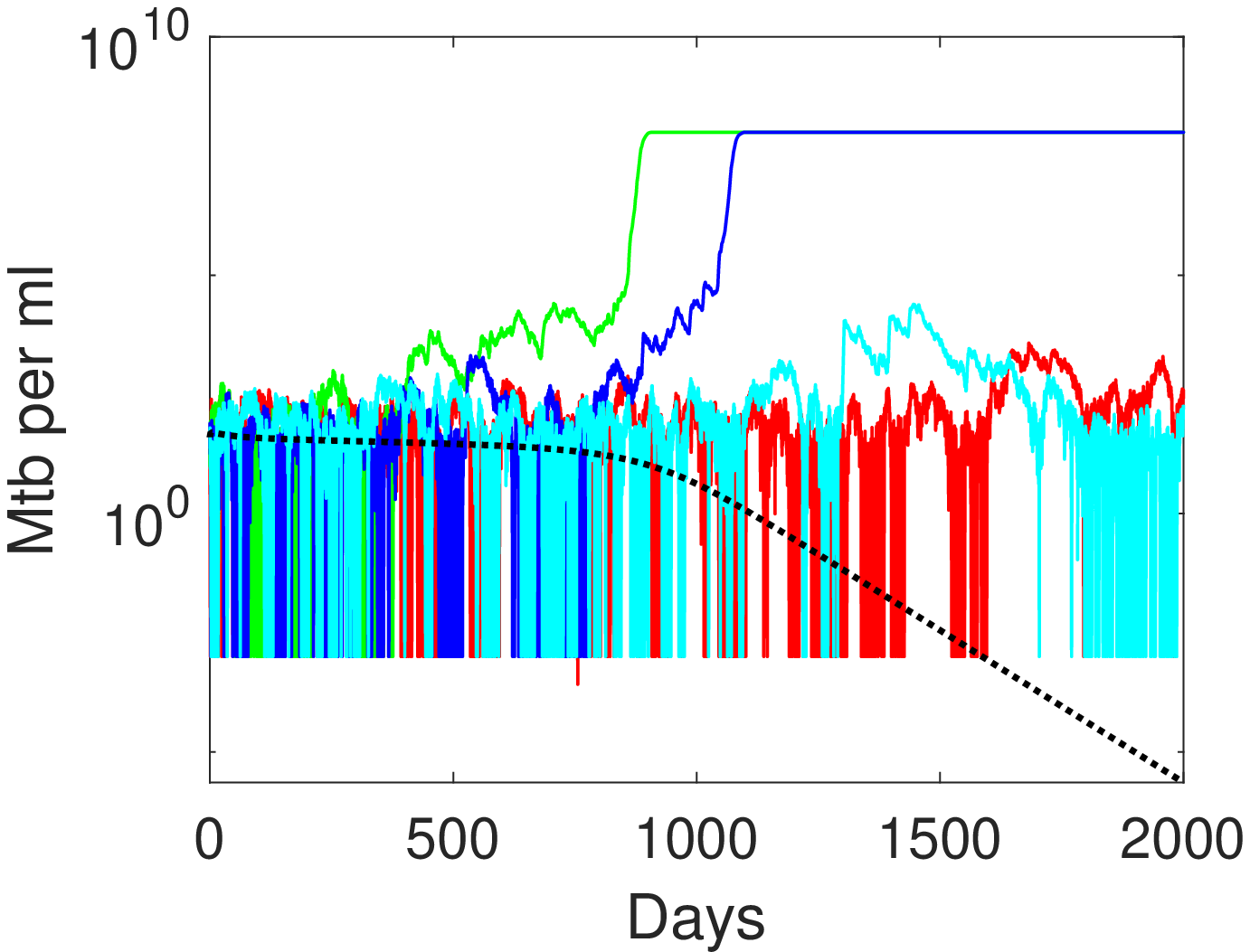}
    \caption{RGN 3, fast return}
    \end{subfigure}
     \begin{subfigure}{0.24\textwidth}
    \includegraphics[width=1\textwidth]{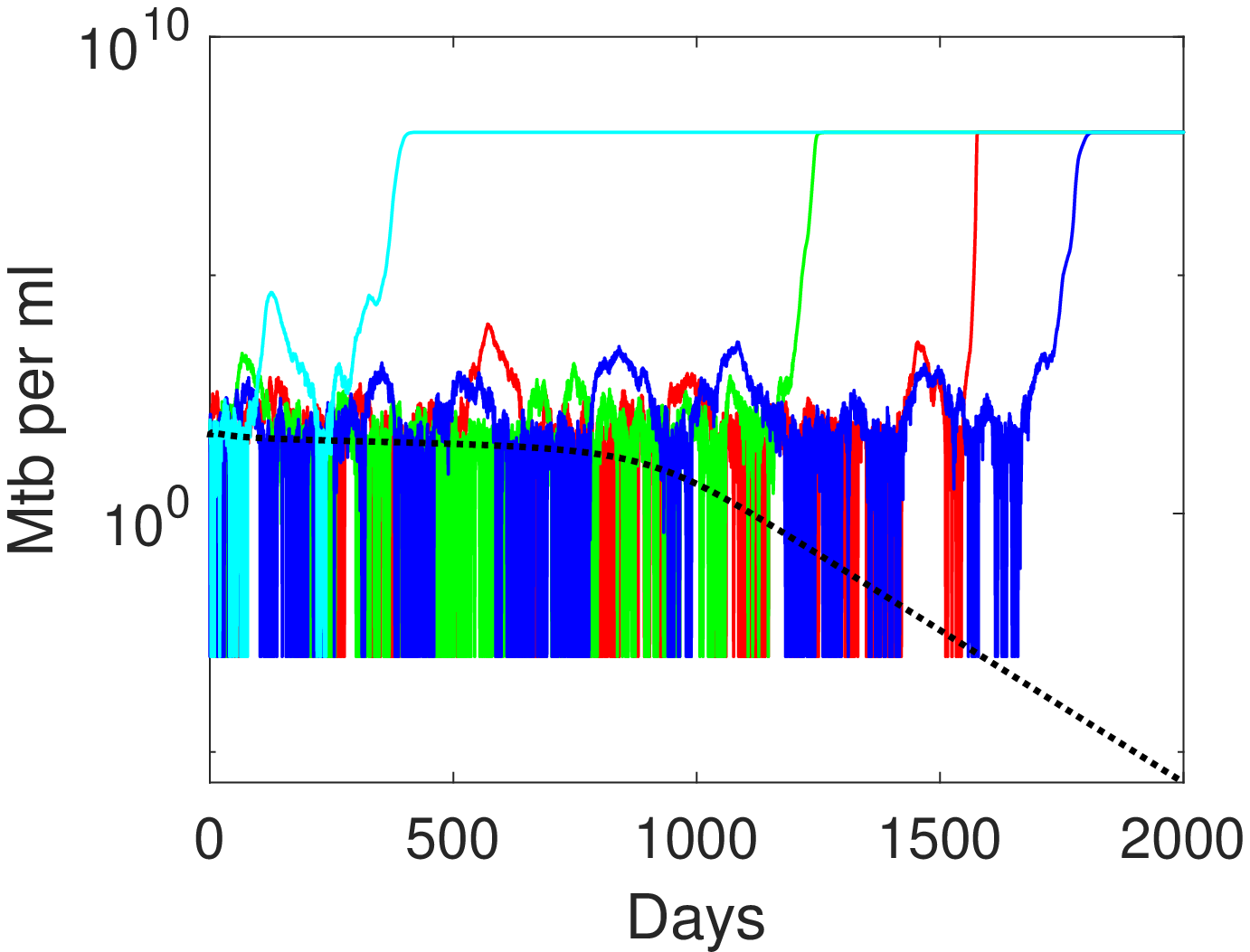}
    \caption{RGN 3, slow return}
    \end{subfigure}
    \caption{Four sample paths for SDEs \eqref{sde2} with demographic and environmental variations is in colored curves. Fast return rate is $\alpha_i=0.5$, slow return rate is $\alpha_i=0.05$, and the return rate satisfies $\alpha_i=2\sigma_i^2$. The ODE solution for model \eqref{eqn1} is in black curve. Parameter values are (1) (in the first row) $\delta =0.2$, $b=0.1$ (Region 2) and $b=0.17$ (Region 3), (2) (in the second row) $\delta =0.2$, $\gamma=1.5$ (Region 2) and $\gamma=1.05$ (Region 3), and (3) (in the third row) $\eta =0.5 \times 10^{-7}$, $\delta=0.1$ (Region 2) and $\delta=0.285$ (Region 3). The other parameter values are from Table \ref{tab1}.}
    \label{fig_sde4}
\end{figure}

Host-directed therapies can modulate host immunity, which in turn tune the associated parameter values. However, these parameter values are influenced by environmental variations due to the patient system's natural regulation and the stability of  adjunct drug concentrations. We then use the second SED model \eqref{sde2} that considers both demographic and environmental variations to evaluate the Mtb concentration. The mean values $\delta_s$, $b_s$, $\gamma_s$, and $\eta_s$ are taken as the values of $\delta$, $b$, $\gamma$ and $\eta$ in Table \ref{tab1}.
We take a return rate relation as $\alpha_i=2\sigma_i^2$ for $i=1,...,4$ to avoid a large variability during the therapy. 
Fast and slow return rates are taken as $\alpha_i=0.5$ and $\alpha_i=0.05$ for simulations of SEDs \eqref{sde2}. 
For the case that the pathogen-directed and host-directed therapy combination affects
the loss rate of infected macrophages $b$ and the bacterial proliferation rate $\delta$, disease outcomes are predicted by bifurcation analyses of the ODE model \eqref{eqn1} in Figure \ref{fig_bif_del_2d} (a), and by the second SDE model \eqref{sde2} in Figure \ref{fig_sde4}. In Region 2, none of the four Mtb sample paths eventually blow up for the fast return rate, but two of them develop to high Mtb level if the return rate is slow. In Region 3, all four Mtb sample paths blows up. But the sample paths take off slower for the fast return rate than the slow return rate.
We do not consider Region 4. Because TB will definitely develop to active disease in this region, we expect therapies will bring parameters away from Region 4.
For the therapies targeting parameters $\gamma$ and $\delta$ combination and $\eta$ and $\delta$ combination, Figures \ref{fig_sde4} indicate similar results.
Our results suggest that a pathogen-directed and host-directed therapy combination can bring parameter values to Region 2  with 
  a fast return rate (i.e., drug concentration can quickly converge to the desired level). 
As a result, the patient can spend more time in the disease clearance or latent infection status. 
This indicates this therapy combination slows down disease progression.

\subsection{A Positive Relationship Between the Bacterial Proliferation and the Slope of Disease Progression}
\begin{figure}
 \centering
    (a)\hspace{8cm}(b)\\
    \includegraphics[width=0.46\textwidth]{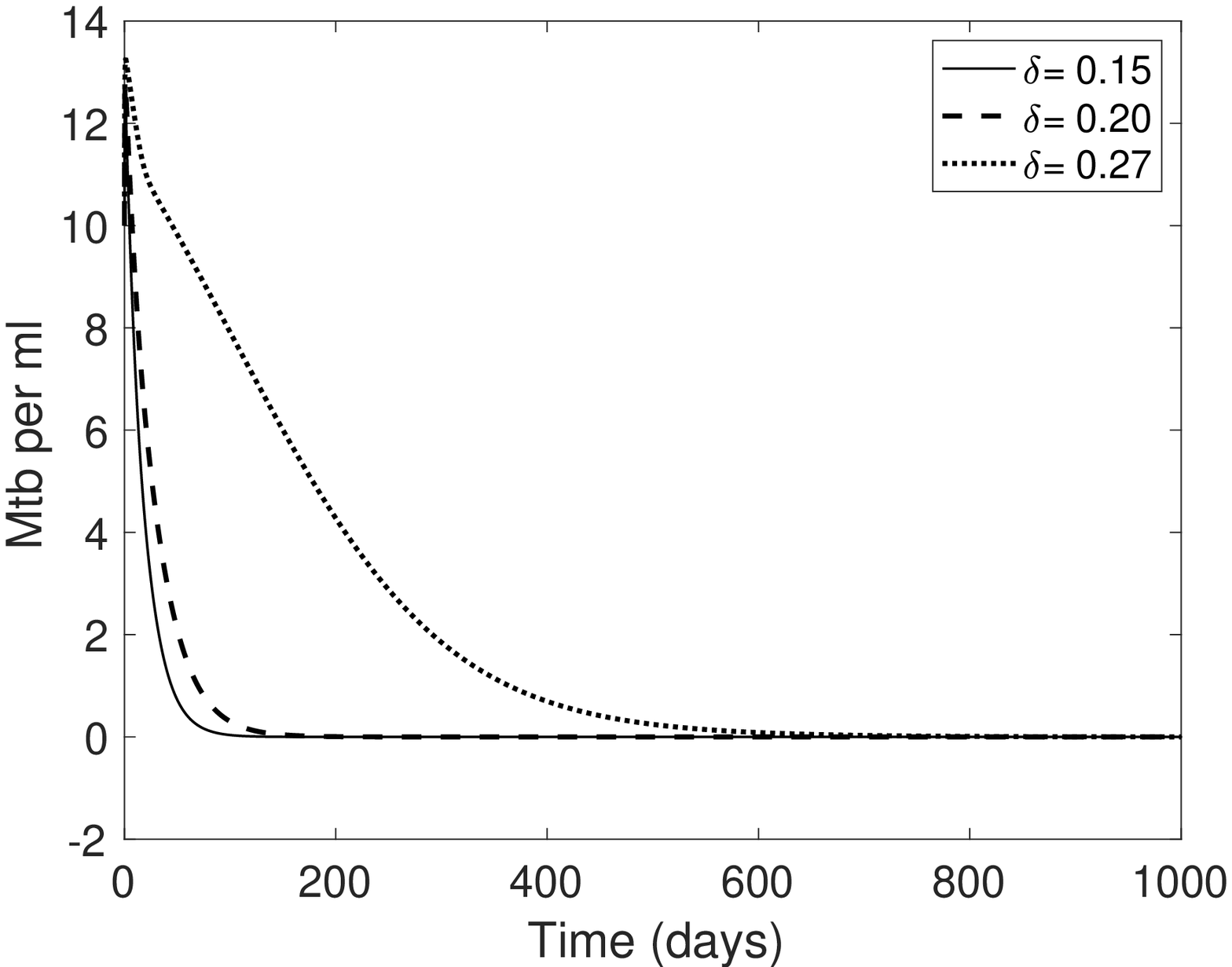}
    \includegraphics[width=0.46\textwidth]{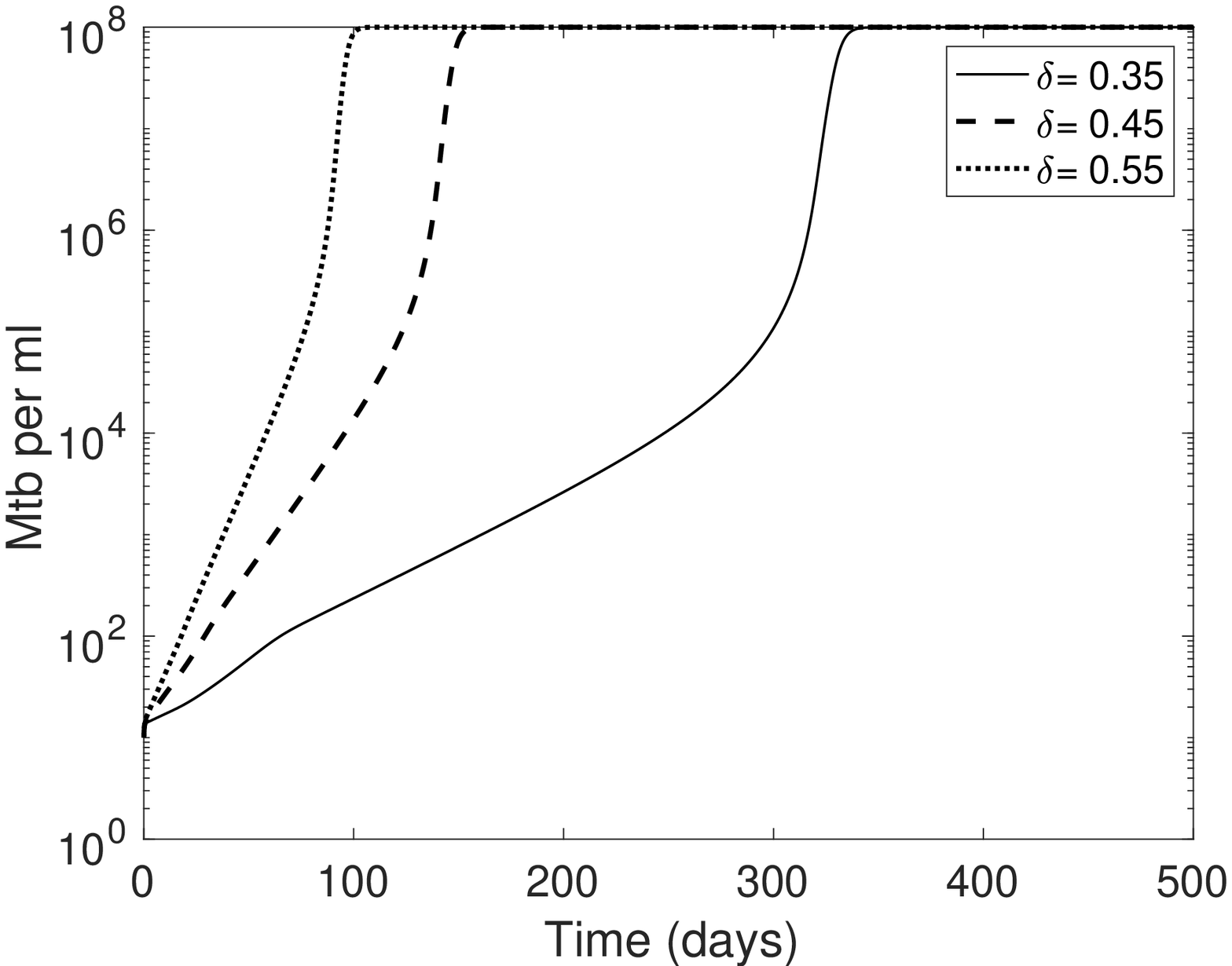}
    \caption{Predicted bacteria concentration decay and growth for a case that a healthy individual ($M_u(0)=\Bar{M}_u=500000 \,cells/ml$ and
    $T(0)=\Bar{T}=1000 \,cells/ml$) inhales  infectious droplet nuclei carrying small amount of Mtb pathogens ($M_i(0)=1\, cells/ml$ and $B(0)=10\, cells/ml$).
    The threshold bacterium growth rate is $\delta_0=0.2956/B/day$.
    $\delta<\delta_0$ results early clearance.  $\delta>\delta_0$ causes disease progression. The slope of the Mtb time history is positive, related to the bacterium growth rate $\delta$. The other parameter values are taken from Table \ref{tab1}.}
    \label{fig_lambda1}
\end{figure}

Figure \ref{fig_lambda1} shows that the bacterial proliferation rate has a positive relationship with the slope of disease clearance and progression after exposure.
More precisely, the slope of bacterium concentration is determined by the dominant eigenvalue. Its linear approximation is a quadratic function of the bacterial proliferation rate $\delta$. 
The following derivation is inspired by methods used by \cite{perelson1999mathematical} and \cite{yu2005closed}.

Tuberculosis infection starts from the inhalation of droplet nuclei carrying Mtb by a healthy individual. The small amount of inhaled tubercle bacilli can be treated as a small positive perturbation of the 
trivial steady state, $(\Bar{M}_{u0},\,\Bar{M}_{i0},\,\Bar{B}_0,\,\Bar{T}_0)$.
The local stable trivial steady state represents the possible scenario
of bacterial clearance, while the unstable trivial steady state means insufficient immune responses fail to clear
the infecting bacteria, and thus disease starts to progress.
Figure \ref{fig_lambda1} shows that the rate of early clearance and disease progression depends on the bacterial proliferation rate $\delta$. 
A higher bacterial proliferation rate results in a slower bacterial decline and a faster disease progress. 
It also suggests the existence of a threshold value
of the bacterial proliferation rate $\delta$, which separates the early
clearance and disease progression in early infection.

Mathematically, the slope of bacteria decline and progression is dependent on the largest eigenvalue of the model \eqref{eqn1} at the trivial steady state $(\Bar{M}_{u0},\,\Bar{M}_{i0},\,\Bar{B}_0,\,\Bar{T}_0)$.
Consider a general system $\Dot{u}=f(u,\,\mu)$, where $u \in \mathbb{R}^n$, $\mu \in \mathbb{R}^m$, and $f \in \mathbb{R}^{n\times m}$ denote state variables, parameters, and
a vector of n functions. Suppose that $u_e=u_e(\mu)$ represents the trivial steady state 
, i.e., $f(u_e(\mu),\,\mu)=0$.
The small amount of invading Mtb pathogens is modeled as a small positive perturbation nearby the trivial steady state, i.e., $x=u-u_e \in \mathbb{R}^n$.
The growth of the perturbation indicates the development of the infection,
while a decay of $x$ corresponds to the disease clearance.
The evolution of disease progression is governed by the perturbation of the equation as follows:
\begin{equation*}
         \Dot{x}=\Dot{u}-\Dot{u}_e = f(x+u_e(\mu),\,\mu)
         = f(u_e(\mu),\,\mu) + x\,\left[ \dfrac{\partial f_i(u_e(\mu),\,\mu)}{\partial x_j}\right]_{n\times n}+
         \;\text{higher}\;\text{order}\;\text{terms.}
\end{equation*}
If the preceding system is hyperbolic, that is all eigenvalues of $A=\left[ \partial f_i(u_e(\mu),\,\mu)/\partial x_j\right]_{n\times n}$ have nonzero real parts, then the system can be approximated by the linear system $\Dot{x}=A\,x$. 
Further, we assume $A$ has $n$ distinct eigenvalues
$\lambda_i$ with associated eigenvector $v_i$,
i.e., $(A-\lambda_i I)v_i=0$ for $i=1,\,\dots,\,n$.
We form a matrix $T=(v_1,\,v_2,\,\dots ,\,v_n)$, and
have $A=T\Lambda T^{-1}$. Then the general solution of the linearized system $\Dot{x}=A\,x$ with a  chosen initial condition $x_0\in \mathbb{R}^n$ is 
a linear combination of $v_1 e^{\lambda_1 t},\,
v_2 e^{\lambda_2 t},\, \dots ,\, v_n e^{\lambda_n t}$.
That is 
$x(t)=e^{At}x_0=T e^{\lambda t} T^{-1} x_0=(v_1 e^{\lambda_1 t},\,
v_2 e^{\lambda_2 t},\, \dots ,\, v_n e^{\lambda_n t})T^{-1} x_0$.
Assuming $\lambda_1$ is the dominant eigenvalue (the spectral radius of the matrix $A$), i.e. $\lambda_1 =\rho(A)= \max \{\lambda_i,\,i=1,\,\dots,\,n \}$, then $e^{\lambda_1 t}$ is the most
influential component in the solution basis set of $x(t)$.
The evolution of the perturbation $x(t)$ is the change of the original variable $u(t)$.
It is determined by the dominant eigenvalue $\lambda_1$.
Moreover, if $\lambda_1<0$, the perturbation dies out ($x(t)\rightarrow 0$).
The infecting bacteria are eliminated by host immune responses. 
If $\lambda_1>0$, the perturbation grows and the infection develops.
The threshold for the infection clearance and establishment is 
$\lambda_1=0$ and $\lambda_i<0$ for $i=2,\,\dots,\,n$.

Taking the general system $\Dot{u}=f(u,\,\mu)$ as the model \eqref{eqn1} and its trivial equilibrium as $u_e=(\Bar{M}_{u0},\,\Bar{M}_{i0},\,\Bar{B}_0,\,\Bar{T}_0)$.
The corresponding eigenvalues are roots of the following characteristic polynomial
\begin{equation}
    \begin{array}{rl}
        P(\lambda) &= (\lambda+\mu_T)(\lambda+\mu_M)(\lambda^2+c_1\lambda +c_2)=0, \quad\text{where} \\[2.0ex]
         c_1&= (b+\gamma)+\dfrac{s_M}{\mu_M}(N_3\beta+\eta)-\delta \quad \text{and}\\[2.0ex]
         c_2&=\left\{\eta\dfrac{s_M}{\mu_M}-\beta\left[(N2-N3)\dfrac{\gamma}{b+\gamma}+(N1-N3)\dfrac{b}{b+\gamma}\right]\dfrac{s_M}{\mu_M}-\delta \right\}(b+\gamma).
    \end{array}
    \label{eqn4}
\end{equation}
That is
\begin{equation}
    \lambda_1=-\dfrac{c_1}{2}+\dfrac{c_1}{2}\sqrt{1-\dfrac{4c_2}{c_1}},\;\lambda_2=-\dfrac{c_1}{2}-\dfrac{c_1}{2}\sqrt{1-\dfrac{4c_2}{c_1}},\;\lambda_3=-\mu_T,\;\lambda_4=-\mu_M,
    \label{eqn5}
\end{equation}
where $\lambda_{3,\,4}<0$ due to the positiveness of all parameter values.
The expression of the preceding equation for $P(\lambda)=0$ and $\Bar{M}_{u0}=s_M/\mu_M$ (\cite{Zhang2020TB}) determine the local stability of the trivial equilibrium. 

Here, we use the eigenvalues associated with the trivial equilibrium to further investigate the disease progression speed.
Extracellular bacteria are introduced into the system by three main ways, extracellular bacterial proliferation, bacterium release by the programmed cell death of infected macrophages, 
and bacterium release from infected macrophages killed by T-cell mediated immune responses.
   In early infection, we assume the loss of extracellular Mtb pathogens is  caused only by macrophage phagocytosis. We omit the loss of infected macrophages by T-cell mediated immune responses.
The per capita  rate of macrophage phagocytosis is $N_3\beta+\eta$ per uninfected macrophage per bacterium. 
Considering the (asymptotic) upper bound of the uninfected macrophage $\Bar{M}_{u0}=s_M/\mu_M$,
the maximum bacterial loss rate is $(N_3\beta+\eta)\Bar{M}_{u0}=(N_3\beta+\eta)s_M/\mu_M$ per bacterium,
which is assumed to be greater than the bacterium proliferation rate $\delta$.
  Further, assuming that all invading pathogens are phagocytized, that is $\delta<(N_3\beta+\eta)s_M/\mu_M$.
This implies $\delta<b+\gamma+(N_3\beta+\eta)s_M/\mu_M$, which is equivalent to $c_1>0$.
Therefore, if $c_2<0$, the square root $\sqrt{1-\frac{4c_2}{c_1}}$ is greater than one, then $\lambda_1$ and $\lambda_2$ are real numbers and have opposite signs. If $c_2>0$, then the real parts of $\lambda_1$ and $\lambda_2$ are all negative.

If $c_2=0$, a zero-eigenvalue bifurcation occurs at $\lambda_1=0$, at which
\begin{equation}
    \delta=\delta_0\overset{\Delta}{=}\eta\,\dfrac{s_M}{\mu_M}-(N2-N3)\,\dfrac{\gamma}{b+\gamma}\,\beta\,\dfrac{s_M}{\mu_M}-(N1-N3)\,\dfrac{b}{b+\gamma}\,\beta\,\dfrac{s_M}{\mu_M}.
\end{equation}
The number of “next generation” infectious Mtb bacteria produced by a single infectious bacterium introduced near the trivial equilibrium is (1) at most $\delta$ by cell proliferation, or (2) $(N_1-N_3)$ and $(N_2-N_3)$ by the  death of infected macrophages and T-cell mediated immune response with the probability of $b/(b+\gamma)$ and $\gamma/(b+\gamma)$, respectively. 
In the meantime, phagocytosis can kill at most $\eta$ engulfed bacteria. We assume the uninfected macrophages is at the maximum available level $s_M/\mu_M$. Therefore, the infection dies out if $\delta<\delta_0$  (i.e., $c_2>0$), but stays if
$\delta>\delta_0$ (i.e. $c_2<0$, $\lambda_1>0$, $\lambda_2<0$). 
The threshold is $\delta=\delta_0$.
Considering the parameter values in Table \ref{tab1},
the bifurcation analysis shown in Figure \ref{fig_bif_sim_B} (a) shows
that the $\delta$ parameter range for latent infected TB (roughly $80\%$ of the total TB infected individuals) is $[0.2621,\,0.2944]$,
which is very close to $\delta_0=0.2957$. That is for $\delta \in [0.2621,\,0.2944]$, $|\delta-\delta_0|\ll 1$ (then $|c_2|\ll 1$), we thus expand the square root in $\lambda_1$ in Equation \eqref{eqn5} at
$\delta=\delta_0$. 
This yields
\begin{equation}
\begin{array}{rl}
    \sqrt{1-\dfrac{4c_2(\delta)}{c_1(\delta)}} &=
    1-2\,\dfrac{1}{c_1(\delta_0)}\,[-(b+\gamma)]\,(\delta-\delta_0) + \cdots ,\quad\text{thus}\\[3.0ex]
    \lambda_1 &= -\dfrac{c_1(\delta)}{2} + \dfrac{c_1(\delta)}{2}
    \left\{ 1-2\,\dfrac{1}{c_1(\delta_0)}\,[-(b+\gamma)]\,(\delta-\delta_0) + \cdots \right\},\\[3.0ex]
    &\approx (b+\gamma)\,\dfrac{c_1(\delta)}{c_1(\delta_0)}\,\left(\delta-\delta_0\right)= -\dfrac{c_1(\delta)c_2(\delta)}{c_1(\delta_0)}.
    \end{array}
    \label{eqn7}
\end{equation}

\begin{figure}
  \begin{center}
   (a)\qquad\qquad\qquad\qquad\qquad\qquad\qquad\qquad\qquad\qquad\qquad(b)\\
    \includegraphics[width=0.48\textwidth]{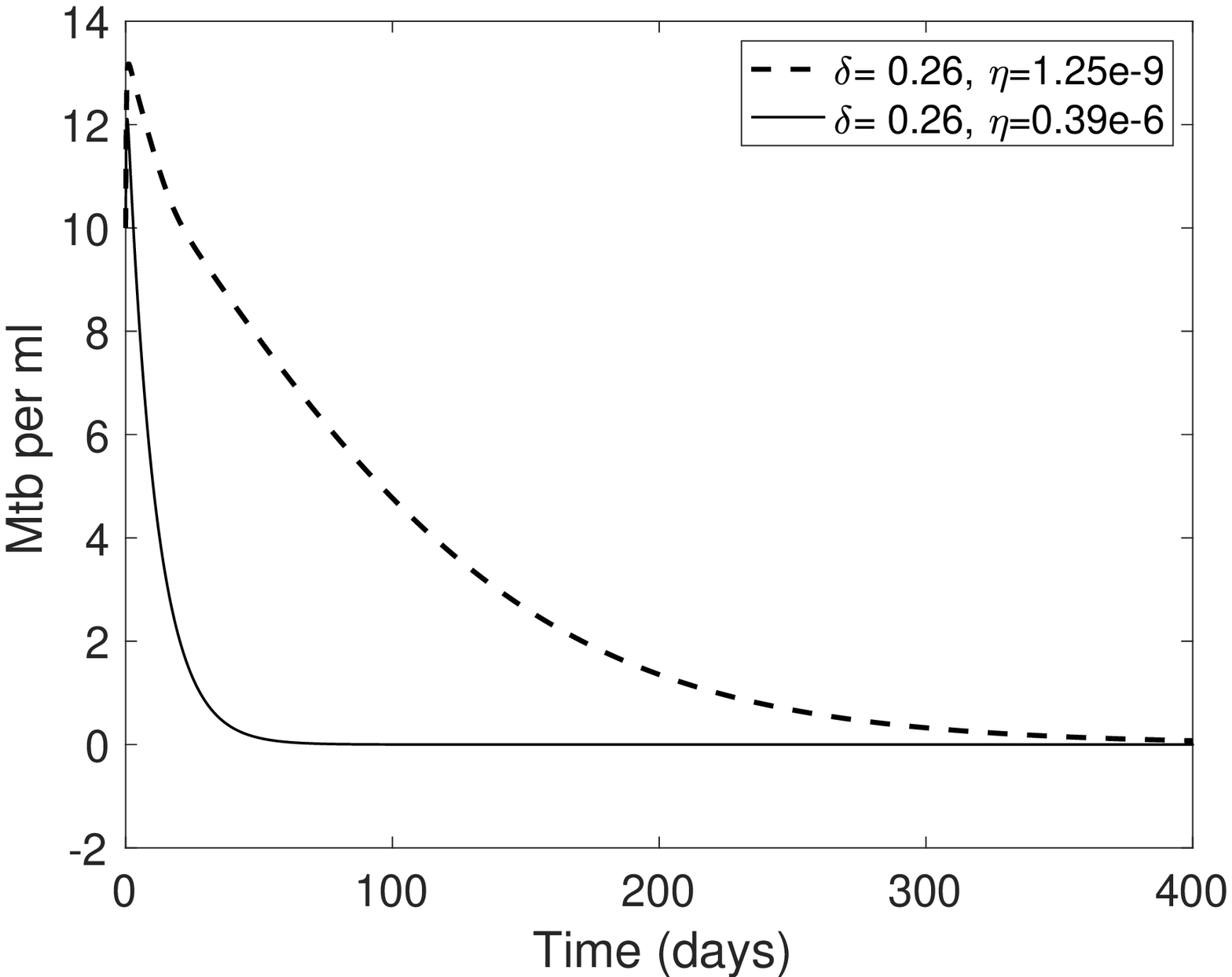}
    \includegraphics[width=0.48\textwidth]{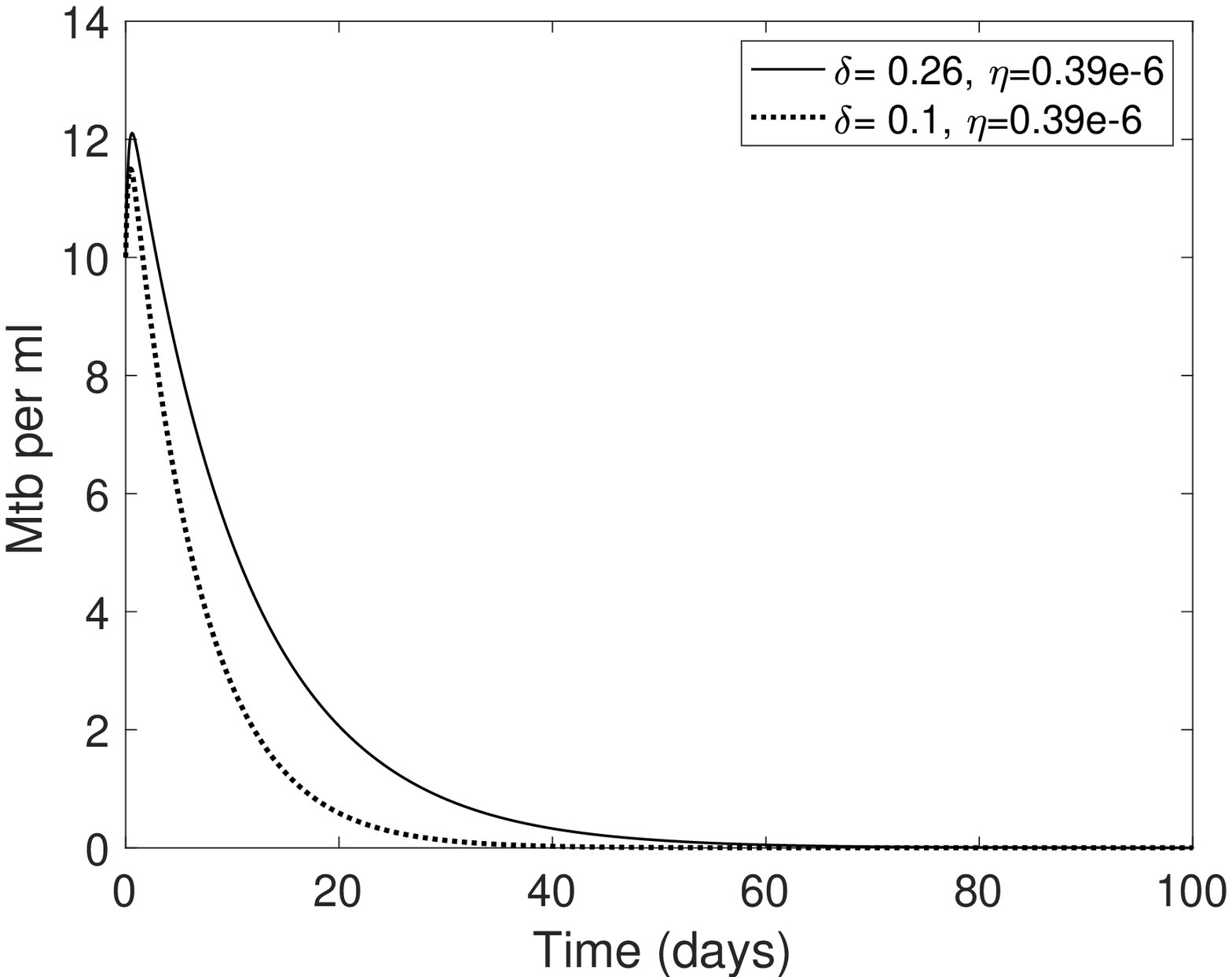}
  \end{center}
  \caption{Bacterial clearance (a) at macrophages killing rate $\eta$ with (solid) and without (dashed) vitamin D boost, and (b) with only vitamin D boost (solid) and with both vitamin D and antibiotic therapies (dotted).   Initial conditions are taken for a case that a healthy individual ($M_u(0)=\Bar{M}_u=500000 \,cells/ml$ and
    $T(0)=\Bar{T}=1000 \,cells/ml$) inhales  infectious droplet nuclei carrying small amount of Mtb pathogens ($M_i(0)=1\, cells/ml$ and $B(0)=10\, cells/ml$).
  The other parameter values are taken from Table \ref{tab1}.}
  \label{fig_eta}
 \end{figure}

Base on the feasible parameter ranges, $c_1>0$ holds.
We approximate the dominant eigenvalue $\lambda_1$ as a quadratic function of the bacterial proliferation rate $\delta$.
We notice that the dominant eigenvalues determines the changing slope of the bacteria concentration. 
This implies that if $\delta<\delta_0$ (or $\delta>\delta_0$), the Mtb bacterial concentration declines (or grows) in a speed with a quadratic relationship of the bacterial proliferation rate $\delta$.
This speed is slowing down if $\delta$ is approaching to its critical value $\delta_0$ from both sides, but is speeding up if $\delta$ is moving away from $\delta_0$. 
This is demonstrated by numerical simulations in Figure \ref{fig_lambda1}.

\subsection{Using Vitamin D as an Adjunctive Therapy}
Vitamin D has been shown to promote macrophage maturation, which in turn inhibits the intracellular Mtb growth and enhance antimicrobial immune response (\cite{de2019immune,arranz2017host,gombart2009vitamin,rosenberger2004interplay}). 
A dose-dependent vitamin D induced reduction in Mtb growth
 can reach $ 75.7\%$ (\cite{martineau2007ifn}).
The strengthened antimicrobial immune response can be described as an enhanced bacteria killing rate by macrophages, $\eta$.
Equations \eqref{eqn4} and \eqref{eqn7} indicate that the dominant eigenvalue  $\lambda_1$  has a negative quadratic relationship with $\eta$.
We experiment with different parameter values for $\eta$ and $\delta$ to demonstrate the effect of vitamin D therapy and the combination therapy of antibiotic and vitamin D.
We adopt the parameter ranges for $\delta \in (0,\,0.35)$ and $\eta \in (1.25 \times 10^{-9},\,1.25 \times 10^{-7})$ (\cite{du2017simple}).
 The simulation in Figure \ref{fig_eta} (a) shows that the bacterial decrease slope is much steeper for the solid curve than the dashed one. 
This supports the idea of using vitamin D as an adjunct therapy for enhanced immune response and obtaining a favorable disease outcome.
 Figure \ref{fig_eta} (b) demonstrates that the dotted curve with a reduced value of $\delta$ shows a faster decay rate. 
This confirms the promising therapy outcomes of adding vitamin D to modulate the host immune response alongside the  antibiotic therapy (\cite{tobin2015host}).

\begin{figure}[h]
    \centering
    \begin{subfigure}{.45\textwidth}
        \centering
        \includegraphics[width=1\textwidth]{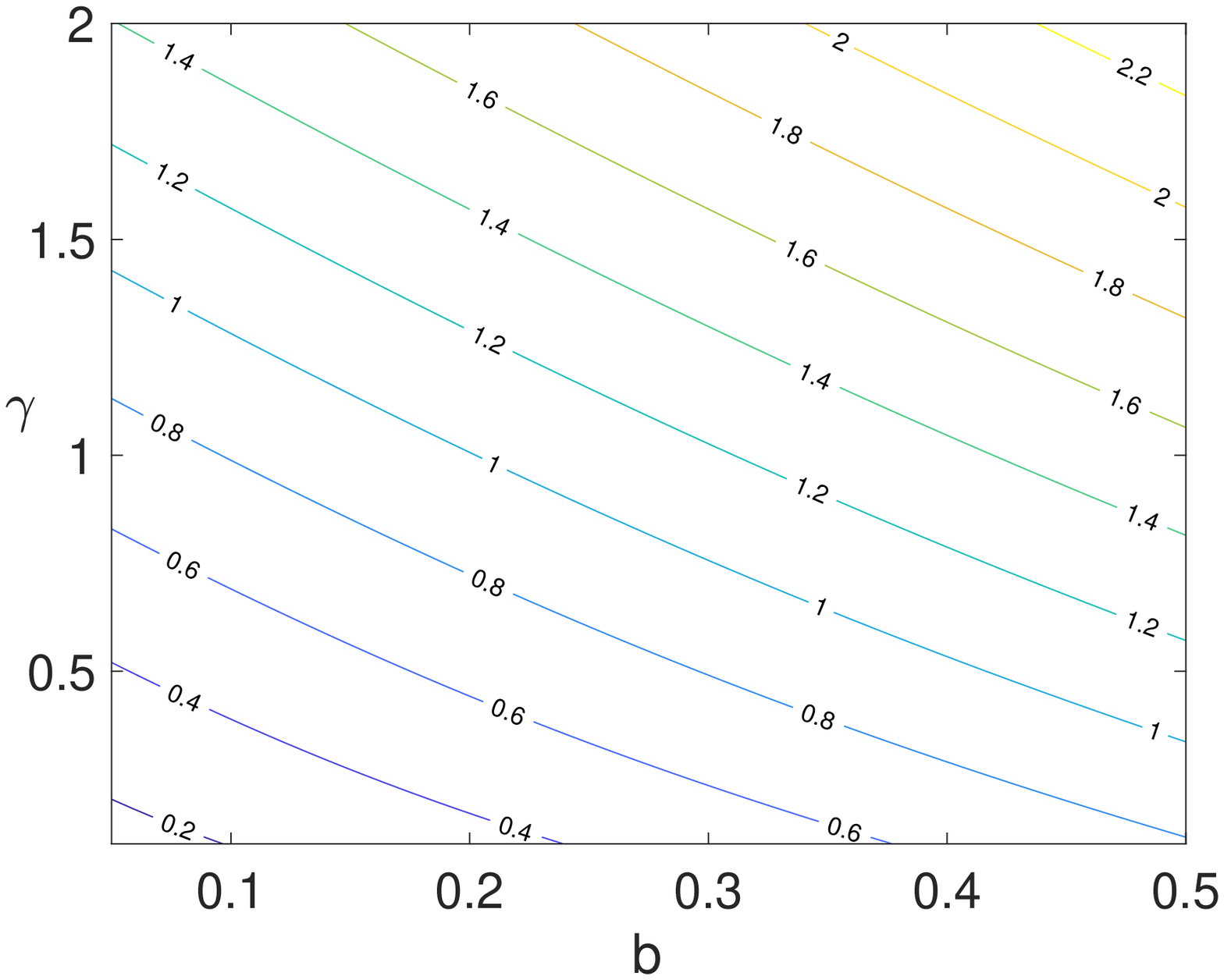}
            \caption{ $N_1>N_3$ and $N_2>N_3$}
    \end{subfigure}
        \begin{subfigure}{.45\textwidth}
            \centering
     \includegraphics[width=1\textwidth]{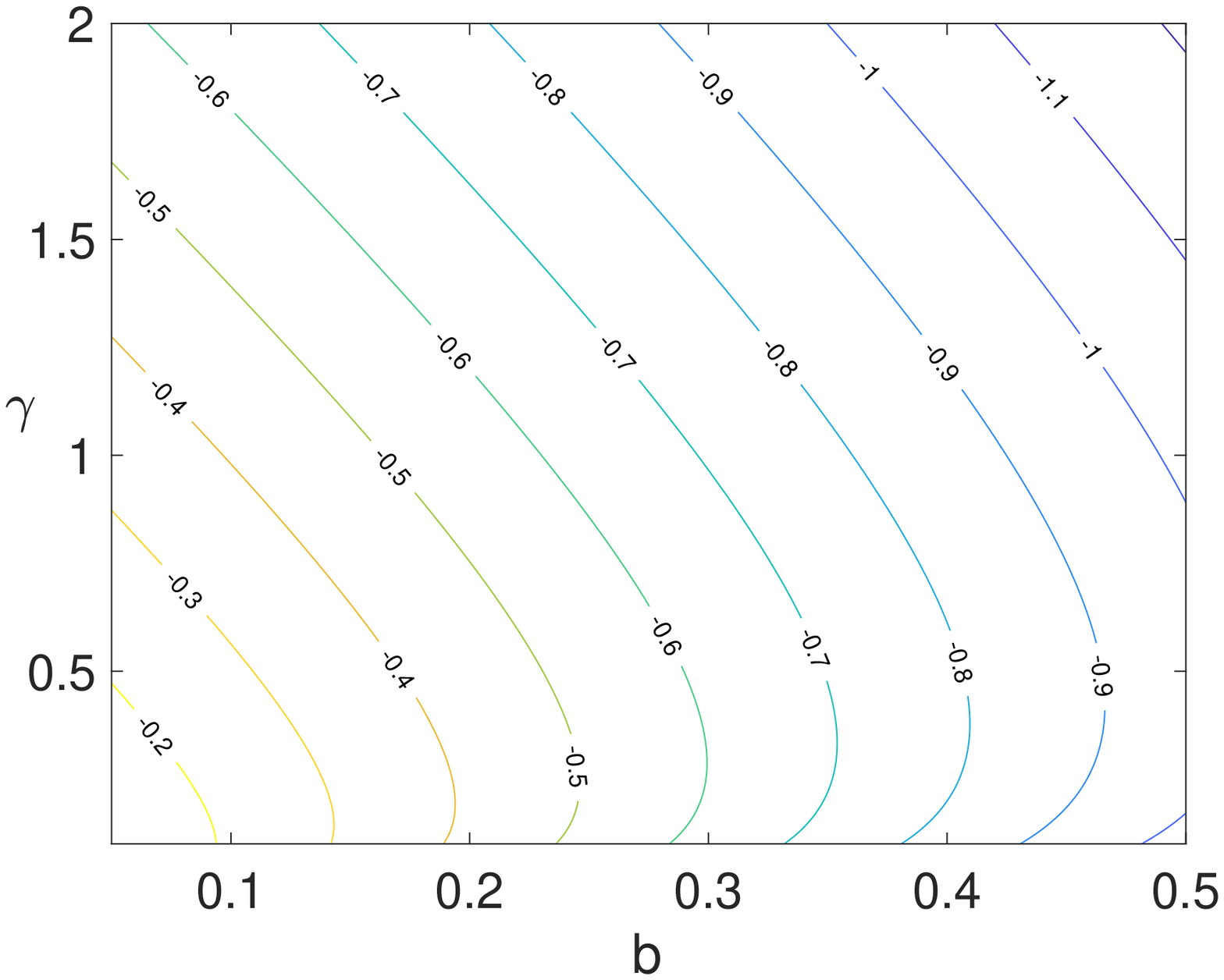}
         \caption{$N_1<N_3$ and $N_2<N_3$}
    \end{subfigure}
           \caption{Contour plots of the dominant eigenvalue $\lambda_1$ as a function of $b$ and $\gamma$. The parameter values for $N_1$, $N_2$ and $N_3$ used in simulations are 
         (a) $N_1=50$, $N_2=30$, and $N_3=25$, (b) $N_1=10$, $N_2=20$ and $N_3=25$. The other parameter values are taken from Table \ref{tab1}.}
    \label{fig_lam1}
\end{figure}

 Moreover, vitamin D has been suggested to play a pivotal role on interferon-$\gamma$ (IFN-$\gamma$)-induced antimicrobial pathway in macrophages against Mtb infection (\cite{fabri2011vitamin}).
In vitamin D-sufficient sera, IFN-$\gamma$ promotes the production of antimicrobial peptides such as Cathelicidin by macrophages (\cite{liu2006toll}).
It in turn enhances the antimicrobial activity by lowing the intracellular survival (\cite{wheelwright2014all}).
  Reducing intracellular Mtb survival induces a reduced level of bacterial release in  T-cell mediated cell death. 
That means a decrease on the values of $N_2$, which denotes average numbers of bacteria released by  T-cell-induced death of infected macrophage.
  It is also reported that vitamin D promotes infected cell apoptosis (\cite{chakraborti2011vitamin}), which induces a decrease on the values of both $N_1$ and $b$.
Note that $N_3$ denotes the average number of bacteria phagocytized by an uninfected macrophage.
The parameter ranges for $N_1$, $N_2$, and $N_3$  from \cite{du2017simple} are $N_1\,\in\,(50,\,100)$, $N_2\,\in\,(20,\,30)$, and $N_3\,\in\,(25,\,50)$.
The loss of infected macrophages caused by  cell death and cell-mediated killing are at the rates of $b$ and $\gamma$. The same paper also provides the parameter ranges for $b\,\in (0.05,\,0.5)$ and $\gamma \,\in (0.1,\,2)$. 
If $N_1>N_3$ and $N_2>N_3$, the simulation in Figure \ref{fig_lam1} (a) shows a positive relation of between the dominant eigenvalue $\lambda_1$ and the parameters $b$ and $\gamma$. 
The positive value of $\lambda_1$ represents the positive slope of the disease progression.
In this case, a large number of intracellular Mtb are released from the death of infected macrophages. It results that an increase in cell death and cell-mediated killing benefits Mtb invasion.
In the case that a small number of intracellular Mtb are released from the death of infected macrophages, we assume that $N_1<N_3$ and $N_2<N_3$.

The simulation in Figure \ref{fig_lam1} (b) shows that $\lambda_1$ is negatively related with the parameters $b$ and $\gamma$.
The negative slope indicates disease clearance.
In this case, most of the intracellular Mtb are killed, and thus the disease is then under control.
Moreover, the magnitude  of $\lambda_1$ is positively related to the infected macrophage loss/bursting rate $b$ and cell-mediated immunity rate $\gamma$. 
This indicates that the stronger the immune responses, the faster the infection is eliminated.
Therefore, immune responses hinder the progression of the infection and benefit the host in this case.
To reduce the intracellular Mtb load, vitamin D supplementation has been suggested as a host-directed adjunctive therapy (\cite{tobin2015host}). 
Vitamin D can help to inhibit intracellular bacterial growth through antimicrobial mechanisms (\cite{martineau2007ifn}), which modulates the immune response to benefit the host defense against the Mtb infection.

\section{Conclusion and Discussion}
Despite the recent advancement of antibiotic TB drugs, antimicrobial resistance and strict medication adherence are still important issues for achieving effective therapeutic outcomes.
As a result, host-directed therapy is proposed as an adjunctive therapy to modulate the host immune response against the Mtb pathogen.
As an emerging and promising approach for attaching the intractable TB disease, a comprehensive understanding of the Mtb-host dynamics is needed to develop treatments that combine pathogen- and host-directed therapies. 

In this paper, we focused on the factors that can be modulated by pathogen- and host-directed therapies to explore the various disease outcomes. 
Being the mainly affected parameter of pathogen-directed therapy, bacterial proliferation rate is shown to have a positive relationship with the disease progression. However, the infected macrophage death rate $b$ and cell-mediated immunity rate $\gamma$ can both benefit and impede the infection, which depends on the relation between the number of bacteria engulfed and released by macrophages,  see Figure \ref{fig_lam1}. 
If macrophages engulf more Mtb bacteria than they release, host immune responses benefit the disease control. 
Otherwise, host immune responses benefit the pathogen development.

There is another example that immune responses lead to both unfavorable and favorable consequences. Increasing the loss of the infected macrophages can cause active disease, which is shown in Figure \ref{fig_bif_del_2d} (d).
However, increasing the macrophage- and T-cell- mediated immune responses benefit the clearance of disease, which is shown in Figure \ref{fig_bif_del_2d} (e) and (f).

The different disease outcomes can be represented as trivial, multiple, and non-trivial steady states in a deterministic ODE model. 
The disease progression between latent TB infection to active disease is a transition between co-existing multiple steady states and single non-trivial steady state. 
The transition critical point is a saddle-node bifurcation. 
Analogously, the fate of the disease after initial exposure is determined by a transcritical bifurcation. 
Therefore, transcritical and saddle-node bifurcations serve as separatrices in the parameter space.
Focusing on the therapy targeting parameters, including the bacterial proliferation rate $\delta$, the infected macrophage loss rate $b$, the cell-mediated immunity rate $\gamma$, and the macrophage killing rate $\eta$, 1- and 2-dimensional bifurcation analyses delimit the parameter range for different disease outcomes. 
We then use the identified parameter regions and consider the demographic variations in Mtb pathogen and host immune cell populations and the environmental variations on host immunity under therapies to study the disease dynamics.

We developed two It\^{o} stochastic models with only demographic variations and with both demographic and environmental variations were developed based on the deterministic in-host Tuberculosis model. 
In Region 2, the ODE model \eqref{eqn1} predicts that the TB infection can be eliminated if the initial condition is a low infection level or can progress to active disease if the initial condition is a high infection level.
With a low initial infection level, indicating initial exposure to Mtb pathogens, the SDE model \eqref{sde1} with demographic variations predicts that the expectation of the uninfected macrophages' population is close to the uninfected macrophage population from the deterministic model. The histogram of the approximate stationary distribution of the uninfected macrophage population is close to a normal distribution. Even though  infected macrophages and bacterial populations have non-zero but small expectations, their peaks are close to zero. 
Therefore,  after initial exposure, the SDE model \eqref{sde1} predicts that uninfected and infected macrophage and Mtb bacterial populations have a low level of variation due to demographic variations in cellular level. 
Interestingly, the stationary distribution of the T-cell population peaks at the uninfected level from the deterministic model prediction, but its mean is relatively large and has a large standard deviation. This large variation implies large stochastic variations applied on T-cell population, see Figure \ref{fig_sde1}.

If initial conditions take high-level infection values, the histograms of approximate stationary distribution for immune cell and bacterial populations are close to normal distribution shapes. Moreover, their means are close to the deterministic model prediction.
However, the Mtb bacterial and T-cell populations have large standard deviations, which implies that the stochastic demographic fluctuations have large influences on those two populations, as seen in Figure \ref{fig_sde2}. 
In Region 3, with low initial infection level, the SDE model prediction shows similar patterns as the prediction in Region 2, as seen in Figure \ref{fig_sde3}.
In Region 4, where the disease will eventually develop to active disease according to the deterministic model prediction, the expectations of the SDE model prediction agree with the deterministic model prediction. However, the standard deviation of uninfected macrophage, Mtb bacterial, and T-cell populations are large. 
This implies large stochastic fluctuations in the cellular level, as shown in Figure \ref{fig_sde4}. 
Interestingly, there exist small bars in Figure \ref{fig_sde4} (c)-(f) that are close to the uninfected equilibrium of the deterministic model prediction.
These small bars imply a small possibility of disease clearance.

In addition to demographic variations, the SDE model \eqref{sde2} considers environmental variations from therapies.
The pathogen-directed therapy mostly influence the pathogen proliferation rate $\delta$. 
The host-directed therapy mainly affects the host-immunity parameters  ($b$, $\gamma$, and $\eta$). 
The simulation results imply that a fast return rate along with host-immunity parameters taken in clearance region of Region 2, are most likely to control the disease progression.
Note that the return rate indicates  the speed of returning to the desirable therapy outcome by drug administration, see Section \ref{sec_sde2}.

The basic reproduction number in the parameter Region 4 is greater than one, but the approximated stationary distribution results in Figure \ref{fig_sde4} (c)-(f) show the possibility of disease clearance. This probability can be studied through continuous-time Markov chain model. This will be the topic for a future project.

Note that the coarse-graining model \eqref{eqn1}, which we used in this study, is based on the computational models constructed by \cite{wigginton2001model}. A feature of the model  \eqref{eqn1} is that it only includes the essential cell populations and biological mechanisms for the  Mtb infection, but demonstrates various disease outcomes. 
This implies that the model \eqref{eqn1} with a smaller number of variables and parameters can still predict the complex disease dynamics.
This feature allows for a rigorous mathematical analysis, which leads to robust qualitative results over large parameter ranges.
The simple structure of this model also  reveals the underlying causal mechanisms for model behaviors (\cite{zhang2014viral}).
The baseline parameter values in Table \ref{tab1} are inherited from the computational model by \cite{wigginton2001model}. 
The parameters, which can be varied by therapies, are provided with their ranges from multiple studies. 
The absence of parameter estimation from experimental data does not weaken the theoretical insights obtained from such models (\cite{alexander2011self}).
Our models have smaller numbers of variables and parameters with equally predicting capability compared to complex computational models. This feature allows advance mathematical analysis to reveal robust dynamical behaviors for the underlying mechanisms in theoretical biology.

  Finally, we recognize the limitations of this work.
It is unknown whether the levels of noise present in the simulation results reflect the levels of demographic and environmental noise in the lung.  Future studies may be needed to examine this and, if the levels are considerably different, then alterations to the model may be necessary. Moreover, in future work, we will consider altering the ODE model, such as adjusting the rates of the infected-induced T cell proliferation to saturate with the concentration of infected macrophages and Mtb bacteria.

\section{Acknowledgment}
The author would like to express her appreciation  to Dr. Linda Allen from Texas Tech University for her comments and suggestions on the stochastic formulation and simulation.  
Special thanks are also extended to Dr. Leif Ellingson from Texas Tech University for his help with carefully editing the paper.
The author also thanks both referees for their
comments and suggestions, which are very helpful for improving the manuscript.
The author acknowledges the generous support from Simons Foundation Collaboration Grants for Mathematicians, award No: A21-0013-001.

\printbibliography
\end{document}